\documentclass[12pt,oneside,reqno]{amsart}
\usepackage{graphicx}
\usepackage{mathrsfs}
\usepackage{stmaryrd}
\usepackage{amsfonts}
\usepackage{enumerate,amsmath,amssymb,amsthm}
\usepackage{booktabs}
\usepackage{diagbox}
\usepackage{amsfonts,color}

\newcommand{\dif}{\mathrm{d}}

\newcommand{\be}{\begin{eqnarray}}
\newcommand{\ee}{\end{eqnarray}}
\newcommand{\ce}{\begin{eqnarray*}}
\newcommand{\de}{\end{eqnarray*}}
\newtheorem{theorem}{Theorem}[section]
\newtheorem{lemma}[theorem]{Lemma}
\newtheorem{remark}[theorem]{Remark}
\newtheorem{definition}[theorem]{Definition}
\newtheorem{proposition}[theorem]{Proposition}
\newtheorem{Examples}[theorem]{Examples}
\newtheorem{corollary}[theorem]{Corollary}
\def\e{\varepsilon}

\def\a{\alpha}
\def\o{\omega}
\def\b{\beta}
\def\d{\delta}
\def\p{\partial}

\def\s{\sigma}

\def\l{\lambda}
\def\la{\langle}
\def\ra{\rangle}
\def\[{{\Big[}}
\def\]{{\Big]}}
\def\<{{\langle}}
\def\>{{\rangle}}
\def\({{\Big(}}
\def\){{\Big)}}

\def\no{\nonumber}
\def\bt{\begin{theorem}}
\def\et{\end{theorem}}
\def\bl{\begin{lemma}}
\def\el{\end{lemma}}
\def\br{\begin{remark}}
\def\er{\end{remark}}
\def\bx{\begin{Examples}}
\def\ex{\end{Examples}}
\def\bd{\begin{definition}}
\def\ed{\end{definition}}
\def\bp{\begin{proposition}}
\def\ep{\end{proposition}}
\def\bc{\begin{corollary}}
\def\ec{\end{corollary}}

\def\cL{{\mathcal L}}

\def\cP{{\mathcal P}}

\def\mE{{\mathbb E}}

\def\mP{{\mathbb P}}

\def\mR{{\mathbb R}}

\def\mW{{\mathbb W}}

\def\sB{{\mathscr B}}
\def\sC{{\mathscr C}}

\def\sF{{\mathscr F}}

\def\sL{{\mathscr L}}

\def\tx{\tilde{x}}
\def\ty{\tilde{y}}
\def\geq{\geqslant}
\def\leq{\leqslant}

\begin{document}

\allowdisplaybreaks
\title{Averaging principles and central limit theorems for multiscale McKean-Vlasov stochastic systems}

\author{Jie Xiang and Huijie Qiao}

\thanks{{\it AMS Subject Classification(2020):} 60H10}

\thanks{{\it Keywords:} Multiscale McKean-Vlasov stochastic systems, averaging principles, the Poisson equation, central limit theorems}

\thanks{This work was supported by NSF of China (No.12071071) and the Jiangsu Provincial Scientific Research Center of Applied Mathematics (No. BK20233002).}

\thanks{Corresponding author: Huijie Qiao, hjqiaogean@seu.edu.cn}

\subjclass{}

\date{}

\dedicatory{School of Mathematics,
Southeast University,\\
Nanjing, Jiangsu 211189, P.R.China\\
}

\begin{abstract}
In this paper, we study a class of multiscale McKean-Vlasov stochastic systems where the entire system depends on the distribution of the fast component. First of all, by the Poisson equation method we prove that the slow component converges to the solution of the averaging equation in the $L^p$ ($p\geq 2$) space with the optimal convergence rate $\frac12$. Then a central limit theorem is established by tightness.
\end{abstract}

\maketitle \rm

\section{Introduction}\label{intro}

McKean-Vlasov stochastic differential equations (SDEs for short), also known as \\distribution-dependent or mean-field SDEs, are characterized by the dependence of their coefficients on the probability law of the solutions. The study of McKean-Vlasov SDEs dates back to the foundational work of McKean \cite{hpm}, which was inspired by Kac's program in Kinetic Theory. Over the years, extensive research on McKean-Vlasov SDEs has yielded significant results. For example, Ding and Qiao \cite{dq1, dq2} explored the well-posedness and stability of solutions for McKean-Vlasov SDEs with non-Lipschitz coefficients. Wang \cite{Wang} established the exponential ergodicity of strong solutions for Landau-type McKean-Vlasov SDEs. For a broader and more detailed overview, references such as \cite{brw, hw, hrw, rz} and the works cited therein provide valuable insights.

The averaging principle, first introduced by Khasminskii in the seminal work \cite{rk}, has proven to be a highly effective and essential tool in the analysis of SDEs for modeling problems across various practical research fields. The core concept of the stochastic averaging principle is to simplify the study of complex stochastic equations by examining corresponding averaging stochastic equations, thereby offering a convenient and efficient approach for investigating many significant properties. In recent years, the stochastic averaging principle has been extended to encompass a broader range of SDEs, as demonstrated in works such as \cite{dsxz, pix, Shen1, Shen2}.

In recent years, the averaging principle for multiscale McKean-Vlasov SDEs has attracted increasing attention (\cite{hll, lx, qw, rsx, xllm}). Let us mention some results related with our work. When the coefficients are independent of the distribution of the fast component, R\"{o}ckner, Sun and Xie \cite{rsx} achieved the optimal convergence order of $\frac{1}{2}$ in the $L^2$ space by leveraging the Poisson equation. Later, Qiao and Wei \cite{qw} considered a type of multiscale McKean-Vlasov systems where the entire system depend on the distribution of the fast component and established the strong convergence rate $\frac{1}{4}$ in the $L^2$ space through a modified Khasminskii time discretization approach. Note that in \cite{qw} the strong convergence rate $\frac{1}{4}$ is not optimal. Thus, in this paper our first aim is to improve this strong convergence rate $\frac{1}{4}$ to the rate $\frac 1 2$. Concretely speaking, consider the following slow-fast system on $\mR^n\times\mR^m$:
\be\left\{\begin{array}{l}
\dif X_t^{\e}=b_1(X_t^{\e}, \sL_{X_t^{\e}}, Y_t^{\e, y_0, \sL_{\xi}}, \sL_{Y_t^{\e, \xi}}) \dif t+\s_1(X_t^{\e}, \sL_{X_t^{\e}}) \dif B_t, \\
X_0^{\e}=\varrho, \quad 0\leq t\leq T, \\
\dif Y_t^{\e, \xi}=\frac{1}{\e} b_2(Y_t^{\e, \xi}, \sL_{Y_t^{\e, \xi}}) \dif t+\frac{1}{\sqrt{\e}} \s_2(Y_t^{\e, \xi}, \sL_{Y_t^{\e, \xi}}) \dif W_t, \\
Y_0^{\e, \xi}=\xi, \quad 0\leq t\leq T, \\
\dif Y_t^{\e, y_0, \sL_{\xi}}=\frac{1}{\e} b_2(Y_t^{\e, y_0, \sL_{\xi}}, \sL_{Y_t^{\e, \xi}}) \dif t+\frac{1}{\sqrt{\e}} \s_2(Y_t^{\e, y_0, \sL_{\xi}}, \sL_{Y_t^{\e, \xi}}) \dif W_t, \\
Y_0^{\e, y_0, \sL_{\xi}}=y_0, \quad 0\leq t\leq T,
\end{array}\right.
\label{orieq}
\ee
where $\left(B_t\right),\left(W_t\right)$ are $d_1$-dimensional and $d_2$-dimensional standard Brownian motions, respectively, defined on the complete filtered probability space $(\Omega,\sF,\{\sF_t\}_{t \in[0, T]}, \mP)$ and are mutually independent. Moreover, these mappings
$b_1: \mR^n \times \cP_2\left(\mR^n\right) \times \mR^m \times \cP_2\left(\mR^m\right) \rightarrow \mR^n$,
$\s_1: \mR^n \times \cP_2\left(\mR^n\right) \rightarrow \mR^{n \times d_1}$,
$b_2: \mR^m \times \cP_2\left(\mR^m\right) \rightarrow \mR^m$,
$\s_2: \mR^m \times \cP_2\left(\mR^m\right) \rightarrow \mR^{m \times d_2}$
are all Borel measurable, and $\varrho, \xi$ are two random variables. We employ the method of the Poisson equation to demonstrate that for some $p\geq 2$,
\ce
\mE\left(\sup_{0 \leq t \leq T}\vert X_t^{\e}- \bar{X}_t\vert^p\right)\leq C_T\e^{\frac p2},
\de
where $X^\e$ is the solution for the slow part of the system (\ref{orieq}) and $\bar{X}$ is the solution of the averaging equation (See Eq.(\ref{barxeq})). Note that in \cite{lx}, for a general multiscale McKean-Vlasov stochastic system depending on the distribution of the fast component, Li and Xie obtained that (See \cite[Theorem 2.1]{lx})
\ce
\sup_{0 \leq t \leq T}\mE|X_t^{\e}- \bar{X}_t|^p\leq C_T\e^{\frac p2}.
\de
It is obvious that our result can not be covered by \cite[Theorem 2.1]{lx}. Moreover, our result is stronger in a certain sense.

Next, it is natural to investigate whether $\frac{X^{\e}-\bar{X}}{\sqrt\e}$ converges and if it does, what is the limit. This type of problem is usually solved by the central limit theorem. The study of the central limit theorem in the context of McKean-Vlasov SDEs has garnered significant attention in recent years (\cite{fyy, flqz, sy}). For a multiscale McKean-Vlasov stochastic system independent of the distribution of the fast component, Hong et al. \cite{hlls} proved a central limit type theorem by the Poisson equation technique. In this paper, we follow the same line as that in \cite{hlls} and establish a central limit theorem for the system (\ref{orieq}). However, the extension is nontrivial. Since the fast component of the system (\ref{orieq}) involves both states and distributions, the associated Poisson equation incorporates Lions derivatives with respect to distributions, making it more complex than those in \cite{hlls} and \cite{rsx}. The first challenge we face lies in establishing the regularity of this Poisson equation under weak assumptions. By leveraging the Markov property of the fast component and employing refined estimates, we are able to overcome this difficulty. The second challenge pertains to proving the weak convergence of the auxiliary process $\vartheta_\cdot^{\varepsilon}$ (see Proposition \ref{weak}). Due to the dependence of $b_1$ on the distributions of both the slow and fast components, demonstrating the weak convergence of certain terms involving $b_1$ becomes highly nontrivial. Fortunately, this issue is resolved through a time discretization technique. Besides, we mention that in \cite{lx}  Li and Xie also observed the central limit theorem for a slow-fast McKean-Vlasov stochastic system, where the entire system depends on the distributions of both the slow motion and the fast component. Comparing our result with \cite[Theorem 2.3]{lx}, we find that our conditions are more general (Also see Section \ref{example}).

The novelty of this paper is reflected in two aspects. Firstly, we achieve the $L^p$ convergence with the optimal strong convergence order of $\frac 1 2$. Secondly,  we establish a central limit theorem for multiscale McKean-Vlasov SDEs where the entire system depends on the distribution of the fast component. It is worth stressing that due to the application of Poisson equation techniques, it is necessary to estimate the moments of both the fast and slow processes to a high degree. 

Finally, we say our motivation about the system (\ref{orieq}). Let us consider the following usual multiscale McKean-Vlasov SDEs where the fast component depends on its distribution:
\be\left\{\begin{array}{l}
\dif X_t^{\e}=b_1(X_t^{\e}, \sL_{X_t^{\e}}, Y_t^{\e}, \sL_{Y_t^{\e}}) \dif t+\s_1(X_t^{\e}, \sL_{X_t^{\e}}) \dif B_t, \\
X_0^{\e}=\varrho, \quad 0\leq t\leq T, \\
\dif Y_t^{\e}=\frac{1}{\e} b_2(Y_t^{\e}, \sL_{Y_t^{\e}}) \dif t+\frac{1}{\sqrt{\e}} \s_2(Y_t^{\e}, \sL_{Y_t^{\e}}) \dif W_t, \\
Y_0^{\e}=\xi, \quad 0\leq t\leq T.
\end{array}
\right.
\label{usuammvsde}
\ee
If we study the strong averaging principle and the central limit theorem of the system (\ref{usuammvsde}) by a time discretization technique, the Markov property of $Y^\e$ is needed (\cite{hlls, rsx}). But we know that $Y^\e$ is {\it not} Markov (\cite{Wang}). It is lucky that in the system (\ref{orieq}) $(Y_t^{\e, y_0, \sL_{\xi}}, \sL_{Y_t^{\e, \xi}})$ is Markov (\cite{rrw}). Therefore, in order to obtain the central limit theorem of the system (\ref{usuammvsde}) under weak conditions, we investigate the central limit theorem of the system (\ref{orieq}) by a time discretization technique. Moreover, the fast component of the system (\ref{orieq}) exhibits several additional important applications (\cite{blpr, cm, rrw}). Hence, it is both meaningful and worthwhile to investigate problems related to $(Y_t^{\e, y_0, \sL_{\xi}}, \sL_{Y_t^{\e, \xi}})$. In our forthcoming work, we will study the asymptotic behavior of the system (\ref{orieq}).

The remainder of this paper is organized as follows. Section \ref{noteassu} introduces the notations and assumptions used throughout this work. Section \ref{main} states the main results. Two main theorems are proved in Sections \ref{xbarxpproo} and \ref{cltthproo}, respectively. In Section \ref{example}, we present an example that demonstrates the applicability of our results. Lastly, the Appendix includes the proofs of several key estimates for clarity and completeness.

The following convention will be used throughout the paper: $C$ with or without indices will denote different positive constants whose values may change from one place to another.

\section{Notations and assumptions}\label{noteassu}

In this section, we will recall some notations and list all assumptions.

\subsection{Notations}\label{nn}
In this subsection, we introduce some notations used in the sequel.

We use $|\cdot|$ and $\|\cdot\|$ for norms of vectors and matrices, respectively. The scalar product in $\mR^d$ is denoted by $\la\cdot,\cdot\ra$, and $A^*$ stands for the transpose of a matrix $A$.

Let $\mathcal{B}_b(\mR^n)$ be the set of all bounded Borel measurable functions on $\mR^n$, and let $C(\mR^n)$ be the collection of all  continuous functions on $\mR^n$. Furthermore, let $\sB(\mR^n)$ be the Borel $\s$-algebra on $\mR^n$ and $\cP(\mR^n)$ be the space of all probability measures defined on $\sB(\mR^n)$ carrying the usual topology of weak convergence. We write $\cP_2(\mR^n)$ for the set of all probability measures $\mu$ on $\sB(\mR^n)$ satisfying
$$
\mu(|\cdot|^2):=\int_{\mR^n}|x|^2\mu(\dif x)<\infty.
$$
It is a Polish space under the $L^2$-Wasserstein distance
$$
 \mW_2(\mu_1,\mu_2):=\inf_{\pi\in\mathscr{C}(\mu_1,\mu_2)}
 \bigg(\int_{\mR^n\times\mR^n}|x-y|^2\pi(\dif x,\dif y)\bigg)^{\frac{1}{2}},\mu_1,\mu_2 \in\cP_2(\mR^n),
 $$
where $\mathscr{C}(\mu_1,\mu_2)$ is the set of all couplings $\pi$ with marginal distributions $\mu_1$ and $\mu_2.$ Moreover, for any $x\in\mR^n$, the Dirac measure $\d_x$ belongs to $\cP_2(\mR^n)$ and if $\mu_1=\sL_X,$ $\mu_2=\sL_Y$ are the corresponding distributions of random variables $X$ and $Y$, respectively, then
$$\mW_{2}(\mu_1,\mu_2)\leq(\mE|X-Y|^{2})^{\frac12},$$
where $\mE$ stands for the expectation with respect to $\mP$.

\subsection{Derivatives for functions on $\cP_2(\mR^n)$}\label{deri}
\quad In this subsection, we recall the definition of $L$-derivative for functions on $\cP_2(\mR^n)$. The definition was first introduced by Lions. Moreover, he used some abstract probability spaces to describe the $L$-derivatives. Here, for the convenience to understand the definition, we apply a straight way to state it (cf. \cite{rw}). Let $I$ be the identity map on $\mR^n$. For $\mu\in\cP_2(\mR^n)$ and $\phi\in L^2(\mR^n,\sB(\mR^n),\mu;\mR^n)$, $\mu(\phi):=\int_{\mR^n}\phi(x)\mu(\dif x)$. Moreover, by simple calculation, it holds that $\mu\circ(I+\phi)^{-1}\in\cP_2(\mR^n)$.

\begin{definition}\label{Lderidef}
(i) A function $h:\cP_2(\mR^n)\rightarrow\mR$ is called $L$-differentiable at $\mu\in\cP_2(\mR^n)$, if the functional
$$
L^2(\mR^n,\sB(\mR^n),\mu;\mR^n)\ni\phi\rightarrow h(\mu\circ(I+\phi)^{-1})
$$
is Fr\'{e}chet differentiable at $0\in L^2(\mR^n,\sB(\mR^n),\mu;\mR^n)$; that is, there exists a unique $\xi\in L^2(\mR^n,\sB(\mR^n),\mu;\mR^n)$ such that
$$
\lim_{\mu(|\phi|^2)\rightarrow0}\frac{h(\mu\circ(I+\phi)^{-1})-h(\mu)-\mu(\la\xi,\phi\ra)}{\sqrt{\mu(|\phi|^2)}}=0.
$$
In the case, we denote $\p_\mu h(\mu)=\xi$ and call it the $L$-derivative of $h$ at $\mu$.\\
(ii) A function $h:\cP_2(\mR^n)\rightarrow\mR$ is called $L$-differentiable on $\cP_2(\mR^n)$ if $L$-derivative $\p_\mu h(\mu)$ exists for all $\mu\in\cP_2(\mR^n)$.
\end{definition}

We say that a matrix-valued function $h(\mu):=(h_{ij}(\mu))$ is $L$-differentiable at $\mu\in\cP_2(\mR^n)$, if all its
components are $L$-differentiable at $\mu$. Set
$$\p_\mu h(\mu):=(\p_\mu h_{ij}(\mu)), \quad \|\p_\mu h(\mu)\|^2_{L^2(\mu)}:=\sum_{i,j}\int_{\mR^n}|\p_\mu h_{ij}(\mu)(\tx)|^2\mu(\dif\tx).$$
Moreover, we say that $\p_\mu h(\mu)(\tx)$ is differentiable at $\tx\in \mR^n$, if all its components are differentiable at $\tx$. Set
$$\p_{\tx}\p_\mu h(\mu)(\tx):=(\p_{\tx}\p_\mu h_{ij}(\mu)(\tx)), \quad \|\p_{\tx}\p_\mu h(\mu)\|^2_{L^2(\mu)}:=\sum_{i,j}\int_{\mR^n}\|\p_{\tx}\p_\mu h_{ij}(\mu)(\tx)\|^2\mu(\dif\tx).$$

\begin{definition}\label{C11def}
For a map $h(\cdot):\cP_2(\mR^n)\rightarrow\mR$, we say $h\in C^{(1,1)}(\cP_2(\mR^n),\mR)$, if this map is $L$-differentiable at any $\mu\in\cP_2(\mR^n)$, its $L$-derivative $\p_\mu h(\mu)(\cdot):\mR^n\rightarrow\mR^n$ is differentiable at any $\tx\in\mR^n$, and the derivatives $\p_\mu h(\mu)(\tx)$ and $\p_{\tx}\p_\mu h(\mu)(\tx)$ are jointly continuous at any $(\mu,\tx)$. We say $h\in C_b^{(1,1)}(\cP_2(\mR^n),\mR)$, if $h\in C^{(1,1)}(\cP_2(\mR^n),\mR)$, and $\sup\limits _{\mu \in \cP_2(\mR^n), \tx\in \mR^n }|\p_\mu h(\mu)(\tx)|<\infty$ and
$\sup\limits_{\mu \in \cP_2(\mR^n), \tx\in \mR^n }\|\p_{\tx}\p_\mu h(\mu)(\tx)\|<\infty$.
For a matrix-valued map $h(\cdot):\cP_2(\mR^n)\rightarrow\mR^{n \times d_1}$, we say $h\in C^{(1,1)}(\cP_2(\mR^n),\mR^{n \times d_1})$ $(resp.~ C_b^{(1,1)}(\cP_2(\mR^n),\mR^{n \times d_1}))$ if all the components belong
to $C^{(1,1)}(\cP_2(\mR^n),\mR)$ $(resp.~ C_b^{(1,1)}(\cP_2(\mR^n),\mR))$.
\end{definition}

\begin{definition}\label{C2def}
For a map $h(\cdot):\mR^n\rightarrow\mR$, we say $h\in C^2(\mR^n,\mR)$, if the derivatives $\p_xh(x)$, $\p_x^2h(x)$ are continuous at any $x$. We say $h\in C_b^2(\mR^n,\mR)$, if $h\in C^2(\mR^n,\mR)$, and the derivatives $\p_xh(x)$, $\p_x^2h(x)$ are bounded at any $x$. For a matrix-valued map $h(\cdot):\mR^n\rightarrow\mR^{n \times d_1}$, we say $h\in C^2(\mR^n,\mR^{n \times d_1})$ $(resp.~C_b^2(\mR^n,\mR^{n \times d_1}))$ if all the components belong to $C^2(\mR^n,\mR)$ $(resp.~C_b^2(\mR^n,\mR))$.
\end{definition}

\begin{definition}\label{C211def}
For a matrix-valued map $h(\cdot,\cdot):\mR^n\times\cP_2(\mR^n)\rightarrow\mR^{n \times d_1}$, we say $h\in C^{2,(1,1)}(\mR^n\times\cP_2(\mR^n),\mR^{n \times d_1})$, if $h(x,\cdot)\in C^{(1,1)}(\cP_2(\mR^n),\mR^{n \times d_1})$ for any $x\in\mR^n$ and $h(\cdot,\mu)\in C^2(\mR^n,\mR^{n \times d_1})$ for any $\mu\in\cP_2(\mR^n)$, and the derivatives $\p_xh(x,\mu)$, $\p_x^2h(x,\mu)$, $\p_\mu h(x,\mu)(\tx)$, $\p_{\tx}\p_\mu h(x,\mu)(\tx)$ are jointly continuous at any $(x,\mu,\tx)$.
We say $h\in C_b^{2,(1,1)}(\mR^n\times\cP_2(\mR^n),\mR^{n \times d_1})$, if $h\in C^{2,(1,1)}(\mR^n\times\cP_2(\mR^n),\mR^{n \times d_1})$, and the derivatives $\p_xh(x,\mu)$, $\p_x^2h(x,\mu)$, $\p_\mu h(x,\mu)(\tx)$, $\p_{\tx}\p_\mu h(x,\mu)(\tx)$ are uniformly bounded w.r.t. $(x,\mu,\tx)$.
\end{definition}

\begin{definition}\label{C211211def}
For a vector-valued map $h(\cdot,\cdot,\cdot,\cdot):\mR^n\times\cP_2(\mR^n)\times\mR^m\times\cP_2(\mR^m)
\rightarrow\mR^n$, we say $h\in C^{2,(1,1),2,(1,1)}(\mR^n\times\cP_2(\mR^n)\times\mR^m\times
\cP_2(\mR^m),\mR^n)$, if $h(\cdot,\cdot,y,\nu)\in C^{2,(1,1)}(\mR^n\times\cP_2(\mR^n),\mR^n)$ for any $(y,\nu)\in\mR^m\times\cP_2(\mR^m)$, and $h(x,\mu,\cdot,\cdot)\in C^{2,(1,1)}(\mR^m\times\cP_2(\mR^m),\mR^n)$ for any $(x,\mu)\in\mR^n\times\cP_2(\mR^n)$, and the partial derivatives 
\ce
&&\p_xh(x,\mu,y,\nu),\quad \p_x^2h(x,\mu,y,\nu), \quad\p_\mu h(x,\mu,y,\nu)(\tx), \quad \p_{\tx}\p_\mu h(x,\mu,y,\nu)(\tx), \\
&&\p_yh(x,\mu,y,\nu),\quad\p_y^2h(x,\mu,y,\nu),\quad \p_x\p_yh(x,\mu,y,\nu),\\
&&\p_\mu\p_y h(x,\mu,y,\nu)(\tx),\quad \p_\nu h(x,\mu,y,\nu)(\ty),\quad \p_{\ty}\p_\nu h(x,\mu,y,\nu)(\ty)
\de
are uniformly continuous on $\mathbb{R}^n\times\mathcal{P}_2(\mathbb{R}^n)\times\mathbb{R}^m\times\mathcal{P}_2(\mathbb{R}^m)
\times\mathbb{R}^n\times\mathbb{R}^m$. We say $h\in C_b^{2,(1,1),2,(1,1)}(\mR^n\times\cP_2(\mR^n)\times\mR^m\times
\cP_2(\mR^m),\mR^n)$, if $h\in C^{2,(1,1),2,(1,1)}(\mR^n\times\cP_2(\mR^n)\times\mR^m\times
\cP_2(\mR^m),\mR^n)$, and all partial derivatives are uniformly bounded w.r.t. $(x,\mu,y,\nu,\tx,\ty)$.
\end{definition}

\subsection{Assumptions}\label{ass}

In this subsection, we present all the assumptions used in the sequel.

\begin{enumerate}[$(\mathbf{H}^1_{b_{1}, \s_{1}})$]
\item
There exists a constant $L_{b_1, \s_1}>0$ such that for $x_i \in \mR^n, \mu_i \in \cP_2\left(\mR^n\right), y_i \in \mR^m$, $\nu_i \in \cP_2\left(\mR^m\right), i=1,2$,
\ce
&&|b_1(x_1, \mu_1, y_1, \nu_1)-b_1(x_2, \mu_2, y_2, \nu_2)|^2+\|\s_1(x_1,\mu_1)-\s_1(x_2, \mu_2)\|^2 \\
&\leq& L_{b_1, \s_1}(|x_1-x_2|^2+\mW_2^2(\mu_1, \mu_2)+|y_1-y_2|^2+\mW_2^2(\nu_1, \nu_2)) .
\de
\end{enumerate}
\begin{enumerate}[$(\mathbf{H}^2_{b_{1}})$]
\item
$b_1\in C_b^{2,(1,1),2,(1,1)}(\mR^n\times\cP_2(\mR^n)\times\mR^m\times\cP_2(\mR^m),\mR^n)$. In addition, for $y_i \in \mR^m$, $\nu_i \in\cP_2(\mR^m), i=1,2$,
$$
\sup _{x \in \mR^n, \mu \in \cP_2(\mR^n)}\|\p_{x} b_1(x, \mu, y_1,\nu_1)-\p_{x} b_1(x, \mu, y_2,\nu_2)\| \leq C(|y_1-y_2|+\mW_2(\nu_1,\nu_2)),
$$

$$
\sup _{x \in \mR^n, \mu \in \cP_2(\mR^n)}\|\p_{\mu} b_1(x, \mu, y_1,\nu_1)-\p_{\mu} b_1(x, \mu, y_2,\nu_2)\|_{L^2(\mu)} \leq C(|y_1-y_2|+\mW_2(\nu_1, \nu_2)),
$$

$$
\sup _{x \in \mR^n, \mu \in \cP_2(\mR^n)}\|\p_{y} b_1(x, \mu, y_1,\nu_1)-\p_{y} b_1(x, \mu, y_2,\nu_2)\| \leq C(|y_1-y_2|+\mW_2(\nu_1, \nu_2)),
$$

$$
\sup _{x \in \mR^n, \mu \in \cP_2(\mR^n)}\|\p_{x}^2 b_1(x, \mu, y_1,\nu_1)-\p_{x}^2 b_1(x, \mu, y_2,\nu_2)\| \leq C(|y_1-y_2|+\mW_2(\nu_1, \nu_2)),
$$

$$
\sup _{x \in \mR^n, \mu \in \cP_2(\mR^n)}\|\p_{\tx}\p_{\mu} b_1(x, \mu, y_1,\nu_1)-\p_{\tx}\p_{\mu} b_1(x, \mu, y_2,\nu_2)\|_{L^2(\mu)} \leq C(|y_1-y_2|+\mW_2(\nu_1, \nu_2)),
$$
\end{enumerate}
\begin{enumerate}[$(\mathbf{H}^1_{b_{2}, \s_{2}})$]
\item
There exists a constant $L_{b_2, \s_2}>0$ such that for $y_i \in \mR^m$, $\nu_i \in \cP_2(\mR^m), i=1,2$,
$$
|b_2(y_1, \nu_1)-b_2(y_2, \nu_2)|^2+\|\s_2(y_1, \nu_1)-\s_2(y_2, \nu_2)\|^2\leq L_{b_2, \s_2}(|y_1-y_2|^2+\mW_2^2(\nu_1, \nu_2)).
$$
\end{enumerate}
\begin{enumerate}[$(\mathbf{H}^2_{b_{2}, \s_{2}})$]
\item
For some $p\geq 2$, there exist two constants $\b_1>0, \b_2>0$ satisfying $\b_1-\b_2>6pL_{b_2, \s_2}$ such that for $y_i \in \mR^m$, $\nu_i \in \cP_2(\mR^m), i=1,2$,
\ce
&&2\la y_1-y_2, b_2(y_1, \nu_1)-b_2(y_2, \nu_2)\ra+(3p-1)\|\s_2(y_1, \nu_1)-\s_2(y_2, \nu_2)\|^2\\
&\leq& -\b_1|y_1-y_2|^2+\b_2\mW_2^2(\nu_1, \nu_2).
\de
\end{enumerate}
\begin{enumerate}[$(\mathbf{H}^3_{b_{2}, \s_{2}})$]
\item
$b_2\in C_b^{2,(1,1)}(\mR^m\times\cP_2(\mR^m),\mR^m)$, $\s_2\in C_b^{2,(1,1)}(\mR^m\times\cP_2(\mR^m),\mR^{m\times d_2})$. In addition, for $y_i \in \mR^m$, $\nu_i \in\cP_2(\mR^m), i=1,2$,
\ce
\|\p_{y} b_2(y_1,\nu_1)-\p_{y} b_2(y_2,\nu_2)\|+\|\p_{y} \s_2(y_1,\nu_1)-\p_{y} \s_2(y_2,\nu_2)\| \leq C(|y_1-y_2|+\mW_2(\nu_1, \nu_2)).
\de
\end{enumerate}

\br\label{condexpl}
(i) $(\mathbf{H}^1_{b_{1}, \s_{1}})$ yields that
there exists a constant $\bar{L}_{b_1, \s_1}>0$ such that for $x \in \mR^n, \mu \in \cP_2\left(\mR^n\right), y \in \mR^m$, $\nu\in \cP_2\left(\mR^m\right)$,
\be
|b_1(x, \mu, y, \nu)|^2+\|\s_1(x, \mu)\|^2\leq \bar{L}_{b_1, \s_1}(1+|x|^2+\mu(|\cdot|^2)+|y|^2+\nu(|\cdot|^2)).\label{b1s1grow}
\ee

(ii) $(\mathbf{H}_{b_2, \s_2}^1)$ implies
there exists a constant $\bar{L}_{b_2, \s_2}>0$ such that for $y\in \mR^m$, $\nu\in \cP_2(\mR^m)$,
\be
|b_2(y, \nu)|^2+\|\s_2(y, \nu)\|^2\leq \bar{L}_{b_2, \s_2}(1+|y|^2+\nu(|\cdot|^2)).\label{b2s2grow}
\ee

(iii) By $(\mathbf{H}_{b_2, \s_2}^1)$ and $(\mathbf{H}_{b_2, \s_2}^2)$, we can obtain
for $y \in \mR^m$, $\nu\in \cP_2\left(\mR^m\right)$,
\be
2\la y, b_2(y, \nu)\ra+(3p-1)\|\s_2(y, \nu)\|^2\leq -\a_1|y|^2+\a_2\nu(|\cdot|^2)+C,\label{b2s2dis}
\ee
where $\a_1:=\b_1-3pL_{b_2, \s_2}$, $\a_2:=\b_2+(3p-1)L_{b_2, \s_2}$, $\a_1-\a_2-L_{b_2, \s_2}>0$, and $C>0$ is a constant.

(iv) By $(\mathbf{H}^2_{b_{2}, \s_{2}})$, we know that for any $y,k\in\mR^m$, $\nu\in \cP_2\left(\mR^m\right)$, and $\d>0$,
\ce
2\la b_2(y+\d k, \nu)-b_2(y, \nu), \d k\ra+(3p-1)\|\s_2(y+\d k, \nu)-\s_2(y, \nu)\|^2\leq -\b_1|\d k|^2,
\de
and
\ce
2\la \d^{-1}(b_2(y+\d k, \nu)-b_2(y, \nu)), k\ra+(3p-1)\|\d^{-1}(\s_2(y+\d k, \nu)-\s_2(y, \nu))\|^2\leq -\b_1|k|^2,
\de
which together with $(\mathbf{H}^3_{b_{2}, \s_{2}})$ yields that for any $y,k\in\mR^m$, $\nu\in \cP_2\left(\mR^m\right)$,
\be
2\la \p_y b_2(y, \nu)\cdot k, k\ra+(3p-1)\|\p_y\s_2(y, \nu)\cdot k\|^2\leq -\b_1|k|^2.
\label{parb2s2}
\ee
\er

\section{Main results}\label{main}

In this section, we provide main results of the paper.

\subsection{The averaging principle for multiscale McKean-Vlasov SDEs}\label{avemain}

In this subsection, we state the averaging principle result for multiscale McKean-Vlasov SDEs.

Let us recall the system \eqref{orieq}, i.e. for fixed $T>0$,
\ce\left\{\begin{array}{l}
\dif X_t^{\e}=b_1(X_t^{\e}, \sL_{X_t^{\e}}, Y_t^{\e, y_0, \sL_{\xi}}, \sL_{Y_t^{\e, \xi}})\dif t+\s_1(X_t^{\e}, \sL_{X_t^{\e}}) \dif B_t, \\
X_0^{\e}=\varrho, \quad 0\leq t\leq T, \\
\dif Y_t^{\e, \xi}=\frac{1}{\e} b_2(Y_t^{\e, \xi}, \sL_{Y_t^{\e, \xi}}) \dif t+\frac{1}{\sqrt{\e}} \s_2(Y_t^{\e, \xi}, \sL_{Y_t^{\e, \xi}})\dif W_t, \\
Y_0^{\e, \xi}=\xi, \quad 0\leq t\leq T, \\
\dif Y_t^{\e, y_0, \sL_{\xi}}=\frac{1}{\e} b_2(Y_t^{\e, y_0, \sL_{\xi}}, \sL_{Y_t^{\e, \xi}}) \dif t+\frac{1}{\sqrt{\e}} \s_2(Y_t^{\e, y_0, \sL_{\xi}}, \sL_{Y_t^{\e, \xi}}) \dif W_t, \\
Y_0^{\e, y_0, \sL_{\xi}}=y_0, \quad 0\leq t\leq T,
\end{array}
\right.
\de
where $\mE|\varrho|^{3p}<\infty$ and $\mE|\xi|^{3p}<\infty$ $(p~\text{is the same as that in}~(\mathbf{H}_{b_2, \s_2}^2))$. Under ($\mathbf{H}_{b_1, \s_1}^1$) and ($\mathbf{H}_{b_2, \s_2}^1$),  it follows from \cite[Theorem 2.1]{Wang} that the system \eqref{orieq} has a unique strong solution $(X_\cdot^{\e},Y_\cdot^{\e, \xi},Y_\cdot^{\e, y_0, \sL_{\xi}})$.

Consider the following system:
\be\left\{\begin{array}{l}
\dif Y_t^{\xi}=b_2(Y_t^{\xi}, \sL_{Y_t^{\xi}}) \dif t+\s_2(Y_t^{\xi}, \sL_{Y_t^{\xi}}) \dif W_t, \\
Y_0^{\xi}=\xi, \quad 0\leq t\leq T, \\
\dif Y_t^{y_0, \sL_{\xi}}= b_2(Y_t^{y_0, \sL_{\xi}}, \sL_{Y_t^{\xi}})\dif t+\s_2(Y_t^{y_0, \sL_{\xi}}, \sL_{Y_t^{\xi}}) \dif W_t, \\
Y_0^{y_0, \sL_{\xi}}=y_0, \quad 0\leq t\leq T.
\end{array}
\right.
\label{frozeq}
\ee
Based on \cite[Theorem 2.1]{Wang}, the above system admits a unique strong solution $(Y_\cdot^{\xi},Y_\cdot^{y_0, \sL_{\xi}})$ under ($\mathbf{H}_{b_2, \s_2}^1$). Observing $(Y_\cdot^{y_0,\sL_{\xi}},\sL_{Y_\cdot^{\xi}})$, \cite[Theorem 4.11]{rrw} ensures that a Markov process with the same distribution can be constructed on a new probability space.
Furthermore, for any $\varphi\in\mathcal{B}_b(\mR^m\times\cP_2\left(\mR^m\right))$, define
\ce
&&P_t(y_0,\sL_{\xi};\cdot,\cdot):=\sL_{Y_t^{y_0, \sL_{\xi}}}\times\d_{\sL_{Y_t^{\xi}}},\quad t\geq0,\\
&&P_t\varphi(y_0,\sL_{\xi}):=\int_{\mR^m\times\cP_2(\mR^m)}\varphi(y,\nu)P_t(y_0,\sL_{\xi};\dif y,\dif\nu)=\mE\varphi(Y_t^{y_0, \sL_{\xi}}, \sL_{Y_t^{\xi}}),
\de
$\{P_t:t\geq0\}$ is a transition semigroup on $\mathcal{B}_b(\mR^m\times\cP_2(\mR^m))$.
Under ($\mathbf{H}_{b_2, \s_2}^2$), from \cite[Theorem 4.12]{rrw}, it can be inferred that there exists a unique invariant probability measure $\eta\times\d_\eta$ for this Markov process, where $\eta$ is the unique invariant probability measure for the first equation of the system \eqref{frozeq} (See \cite[Theorem 3.1]{Wang}). We then construct the averaging equation on $(\Omega,\sF,\left\{\sF_t\right\}_{t \in[0, T]}, \mP)$ as follows:
\be\left\{\begin{array}{l}
\mathrm{d} \bar{X}_t=\bar{b}_1(\bar{X}_t, \sL_{\bar{X}_t}) \dif t+\s_1(\bar{X}_t, \sL_{\bar{X}_t}) \dif B_t, \\
\bar{X}_0=\varrho,
\end{array}
\right.
\label{barxeq}
\ee
where $\bar{b}_1(x,\mu)=\int_{\mR^m\times\cP_2(\mR^m)}b_1(x,\mu,y,\nu)(\eta\times\d_\eta)(\dif y, \dif\nu)$.

Now, it is the position to state the main result in this section.

\bt\label{xbarxp}
Supposed that $(\mathbf{H}_{b_1, \s_1}^1)$, $(\mathbf{H}^2_{b_{1}})$, $(\mathbf{H}_{b_2, \s_2}^1)$-$(\mathbf{H}_{b_2, \s_2}^3)$ hold. Then there exists a positive constant $C_T>0$ such that
$$
\mE\left(\sup_{0 \leq t \leq T}\vert X_t^{\e}- \bar{X}_t\vert^p\right)\leq C_T\e^{\frac p2}(1+\mE\vert\varrho\vert^{3p}+\vert
y_0\vert^{3p}+\mE\vert\xi\vert^{3p}),
$$
where $\bar{X}$ is a solution of Eq.\eqref{barxeq}.
\et

The proof of Theorem \ref{xbarxp} is placed in Section \ref{xbarxpproo}.

\subsection{The central limit theorem for multiscale McKean-Vlasov SDEs}\label{cltmain}

In this subsection, we take $p=4$ and state the central limit theorem result for multiscale McKean-Vlasov SDEs.

First of all, we introduce some notations. Define the operator $\cL$ as follows: for every $\psi\in C_b^{2,(1,1)}(\mR^m\times\cP_2(\mR^m),\mR^n)$,
\ce
\cL\psi(y,\nu)
&:=&b_2(y,\nu)\cdot\p_y\psi(y,\nu)+\frac12Tr[\s_2\s_2^*(y,\nu)\cdot\p^2_{y}\psi(y,\nu)]\\
&&+\int_{\mR^m}\left[b_2(\ty,\nu)\cdot\p_\nu\psi(y,\nu)(\ty)+\frac12Tr[\s_2\s_2^*(\ty,\nu)\cdot\p_{\ty}\p_{\nu}\psi(y,\nu)(\ty)]\right]\nu(\dif \ty),
\de
and it is easy to see that $\cL$ is the infinitesimal generator of $(Y_t^{y,\nu},\sL_{Y_t^{\xi}})$ with $\nu=\sL_{\xi}$. Then consider the following Poisson equation:
\be
-\cL\Psi(x,\mu,y,\nu)=b_1(x,\mu,y,\nu)-\bar{b}_1(x,\mu).
\label{b1poieq}
\ee
Under the assumptions $(\mathbf{H}_{b_1, \s_1}^1)$, $(\mathbf{H}^2_{b_{1}})$ and $(\mathbf{H}_{b_2, \s_2}^1)$-$(\mathbf{H}_{b_2, \s_2}^3)$,
Eq.\eqref{b1poieq} admits a solution $\Psi(x,\mu,y,\nu)$ belonging to $C^{2,(1,1),2,(1,1)}(\mR^n\times\cP_2(\mR^n)\times\mR^m\times\cP_2(\mR^m),\mR^n)$ (cf. Section \ref{poieqregu}). Finally, set
\ce
&&\p_y\Psi_{\s_2}(x,\mu,y,\nu):=\p_y\Psi(x,\mu,y,\nu)\cdot\s_2(y,\nu),\\
&&\overline{(\p_y\Psi_{\s_2})(\p_y\Psi_{\s_2})^*}(x,\mu):=\int_{\mR^m\times\cP_2(\mR^m)}
(\p_y\Psi_{\s_2}(x,\mu,y,\nu))(\p_y\Psi_{\s_2}(x,\mu,y,\nu))^*(\eta\times\d_\eta)(\dif y, \dif\nu),\\
&&\Upsilon(x,\mu):=\left(\overline{(\p_y\Psi_{\s_2})(\p_y\Psi_{\s_2})^*}(x,\mu)\right)^{\frac12}.
\de

The following central limit theorem is the main result of this section.

\bt\label{cltth}
Assume that $(\mathbf{H}_{b_1, \s_1}^1)$, $(\mathbf{H}^2_{b_{1}})$, $(\mathbf{H}_{b_2, \s_2}^1)$-$(\mathbf{H}_{b_2, \s_2}^3)$ hold and $\s_1\in C_b^{2,(1,1)}(\mR^n\times\cP_2(\mR^n),\mR^{n \times d_1})$. Then the process
$$
U^{\e}:=\frac{X^{\e}- \bar{X}}{\sqrt\e}
$$
converges weakly to $U$ in $C([0,T],\mR^n)$ as $\e\rightarrow0$, where the process $U$ is the solution to the following equation
\be\left\{\begin{array}{l}
\mathrm{d}U_t=\p_x\bar{b}_1(\bar{X}_t, \sL_{\bar{X}_t})U_t\dif t+\mE[\p_\mu\bar{b}_1(u, \sL_{\bar{X}_t})(\bar{X}_t)U_t]|_{u=\bar{X}_t}\dif t
+\p_x\s_1(\bar{X}_t, \sL_{\bar{X}_t})U_t\dif B_t \\
\qquad\quad+\mE[\p_\mu\s_1(u, \sL_{\bar{X}_t})(\bar{X}_t)U_t]|_{u=\bar{X}_t}\dif B_t+\Upsilon(\bar{X}_t, \sL_{\bar{X}_t})\mathrm{d} V_t,\\
U_0=0,
\end{array}
\right.
\label{uteq}
\ee
where $V$ is a $n$-dimensional standard Brownian motion independent of $B$.
\et

We postpone the proof of Theorem \ref{cltth} to Section \ref{cltthproo}.

\section{Proof of Theorem \ref{xbarxp}}\label{xbarxpproo}

In this section, we present the proof of Theorem \ref{xbarxp}. Our proof is divided into four parts. The first part (Subsection \ref{estiforxtyt}) and the second part (Subsection \ref{estiforymu}) provide estimates for systems \eqref{orieq} and \eqref{frozeq}, respectively. In the third part (Subsection \ref{poieqregu}), we analyze the regularity of the solution to the Poisson equation. Finally, in the fourth part (Subsection \ref{sec4.4}), we conclude the proof of Theorem 3.1 by integrating the results obtained in the preceding three subsections.

\subsection{Some estimates for the system \eqref{orieq}}\label{estiforxtyt}

\bl\label{xtyt}
Under assumptions $(\mathbf{H}_{b_1, \s_1}^1)$, $(\mathbf{H}_{b_2, \s_2}^1)$ and $(\mathbf{H}_{b_2, \s_2}^2)$, there exist positive constants $C$ and $C_T$ independent of $\e$ such that
\be
&&\mE\left(\sup_{0 \leq t \leq T}\vert X_t^{\e}\vert^{3p}\right)\leq C_T(1+\mE\vert\varrho\vert^{3p}+\vert y_0\vert^{3p}+\mE\vert\xi\vert^{3p}),\label{xeboun}\\
&&\sup_{0 \leq t \leq T}\mE\vert Y_t^{\e, \xi}\vert^{3p}\leq\mE\vert\xi\vert^{3p}+C,\label{yexibounw}\\
&&\sup_{0 \leq t \leq T}\mE\vert Y_t^{\e, y_0, \sL_{\xi}}\vert^{3p}\leq \vert y_0\vert^{3p}+C(1+\mE\vert\xi\vert^{3p}),\label{yey0xibounw}\\
&&\mE\left(\sup_{0 \leq t \leq T}|Y_t^{\e, \xi}|^{3p}\right)\leq2\mE\vert\xi\vert^{3p}+\frac{C_T(1+\mE|\xi|^{3p})}{\e},\label{yexibounn}\\
&&\mE\left(\sup_{0 \leq t \leq T}|Y_t^{\e, y_0, \sL_{\xi}}|^{3p}\right)\leq2\vert y_0\vert^{3p}
+\frac{C_T(1+\mE|\xi|^{3p}+|y_0|^{3p})}{\e}.\label{yey0xibounn}
\ee
\el
\begin{proof}
By the same deduction in \cite[Lemma 4.1]{qw}, we can derive the estimates \eqref{xeboun}-\eqref{yey0xibounw}. In the following, we will establish the estimates \eqref{yexibounn} and \eqref{yey0xibounn}.

For $Y_t^{\e, \xi}$, applying the It\^{o} formula to $|Y_t^{\e, \xi}|^{3p}$, we deduce that
\be
|Y_t^{\e, \xi}|^{3p}
&=&\vert\xi\vert^{3p}+\frac{3p}{\e}\int_0^t|Y_s^{\e, \xi}|^{3p-2}\la Y_s^{\e, \xi},b_2(Y_s^{\e, \xi},\sL_{Y_s^{\e, \xi}})\ra\dif s\no\\
&&+\frac{3p}{\sqrt{\e}}\int_0^t|Y_s^{\e, \xi}|^{3p-2}\la Y_s^{\e, \xi},\s_2(Y_s^{\e, \xi},\sL_{Y_s^{\e, \xi}})\dif W_s\ra\no\\
&&+\frac{3p(3p-2)}{2\e}\int_0^t|Y_s^{\e, \xi}|^{3p-4}\|\s_2(Y_s^{\e, \xi},\sL_{Y_s^{\e, \xi}})Y_s^{\e, \xi}\|^2\dif s\no\\
&&+\frac{3p}{2\e}\int_0^t|Y_s^{\e, \xi}|^{3p-2}\|\s_2(Y_s^{\e, \xi},\sL_{Y_s^{\e, \xi}})\|^2\dif s.
\label{yexiito}
\ee
Note that \eqref{b2s2dis}, \eqref{yexiito} and \eqref{yexibounw} imply that
\ce
|Y_t^{\e, \xi}|^{3p}
&\leq&\vert\xi\vert^{3p}+\frac{3p}{2\e}\int_0^t|Y_s^{\e, \xi}|^{3p-2}\Big[2\la Y_s^{\e, \xi},b_2(Y_s^{\e, \xi},\sL_{Y_s^{\e, \xi}})\ra\\
&&\qquad\qquad\qquad\qquad\qquad\qquad+(3p-1)\|\s_2(Y_s^{\e, \xi},\sL_{Y_s^{\e, \xi}})\|^2\Big]\dif s\\
&&+\frac{3p}{\sqrt{\e}}\int_0^t|Y_s^{\e, \xi}|^{3p-2}\la Y_s^{\e, \xi},\s_2(Y_s^{\e, \xi},\sL_{Y_s^{\e, \xi}})\dif W_s\ra\\
&\leq&\vert\xi\vert^{3p}+\frac{3p}{2\e}\int_0^t|Y_s^{\e, \xi}|^{3p-2}\left[-\a_1|Y_s^{\e, \xi}|^{2}+\a_2\mE|Y_s^{\e, \xi}|^{2}+C \right]\dif s\\
&&+\frac{3p}{\sqrt{\e}}\int_0^t|Y_s^{\e, \xi}|^{3p-2}\la Y_s^{\e, \xi},\s_2(Y_s^{\e, \xi},\sL_{Y_s^{\e, \xi}})\dif W_s\ra\\
&\leq&\vert\xi\vert^{3p}+\frac{3p}{2\e}\int_0^t\left[-\left(\a_1-\frac{3p-2}{3p}\a_2-L_{b_2, \s_2}\right)|Y_s^{\e, \xi}|^{3p}
+\frac2{3p}\a_2\mE|Y_s^{\e, \xi}|^{3p}+C \right]\dif s\\
&&+\frac{3p}{\sqrt{\e}}\int_0^t|Y_s^{\e, \xi}|^{3p-2}\la Y_s^{\e, \xi},\s_2(Y_s^{\e, \xi},\sL_{Y_s^{\e, \xi}})\dif W_s\ra\\
&\leq&\vert\xi\vert^{3p}+\frac{C_T(1+\mE|\xi|^{3p})}{\e}+\frac{3p}{\sqrt{\e}}\int_0^t|Y_s^{\e, \xi}|^{3p-2}
\la Y_s^{\e, \xi},\s_2(Y_s^{\e, \xi},\sL_{Y_s^{\e, \xi}})\dif W_s\ra.
\de
By the Burkholder-Davis-Gundy inequality, the Young inequality and \eqref{b2s2grow}, we get
\ce
&&\mE\left(\sup_{0 \leq t \leq T}|Y_t^{\e, \xi}|^{3p}\right)\\
&\leq&\mE\vert\xi\vert^{3p}+\frac{C_T(1+\mE|\xi|^{3p})}{\e}+\frac{C}{\sqrt{\e}}\mE\left[\int_0^T|Y_s^{\e, \xi}|^{6p-2}
\|\s_2(Y_s^{\e, \xi},\sL_{Y_s^{\e, \xi}})\|^2\dif s\right]^{\frac12}\\
&\leq&\mE\vert\xi\vert^{3p}+\frac{C_T(1+\mE|\xi|^{3p})}{\e}\\
&&+\frac{C}{\sqrt{\e}}\mE\left[\sup_{0 \leq t \leq T}|Y_t^{\e, \xi}|^{3p}
\int_0^T|Y_s^{\e, \xi}|^{3p-2}(1+|Y_s^{\e, \xi}|^{2}+\mE|Y_s^{\e, \xi}|^{2})\dif s\right]^{\frac12}\\
&\leq&\mE\vert\xi\vert^{3p}+\frac{C_T(1+\mE|\xi|^{3p})}{\e}+\frac12\mE\left(\sup_{0 \leq t \leq T}|Y_t^{\e, \xi}|^{3p}\right)
+\frac{C}{\e}\int_0^T(1+\mE|Y_s^{\e, \xi}|^{3p})\dif s\\
&\leq&\mE\vert\xi\vert^{3p}+\frac{C_T(1+\mE|\xi|^{3p})}{\e}+\frac12\mE\left(\sup_{0 \leq t \leq T}|Y_t^{\e, \xi}|^{3p}\right),
\de
which implies that
\ce
\mE\left(\sup_{0 \leq t \leq T}|Y_t^{\e, \xi}|^{3p}\right)\leq2\mE\vert\xi\vert^{3p}+\frac{C_T(1+\mE|\xi|^{3p})}{\e}.
\de

Finally, for $Y_t^{\e, y_0, \sL_{\xi}}$, using a similar argument as for $Y_t^{\e, \xi}$, we can obtain the estimate \eqref{yey0xibounn}. The proof is complete.
\end{proof}

\subsection{Some estimates for the system \eqref{frozeq}}\label{estiforymu}

\bl\label{ymul2}
Under assumptions $(\mathbf{H}_{b_2, \s_2}^1)$ and $(\mathbf{H}_{b_2, \s_2}^2)$, there exists a positive constant $C$ independent of $\e$ such that
$$
\mE\vert Y_t^{y_0, \sL_{\xi}}\vert^2\leq|y_0|^2e^{-\a_1 t}+e^{-(\a_1-\a_2)t}\mE|\xi|^2+C,
$$
and for $p\geq2$,
$$
\mE\vert Y_t^{y_0, \sL_{\xi}}\vert^{2p}\leq|y_0|^{2p}e^{-p(\a_1-\a_2-L_{b_2, \s_2}) t}+C\mE|\xi|^{2p}+C.
$$
\el
\begin{proof}
First of all, we estimate $Y_t^{\xi}$.
Applying the It\^{o} formula to $|Y_t^{\xi}|^{2p}$ and taking expectation, we obtain
\ce
\mE|Y_t^{\xi}|^{2p}
&=&\mE|\xi|^{2p}+{2p}\mE\int_0^t|Y_s^{\xi}|^{2p-2}\la Y_s^{\xi},b_2(Y_s^{\xi},\sL_{Y_s^{\xi}})\ra\dif s\\
&&+2p(p-1)\mE\int_0^t|Y_s^{\xi}|^{2p-4}\|\s_2(Y_s^{\xi},\sL_{Y_s^{\xi}})Y_s^{\xi}\|^2\dif s\\
&&+p\mE\int_0^t|Y_s^{\xi}|^{2p-2}\|\s_2(Y_s^{\xi},\sL_{Y_s^{\xi}})\|^2\dif s.\\
\de
By \eqref{b2s2dis}, we have
\ce
\frac{\dif}{\dif t}\mE|Y_t^{\xi}|^{2p}
&=&2p\mE\Big[|Y_t^{\xi}|^{2p-2}\la Y_t^{\xi},b_2(Y_t^{\xi},\sL_{Y_t^{\xi}})\ra\Big]
+2p(p-1)\mE\Big[|Y_t^{\xi}|^{2p-4}\|\s_2(Y_t^{\xi},\sL_{Y_t^{\xi}})Y_t^{\xi}\|^2\Big]\\
&&+p\mE\Big[|Y_t^{\xi}|^{2p-2}\|\s_2(Y_t^{\xi},\sL_{Y_t^{\xi}})\|^2\Big]\\
&\leq&p\mE\Big[|Y_t^{\xi}|^{2p-2}\Big(2\la Y_t^{\xi},b_2(Y_t^{\xi},\sL_{Y_t^{\xi}})\ra+(2p-1)\|\s_2(Y_t^{\xi},\sL_{Y_t^{\xi}})\|^2\Big)\Big]\\
&\leq&p\mE\left[|Y_t^{\xi}|^{2p-2}(-\a_1|Y_t^{\xi}|^2+\a_2\sL_{Y_t^{\xi}}(|\cdot|^2)+C)\right]\\
&\leq&-p(\a_1-\a_2-L_{b_2, \s_2})\mE|Y_t^{\xi}|^{2p}+C.
\de
The comparison theorem implies
\be
\mE|Y_t^{\xi}|^{2p}\leq\mE|\xi|^{2p}e^{-p(\a_1-\a_2-L_{b_2, \s_2})t}+C.\label{yxiboun}
\ee

Next, we deal with $Y_t^{y_0, \sL_{\xi}}$. For any $\l>0$, the It\^{o} formula yields that
\ce
&&\mE|Y_t^{y_0, \sL_{\xi}}|^{2p}e^{\l t}\\
&=&|y_0|^{2p}+\l\mE\int_0^t e^{\l s}|Y_s^{y_0, \sL_{\xi}}|^{2p}\dif s+2p\mE\int_0^te^{\l s}|Y_s^{y_0, \sL_{\xi}}|^{2p-2}
\la Y_s^{y_0, \sL_{\xi}},b_2(Y_s^{y_0, \sL_{\xi}},\sL_{Y_s^{\xi}})\ra\dif s\\
&&+2p(p-1)\mE\int_0^te^{\l s}|Y_s^{y_0, \sL_{\xi}}|^{2p-4}\|\s_2(Y_s^{y_0, \sL_{\xi}},\sL_{Y_s^{\xi}})Y_s^{y_0, \sL_{\xi}}\|^2\dif s\\
&&+p\mE\int_0^te^{\l s}|Y_s^{y_0, \sL_{\xi}}|^{2p-2}\|\s_2(Y_s^{y_0, \sL_{\xi}},\sL_{Y_s^{\xi}})\|^2\dif s\\
&\leq&|y_0|^{2p}+\l\mE\int_0^t e^{\l s}|Y_s^{y_0, \sL_{\xi}}|^{2p}\dif s\\
&&+p\mE\int_0^te^{\l s}|Y_s^{y_0, \sL_{\xi}}|^{2p-2}\left[2\la Y_s^{y_0, \sL_{\xi}},b_2(Y_s^{y_0, \sL_{\xi}},\sL_{Y_s^{\xi}})\ra
+(2p-1)\|\s_2(Y_s^{y_0, \sL_{\xi}},\sL_{Y_s^{\xi}})\|^2\right]\dif s\\
&\leq&|y_0|^{2p}+\l\mE\int_0^t e^{\l s}|Y_s^{y_0, \sL_{\xi}}|^{2p}\dif s\\
&&+p\int_0^te^{\l s}\left[-(\a_1-\a_2-L_{b_2, \s_2})\mE|Y_s^{y_0, \sL_{\xi}}|^{2p}+C\mE|Y_s^{\xi}|^{2p}+C\right]\dif s.
\de
Letting $\l=p(\a_1-\a_2-L_{b_2, \s_2})$, we have
\ce
\mE|Y_t^{y_0, \sL_{\xi}}|^{2p}e^{p(\a_1-\a_2-L_{b_2, \s_2}) t}
&\leq&|y_0|^{2p}+C\int_0^te^{p(\a_1-\a_2-L_{b_2, \s_2}) s}(\mE|Y_s^{\xi}|^{2p}+1)\dif s\\
&\leq&|y_0|^{2p}+C\int_0^te^{p(\a_1-\a_2-L_{b_2, \s_2}) s}(\mE|\xi|^{2p}+1)\dif s.
\de
Then, a simple calculation implies that
\ce
\mE|Y_t^{y_0, \sL_{\xi}}|^{2p}\leq|y_0|^{2p}e^{-p(\a_1-\a_2-L_{b_2, \s_2})t}+C\mE|\xi|^{2p}+C.
\de

Finally,  by similar arguments as above, we get
\ce
\mE|Y_t^{\xi}|^2\leq\mE|\xi|^2e^{-(\a_1-\a_2)t}+C,
\de
and
\ce
\mE\vert Y_t^{y_0, \sL_{\xi}}\vert^2\leq|y_0|^2e^{-\a_1 t}+\mE|\xi|^2e^{-(\a_1-\a_2)t}+C.
\de
The proof is complete.
\end{proof}

The following two results can be obtained by arguments similar to those in \cite[Lemma 4.5]{qw} and \cite[Lemma 4.6]{qw}.
\bl\label{y1nu1y2nu2}
Suppose that assumptions $(\mathbf{H}_{b_2, \s_2}^1)$ and $(\mathbf{H}_{b_2, \s_2}^2)$ hold. Then it holds that for any $t\geq0$, $y_i \in \mR^m$, $\zeta_i\in L^2(\Omega, \sF_0, \mP; \mR^m), i=1,2$,
\ce
&&\mE|Y_t^{\zeta_1}-Y_t^{\zeta_2}|^2\leq\mE|\zeta_1-\zeta_2|^2e^{-(\b_1-\b_2)t},\\
&&\mE|Y_t^{y_1, \sL_{\zeta_1}}-Y_t^{y_2, \sL_{\zeta_2}}|^2\leq|y_1-y_2|^2e^{-\b_1t}+\mW_2^2(\sL_{\zeta_1}, \sL_{\zeta_2})e^{-(\b_1-\b_2)t},\\
&&\mE|Y_t^{y_1, \sL_{\zeta_1}}-Y_t^{y_2, \sL_{\zeta_2}}|^4\leq|y_1-y_2|^4e^{-2(\b_1-\b_2)t}+C\mW_2^4(\sL_{\zeta_1}, \sL_{\zeta_2})te^{-2(\b_1-\b_2)t}.
\de
\el

\bl\label{meb1}
Suppose that assumptions $(\mathbf{H}_{b_1, \s_1}^1)$, $(\mathbf{H}_{b_2, \s_2}^1)$ and $(\mathbf{H}_{b_2, \s_2}^2)$ hold. Then it holds that for any $t\geq0$, $x\in \mR^n$, $\mu\in\cP_2\left(\mR^n\right)$, there exists a constant $C>0$ such that
$$
|\mE b_1(x,\mu, Y_t^{y_0, \sL_{\xi}}, \sL_{Y_t^{\xi}})-\bar{b}_1(x,\mu)|^2\leq Ce^{-(\b_1-\b_2)t}(1+|y_0|^2+\sL_{\xi}(|\cdot|^2)).
$$
\el

\subsection{Poisson equation}\label{poieqregu}

In this subsection, we analyze the regularity of the solution to the Poisson equation.

\bp\label{poieqesti}
Suppose that $(\mathbf{H}_{b_1, \s_1}^1)$, $(\mathbf{H}^2_{b_{1}})$, $(\mathbf{H}_{b_2, \s_2}^1)$-$(\mathbf{H}_{b_2, \s_2}^3)$ hold. Set
\be
\Psi(x,\mu,y,\nu):=\int_0^\infty\Big(\mE b_1(x,\mu, Y_s^{y,\nu}, \sL_{Y_s^{\xi}})-\bar{b}_1(x,\mu)\Big)\dif s, \quad \nu=\sL_{\xi}.\label{psi}
\ee
Then $\Psi(x,\mu,y,\nu)$ belongs to $C^{2,(1,1),2,(1,1)}(\mR^n\times\cP_2(\mR^n)\times\mR^m\times\cP_2(\mR^m),\mR^n)$ and is the unique solution to Eq.\eqref{b1poieq}. Moreover, it holds that for any
$x\in \mR^n$, $\mu\in\cP_2\left(\mR^n\right)$, $y\in \mR^m$, $\nu\in\cP_2\left(\mR^m\right)$,
\be
&&\sup _{x \in \mR^n, \mu \in \cP_2(\mR^n)}\max\Big\{|\Psi(x,\mu,y,\nu)|,\|\p_x\Psi(x,\mu,y,\nu)\|,\|\p_\mu\Psi(x,\mu,y,\nu)\|_{L^2(\mu)},\no\\
&&\qquad\qquad\|\p_{x}^2\Psi(x,\mu,y,\nu)\|,\|\p_{\tx}\p_\mu\Psi(x,\mu,y,\nu)(\cdot)\|_{L^2(\mu)}\Big\}\no\\
&\leq& C(1+|y|^2+\nu(|\cdot|^2))^{1/2},\label{parphiynv}
\ee
and
\be
&&\sup _{x \in \mR^n, \mu\in\cP_2(\mR^n),y\in \mR^m, \nu\in\cP_2(\mR^m)}\max\Big\{\|\p_\nu\Psi(x,\mu,y,\nu)\|_{L^2(\nu)}, \|\p_y\Psi(x,\mu,y,\nu)\|, \no\\
&&\qquad\qquad\qquad\qquad\|\p_x\p_y\Psi(x,\mu,y,\nu)\|, \|\p_\mu\p_y\Psi(x,\mu,y,\nu)\|_{L^2(\mu)},\no\\
&&\qquad\qquad\qquad\qquad\|\p_{y}^2\Psi(x,\mu,y,\nu)\|,\|\p_\nu\p_y\Psi(x,\mu,y,\nu)(\cdot)\|_{L^2(\nu)}\Big\}
\leq C.\label{parphic}
\ee
\ep

\begin{proof}
The proof is divided into three steps. First, we show that $\Psi(x,\mu,y,\nu)$ belongs to $C^{2,(1,1),2,(1,1)}(\mR^n\times\cP_2(\mR^n)\times\mR^m\times\cP_2(\mR^m),\mR^n)$ and serves as the unique solution to Eq.\eqref{b1poieq}. Next, in the second and third steps, we provide the verification for the regularity estimates \eqref{parphiynv} and \eqref{parphic}, respectively.

{\bf Step 1.} We prove that $\Psi(x,\mu,y,\nu)$ belongs to $C^{2,(1,1),2,(1,1)}(\mR^n\times\cP_2(\mR^n)\times\mR^m\times\cP_2(\mR^m),\mR^n)$ and is the unique solution to Eq.\eqref{b1poieq}.

Since $\eta$ is independent of $x$, we can take the derivative of $\bar{b}(x,\mu)$ directly with respect to $x$, yielding
\ce
\p_x^i\bar{b}_1(x,\mu)=\int_{\mR^m\times\cP_2(\mR^m)}\p_x^ib_1(x,\mu,y,\nu)(\eta\times\d_\eta)(\dif y, \dif\nu), \quad i=1,2,
\de
which implies that $\bar{b}_1(\cdot,\mu)\in C_b^{2}(\mR^n,\mR^n)$. By the same argument, it follows that $\bar{b}_1(x,\cdot)\in C_b^{(1,1)}(\cP_2(\mR^n),\mR^n)$, leading to the conclusion that $\bar{b}_1\in C_b^{2,(1,1)}(\mR^n\times\cP_2(\mR^n),\mR^n)$.

Under the assumptions $(\mathbf{H}^2_{b_{1}})$ and $(\mathbf{H}^3_{b_{2}, \s_{2}})$, and by the similar argument as in \cite[Lemma 6.2]{blpr}, we conclude that $\Psi(x,\mu,y,\nu)$ belongs to $C^{2,(1,1),2,(1,1)}(\mR^n\times\cP_2(\mR^n)\times\mR^m\times\cP_2(\mR^m),\mR^n)$.
Then by acting the generator $\cL$ on $(Y_t^{y, \nu}, \sL_{Y_t^{\xi}})$, $\nu=\sL_{\xi}$, it holds that

\ce
\cL\Psi(x,\mu,y,\nu)
&=&\cL\int_0^\infty \Big(\mE b_1(x,\mu, Y_s^{y,\nu}, \sL_{Y_s^{\xi}})-\bar{b}_1(x,\mu)\Big)\dif s\\
&=&\int_0^\infty (\cL P_s)\Big(b_1(x,\mu,y,\nu)-\bar{b}_1(x,\mu)\Big)\dif s\\
&=&\int_0^\infty \frac{\dif P_s\Big(b_1(x,\mu,y,\nu)-\bar{b}_1(x,\mu)\Big)}{\dif s}\dif s\\
&=&P_s\Big(b_1(x,\mu,y,\nu)-\bar{b}_1(x,\mu)\Big)\Big|_0^\infty\\
&=&\lim_{s\rightarrow\infty}P_s\Big(b_1(x,\mu,y,\nu)\Big)-\bar{b}_1(x,\mu)-\Big(b_1(x,\mu,y,\nu)-\bar{b}_1(x,\mu)\Big)\\
&=&-\Big(b_1(x,\mu,y,\nu)-\bar{b}_1(x,\mu)\Big),
\de
which yields that $\Psi(x,\mu,y,\nu)$ is a solution to Eq.\eqref{b1poieq}. According to the It\^{o} formula, the uniqueness can be achieved through a precise calculation.

{\bf Step 2.} We establish the estimate \eqref{parphiynv}.
\ce
|\Psi(x,\mu,y,\nu)|
&\leq&\int_0^\infty\Big|\mE b_1(x,\mu, Y_s^{y,\nu}, \sL_{Y_s^{\xi}})-\bar{b}_1(x,\mu)\Big|\dif s\\
&\leq&\int_0^\infty Ce^{\frac{-(\b_1-\b_2)s}2}(1+|y|^2+\nu(|\cdot|^2))^{1/2}\dif s\\
&\leq&C(1+|y|^2+\nu(|\cdot|^2))^{1/2}.
\de

Let $\check{b}_1(x,\mu,y,\nu,s)=\mE b_1(x,\mu, Y_s^{y,\nu}, \sL_{Y_s^{\xi}})-\bar{b}_1(x,\mu),\quad \nu=\sL_{\xi}$. For any $s_0>0$, we define
$$\tilde{b}_{1,s_0}(x,\mu,y,\nu,s):=\hat{b}_1(x,\mu,y,\nu,s)-\hat{b}_1(x,\mu,y,\nu,s+s_0),$$
where
$\hat{b}_1(x,\mu,y,\nu,s):=\mE b_1(x,\mu, Y_s^{y,\nu}, \sL_{Y_s^{\xi}}).$ Hence,
\ce
\tilde{b}_{1,s_0}(x,\mu,y,\nu,s)
&=&\mE b_1(x,\mu, Y_s^{y,\nu}, \sL_{Y_s^{\xi}})-\mE b_1(x,\mu, Y_{s+s_0}^{y,\nu}, \sL_{Y_{s+s_0}^{\xi}})\\
&=&\mE b_1(x,\mu, Y_s^{y,\nu}, \sL_{Y_s^{\xi}})-\bar{b}_1(x,\mu)\\
&&-[\mE b_1(x,\mu, Y_{s+s_0}^{y,\nu}, \sL_{Y_{s+s_0}^{\xi}})-\bar{b}_1(x,\mu)].
\de
Based on Lemma \ref{meb1}, we know that
$$
\lim_{s_0\rightarrow\infty}\tilde{b}_{1,s_0}(x,\mu,y,\nu,s)=\check{b}_1(x,\mu,y,\nu,s).
$$

In order to prove
$$\|\p_x\Psi(x,\mu,y,\nu)\|\leq C(1+|y|^2+\nu(|\cdot|^2))^{1/2},$$
it suffices to show that for any $s_0>0$, $s\geq0$, $x\in \mR^n$, $\mu\in\cP_2\left(\mR^n\right)$, $y\in \mR^m$, $\nu\in\cP_2\left(\mR^m\right)$,
\be
\|\p_x\tilde{b}_{1,s_0}(x,\mu,y,\nu,s)\|\leq Ce^{-\frac{\b_1-\b_2}2s}(1+|y|^2+\nu(|\cdot|^2))^{1/2}.\label{parxb1s0}
\ee
We prove \eqref{parxb1s0} in the Appendix.

By similar arguments as above, we get
\ce
&&\|\p_\mu\tilde{b}_{1,s_0}(x,\mu,y,\nu,s)\|_{L^2(\mu)}\leq Ce^{-\frac{\b_1-\b_2}2s}(1+|y|^2+\nu(|\cdot|^2))^{1/2},\\
&&\|\p_{x}^2\tilde{b}_{1,s_0}(x,\mu,y,\nu,s)\|\leq Ce^{-\frac{\b_1-\b_2}2s}(1+|y|^2+\nu(|\cdot|^2))^{1/2},\\
&&\|\p_{\tx}\p_\mu\tilde{b}_{1,s_0}(x,\mu,y,\nu,s)\|_{L^2(\mu)}
\leq Ce^{-\frac{\b_1-\b_2}2s}(1+|y|^2+\nu(|\cdot|^2))^{1/2},
\de
which yields that
\ce
&&\|\p_\mu\Psi(x,\mu,y,\nu)\|_{L^2(\mu)}\leq C(1+|y|^2+\nu(|\cdot|^2))^{1/2},\\
&&\|\p_{x}^2\Psi(x,\mu,y,\nu)\|\leq C(1+|y|^2+\nu(|\cdot|^2))^{1/2},\\
&&\|\p_{\tx}\p_\mu\Psi(x,\mu,y,\nu)\|_{L^2(\mu)}\leq C(1+|y|^2+\nu(|\cdot|^2))^{1/2}.
\de

{\bf Step 3.} We establish the estimate \eqref{parphic}.

Relying on \eqref{psi} and ($\mathbf{H}_{b_1, \s_1}^1$), we conclude that
\ce
&&|\Psi(x,\mu,y,\nu_1)-\Psi(x,\mu,y,\nu_2)|\no\\
&\leq& C\int_0^\infty\mE\Big(|Y_s^{y,\nu_1}-Y_s^{y,\nu_2}|+\mW_2(\sL_{Y_s^{\xi_1}}, \sL_{Y_s^{\xi_2}})\Big)\dif s,
\de
where $\nu_1=\sL_{\xi_1}$, $\nu_2=\sL_{\xi_2}$. Applying Lemma \ref{y1nu1y2nu2}, we obtain
\ce
\mW^2_2(\sL_{Y_t^{\xi_1}}, \sL_{Y_t^{\xi_2}})=\mW^2_2(\sL_{Y_t^{\xi_1^\prime}}, \sL_{Y_t^{\xi_2^\prime}})
\leq\mE|Y_t^{\xi_1^\prime}-Y_t^{\xi_2^\prime}|^2\leq\mE|\xi_1^\prime-\xi_2^\prime|^2e^{-(\b_1-\b_2)t},
\de
for all $\xi_1^\prime$, $\xi_2^\prime$ with $\sL_{\xi_1^\prime}=\sL_{\xi_1}$, $\sL_{\xi_2^\prime}=\sL_{\xi_2}$. Taking the infimum on the right-hand side over all such $\xi_1^\prime$, $\xi_2^\prime$ and invoking the characterization of the $L^2$-Wasserstein metric, we deduce that
\be
\mW^2_2(\sL_{Y_t^{\xi_1}}, \sL_{Y_t^{\xi_2}})\leq\mW_2^2(\sL_{\xi_1}, \sL_{\xi_2})e^{-(\b_1-\b_2)t}.\label{w2lyxi1lyxi2}
\ee
Therefore,
\ce
|\Psi(x,\mu,y,\nu_1)-\Psi(x,\mu,y,\nu_2)|\leq C\int_0^\infty  e^{-\frac{\b_1-\b_2}2s}\mW_2(\nu_1, \nu_2)\dif s\leq C\mW_2(\nu_1, \nu_2),
\de
which implies that
$$
\sup _{x \in \mR^n, \mu\in\cP_2(\mR^n),y\in \mR^m, \nu\in\cP_2(\mR^m)}\|\p_\nu\Psi(x,\mu,y,\nu)\|_{L^2(\nu)}\leq C.
$$

Note that $\p_y Y_{t}^{y,\nu}$  and $\p_{y}^2 Y_{t}^{y,\nu}$ satisfy that
\ce\left\{\begin{array}{l}
\dif \p_y Y_{t}^{y,\nu}=\p_y b_2(Y_t^{y, \nu},\sL_{Y_t^{\xi}})\cdot\p_y Y_{t}^{y,\nu}\dif t
+\p_y \s_2(Y_t^{y, \nu},\sL_{Y_t^{\xi}})\cdot\p_y Y_{t}^{y,\nu}\dif W_t, \\
\p_y Y_{0}^{y,\nu}=I,
\end{array}\right.
\de
and
\ce\left\{\begin{array}{l}
\dif \p_{y}^2 Y_{t}^{y,\nu}
=\big[\p_{y}^2 b_2(Y_t^{y, \nu},\sL_{Y_t^{\xi}})\cdot\p_y Y_{t}^{y,\nu}\cdot\p_y Y_{t}^{y,\nu}+\p_{y} b_2(Y_t^{y, \nu},\sL_{Y_t^{\xi}})\cdot\p_{y}^2 Y_{t}^{y,\nu}\big]\dif t\\
\qquad\qquad\quad+\big[\p_{y}^2\s_2(Y_t^{y, \nu},\sL_{Y_t^{\xi}})\cdot\p_y Y_{t}^{y,\nu}\cdot\p_y Y_{t}^{y,\nu}+\p_{y}\s_2(Y_t^{y, \nu},\sL_{Y_t^{\xi}})\cdot\p_{y}^2 Y_{t}^{y,\nu}\big]\dif W_t, \\
\p_{y}^2 Y_{0}^{y,\nu}=0.
\end{array}
\right.
\de
Applying the It\^{o} formula, we deduce that
\ce
\frac{\dif}{\dif t}\mE\|\p_y Y_{t}^{y,\nu}\|^4
&\leq&4\mE\big[\|\p_y Y_{t}^{y,\nu}\|^2\la\p_y b_2(Y_t^{y, \nu},\sL_{Y_t^{\xi}})\cdot\p_y Y_{t}^{y,\nu},\p_y Y_{t}^{y,\nu}\ra\big]\\
&&+6\mE\big[\|\p_y Y_{t}^{y,\nu}\|^2\|\p_y \s_2(Y_t^{y, \nu},\sL_{Y_t^{\xi}})\cdot\p_y Y_{t}^{y,\nu}\|^2\big]\\
&\leq&-2\b_1\mE\|\p_y Y_{t}^{y,\nu}\|^4,
\de
where \eqref{parb2s2} is used in the second inequality. Hence,
\be
\mE\|\p_y Y_{t}^{y,\nu}\|^4\leq e^{-2\b_1 t}.\label{paryYynu}
\ee
By analogous arguments, one can also derive that
\be
\mE\|\p_{y}^2 Y_{t}^{y,\nu}\|^2\leq Ce^{-\frac{\b_1}2 t}.\label{pary2Yynu}
\ee
Moreover, using the It\^{o} formula, by \eqref{parb2s2}, $(\mathbf{H}^3_{b_{2}, \s_{2}})$, \eqref{w2lyxi1lyxi2}, \eqref{paryYynu} and Lemma \ref{y1nu1y2nu2}, we obtain that
\ce
&&\frac{\dif}{\dif t}\mE\|\p_y Y_{t}^{y,\nu_1}-\p_y Y_{t}^{y,\nu_2}\|^2\\
&=&2\mE\big\la\p_y b_2(Y_t^{y, \nu_1},\sL_{Y_t^{\xi_1}})\cdot\p_y Y_{t}^{y,\nu_1}
-\p_y b_2(Y_t^{y, \nu_2},\sL_{Y_t^{\xi_2}})\cdot\p_y Y_{t}^{y,\nu_2},\p_y Y_{t}^{y,\nu_1}-\p_y Y_{t}^{y,\nu_2}\big\ra\\
&&+\mE\|\p_y \s_2(Y_t^{y, \nu_1},\sL_{Y_t^{\xi_1}})\cdot\p_y Y_{t}^{y,\nu_1}-\p_y \s_2(Y_t^{y, \nu_2},\sL_{Y_t^{\xi_2}})\cdot\p_y Y_{t}^{y,\nu_2}\|^2\\
&\leq&2\mE\big\la\p_y b_2(Y_t^{y, \nu_1},\sL_{Y_t^{\xi_1}})\cdot\big(\p_y Y_{t}^{y,\nu_1}-\p_y Y_{t}^{y,\nu_2}\big),
\p_y Y_{t}^{y,\nu_1}-\p_y Y_{t}^{y,\nu_2}\big\ra\\
&&+2\mE\big\la\big(\p_y b_2(Y_t^{y, \nu_1},\sL_{Y_t^{\xi_1}})-\p_y b_2(Y_t^{y, \nu_2},\sL_{Y_t^{\xi_2}})\big)\cdot\p_y Y_{t}^{y,\nu_2},\p_y Y_{t}^{y,\nu_1}-\p_y Y_{t}^{y,\nu_2}\big\ra\\
&&+2\mE\big\|\p_y \s_2(Y_t^{y, \nu_1},\sL_{Y_t^{\xi_1}})\cdot\big(\p_y Y_{t}^{y,\nu_1}-\p_y Y_{t}^{y,\nu_2}\big)\big\|^2\\
&&+2\mE\big\|\big(\p_y \s_2(Y_t^{y, \nu_1},\sL_{Y_t^{\xi_1}})-\p_y \s_2(Y_t^{y, \nu_2},\sL_{Y_t^{\xi_2}})\big)\cdot\p_y Y_{t}^{y,\nu_2}\big\|^2\\
&\leq&-\frac{\b_1}2\mE\|\p_y Y_{t}^{y,\nu_1}-\p_y Y_{t}^{y,\nu_2}\|^2\\
&&+C\mE\big\|\big(\p_y b_2(Y_t^{y, \nu_1},\sL_{Y_t^{\xi_1}})-\p_y b_2(Y_t^{y, \nu_2},\sL_{Y_t^{\xi_2}})\big)\cdot\p_y Y_{t}^{y,\nu_2}\big\|^2\\
&&+2\mE\big\|\big(\p_y \s_2(Y_t^{y, \nu_1},\sL_{Y_t^{\xi_1}})-\p_y \s_2(Y_t^{y, \nu_2},\sL_{Y_t^{\xi_2}})\big)\cdot\p_y Y_{t}^{y,\nu_2}\big\|^2\\
&\leq&-\frac{\b_1}2\mE\|\p_y Y_{t}^{y,\nu_1}-\p_y Y_{t}^{y,\nu_2}\|^2\\
&&+C\mE\big[\big(|Y_t^{y, \nu_1}-Y_t^{y, \nu_2}|^2+\mW_2^2(\sL_{Y_t^{\xi_1}},\sL_{Y_t^{\xi_2}}) \big)\cdot\|\p_y Y_{t}^{y,\nu_2}\|^2\big]\\
&\leq&-\frac{\b_1}2\mE\|\p_y Y_{t}^{y,\nu_1}-\p_y Y_{t}^{y,\nu_2}\|^2+C(t+1)e^{-(2\b_1-\b_2)t}\mW_2^2(\nu_1, \nu_2),
\de
where $\nu_1=\sL_{\xi_1}$, $\nu_2=\sL_{\xi_2}$.
Hence,
\be
\mE\|\p_y Y_{t}^{y,\nu_1}-\p_y Y_{t}^{y,\nu_2}\|^2\leq Ce^{-\frac{\b_1}2t}\mW_2^2(\nu_1, \nu_2).\label{paryYynu1nu2}
\ee

In view of \eqref{psi}, by ($\mathbf{H}_{b_1}^2$) and \eqref{paryYynu}, we deduce that
\ce
\|\p_y\Psi(x,\mu,y,\nu)\|
&\leq&\int_0^\infty\mE\big[\|\p_y b_1(x,\mu, Y_s^{y,\nu}, \sL_{Y_s^{\xi}})\|\cdot\|\p_y Y_s^{y,\nu}\|\big]\dif s\\
&\leq& C\int_0^\infty\mE\|\p_y Y_s^{y,\nu}\|\dif s
\leq C,
\de
and
\ce
\|\p_x\p_y\Psi(x,\mu,y,\nu)\|\leq \int_0^\infty\mE\big[\|\p_x\p_y b_1(x,\mu, Y_s^{y,\nu}, \sL_{Y_s^{\xi}})\|\cdot\|\p_y Y_s^{y,\nu}\|\big]
\dif s\leq C.
\de
Similarly, we have
$$
\|\p_\mu\p_y\Psi(x,\mu,y,\nu)\|_{L^2(\mu)}\leq C.
$$

The term $\p_{y}^2\Psi(x,\mu,y,\nu)$ can be represented as
\ce
\p_{y}^2\Psi(x,\mu,y,\nu)
&=&\int_0^\infty\mE\big[\p_{y}^2 b_1(x,\mu, Y_s^{y,\nu}, \sL_{Y_s^{\xi}})\cdot\p_y Y_s^{y,\nu}\cdot\p_y Y_s^{y,\nu}\big]\dif s\\
&&+\int_0^\infty\mE\big[\p_{y} b_1(x,\mu, Y_s^{y,\nu}, \sL_{Y_s^{\xi}})\cdot\p_{y}^2 Y_s^{y,\nu}\big]\dif s.
\de
Combining ($\mathbf{H}_{b_1}^2$), \eqref{paryYynu} and \eqref{pary2Yynu}, we have
\ce
\|\p_{y}^2\Psi(x,\mu,y,\nu)\|\leq C\int_0^\infty\mE\|\p_y Y_s^{y,\nu}\|^2\dif s+C\int_0^\infty\mE\|\p_{y}^2 Y_s^{y,\nu}\|\dif s\leq C.
\de

In the following, we estimate $\p_\nu\p_y\Psi(x,\mu,y,\nu)$.
Using ($\mathbf{H}_{b_1}^2$), \eqref{w2lyxi1lyxi2}, \eqref{paryYynu}, \eqref{paryYynu1nu2} and Lemma \ref{y1nu1y2nu2}, we get
\ce
&&\|\p_y\Psi(x,\mu,y,\nu_1)-\p_y\Psi(x,\mu,y,\nu_2)\|\no\\
&\leq&\int_0^\infty\mE\Big[\|\p_y b_1(x,\mu, Y_s^{y,\nu_1}, \sL_{Y_s^{\xi_1}})-\p_y b_1(x,\mu, Y_s^{y,\nu_2}, \sL_{Y_s^{\xi_2}})\|
\cdot\|\p_y Y_s^{y,\nu_1}\|\no\\
&&\qquad\qquad+\|\p_y b_1(x,\mu, Y_s^{y,\nu_2}, \sL_{Y_s^{\xi_2}})\|\cdot\|\p_y Y_s^{y,\nu_1}-\p_y Y_s^{y,\nu_2}\|\Big]\dif s\no\\
&\leq&C\int_0^\infty\Big[\Big(\mE|Y_s^{y,\nu_1}-Y_s^{y,\nu_2}|^2\Big)^{\frac12}+\mW_2(\sL_{Y_s^{\xi_1}},\sL_{Y_s^{\xi_2}})\Big]\cdot
\Big(\mE\|\p_y Y_s^{y,\nu_1}\|^2\Big)^{\frac12}\dif s\no\\
&&+C\int_0^\infty\Big(\mE\|\p_y Y_s^{y,\nu_1}-\p_y Y_s^{y,\nu_2}\|^2\Big)^{\frac12}\dif s\no\\
&\leq&C\mW_2(\nu_1, \nu_2),
\de
which implies that $$\sup _{x \in \mR^n, \mu\in\cP_2(\mR^n),y\in \mR^m, \nu\in\cP_2(\mR^m)}\|\p_\nu\p_y\Psi(x,\mu,y,\nu)\|_{L^2(\nu)}\leq C.$$\\
The proof is complete.
\end{proof}

\subsection{The proof of Theorem \ref{xbarxp}}\label{sec4.4}
First of all, we prove the Lipschitz continuity of the coefficient $\bar{b}_1$.

By Lemma \ref{meb1} and ($\mathbf{H}_{b_1, \s_1}^1$), for any $x_i\in \mR^n$, $\mu_i\in\cP_2\left(\mR^n\right)$, $i=1,2$, we have
\ce
&&|\bar{b}_1(x_1,\mu_1)-\bar{b}_1(x_2,\mu_2)|^2\\
&\leq&3|\bar{b}_1(x_1,\mu_1)-\mE b_1(x_1,\mu_1, Y_t^{y_0, \sL_{\xi}}, \sL_{Y_t^{\xi}})|^2\\
&&+3|\mE b_1(x_2,\mu_2, Y_t^{y_0, \sL_{\xi}}, \sL_{Y_t^{\xi}})-\bar{b}_1(x_2,\mu_2)|^2\\
&&+3|\mE b_1(x_1,\mu_1, Y_t^{y_0, \sL_{\xi}}, \sL_{Y_t^{\xi}})-\mE b_1(x_2,\mu_2, Y_t^{y_0, \sL_{\xi}}, \sL_{Y_t^{\xi}})|^2\\
&\leq&Ce^{-(\b_1-\b_2)t}(1+|y_0|^2+\sL_{\xi}(|\cdot|^2))+C(|x_1-x_2|^2+\mW_2^2(\mu_1, \mu_2)).
\de
Letting $t\rightarrow\infty$, we obtain
\be
|\bar{b}_1(x_1,\mu_1)-\bar{b}_1(x_2,\mu_2)|^2\leq C(|x_1-x_2|^2+\mW_2^2(\mu_1, \mu_2)).\label{barb1lip}
\ee
By $(\mathbf{H}^1_{b_{1}, \s_{1}})$, we conclude that Eq.\eqref{barxeq} admits a unique strong solution $\bar{X}$.

Note that
\ce
X_t^{\e}-\bar{X}_t
&=&\int_0^t\Big(b_1(X_s^{\e}, \sL_{X_s^{\e}}, Y_s^{\e, y_0, \sL_{\xi}}, \sL_{Y_s^{\e, \xi}})-\bar{b}_1(\bar{X}_s, \sL_{\bar{X}_s})\Big)\dif s\\
&&+\int_0^t\Big(\s_1(X_s^{\e}, \sL_{X_s^{\e}})-\s_1(\bar{X}_s, \sL_{\bar{X}_s})\Big) \dif B_s\\
&=&\int_0^t\Big(b_1(X_s^{\e}, \sL_{X_s^{\e}},Y_s^{\e, y_0, \sL_{\xi}}, \sL_{Y_s^{\e, \xi}})-\bar{b}_1(X_s^{\e}, \sL_{X_s^{\e}})\Big)\dif s\\
&&+\int_0^t\Big(\bar{b}_1(X_s^{\e}, \sL_{X_s^{\e}})-\bar{b}_1(\bar{X}_s, \sL_{\bar{X}_s})\Big)\dif s\\
&&+\int_0^t\Big(\s_1(X_s^{\e}, \sL_{X_s^{\e}})-\s_1(\bar{X}_s, \sL_{\bar{X}_s})\Big) \dif B_s.
\de
From the H\"{o}lder inequality, the Burkholder-Davis-Gundy inequality, ($\mathbf{H}_{b_1, \s_1}^1$) and \eqref{barb1lip}, it follows that for any $p\geq2$,
\ce
&&\mE\left(\sup_{0\leq t\leq T}|X_t^{\e}-\bar{X}_t|^p\right)\\
&\leq&C\mE\left(\sup_{0\leq t\leq T}\left|\int_0^t\Big(b_1(X_s^{\e}, \sL_{X_s^{\e}},Y_s^{\e, y_0, \sL_{\xi}}, \sL_{Y_s^{\e, \xi}})
 -\bar{b}_1(X_s^{\e}, \sL_{X_s^{\e}})\Big)\dif s\right|^p\right)\\
&&+CT^{p-1}\mE\int_0^T\left|\bar{b}_1(X_s^{\e}, \sL_{X_s^{\e}})-\bar{b}_1(\bar{X}_s, \sL_{\bar{X}_s})\right|^p\dif s\\
&&+C\mE\left[\int_0^T\left\|\s_1(X_s^{\e}, \sL_{X_s^{\e}})-\s_1(\bar{X}_s, \sL_{\bar{X}_s})\right\|^2 \dif s\right]^{\frac p2}\\
&\leq&C\mE\left(\sup_{0\leq t\leq T}\left|\int_0^t\Big(b_1(X_s^{\e}, \sL_{X_s^{\e}},Y_s^{\e, y_0, \sL_{\xi}}, \sL_{Y_s^{\e, \xi}})
 -\bar{b}_1(X_s^{\e}, \sL_{X_s^{\e}})\Big)\dif s\right|^p\right)\\
&&+C_T\int_0^T\mE|X_s^{\e}-\bar{X}_s|^p\dif s.
\de
The Grownall inequality yields that
\be
&&\mE\left(\sup_{0\leq t\leq T}|X_t^{\e}-\bar{X}_t|^p\right)\no\\
&\leq&C_T\mE\left(\sup_{0\leq t\leq T}\left|\int_0^t\Big(b_1(X_s^{\e}, \sL_{X_s^{\e}},Y_s^{\e, y_0, \sL_{\xi}}, \sL_{Y_s^{\e, \xi}})
 -\bar{b}_1(X_s^{\e}, \sL_{X_s^{\e}})\Big)\dif s\right|^p\right).\label{xebarxgrow}
\ee

In what follows, we aim to estimate the right-hand side of the inequality above. By Proposition \ref{poieqesti}, there exists a $\Psi(x,\mu,y,\nu)$ belonging to $C^{2,(1,1),2,(1,1)}(\mR^n\times\cP_2(\mR^n)\times\mR^m\times\cP_2(\mR^m),\mR^n)$  such that
$$
-\cL\Psi(x,\mu,y,\nu)=b_1(x,\mu,y,\nu)-\bar{b}_1(x,\mu).
$$
Applying the It\^{o} formula to $\Psi(X_t^{\e}, \sL_{X_t^{\e}},Y_t^{\e, y_0, \sL_{\xi}}, \sL_{Y_t^{\e, \xi}})$, we have
\be
&&\qquad\Psi(X_t^{\e}, \sL_{X_t^{\e}},Y_t^{\e, y_0, \sL_{\xi}}, \sL_{Y_t^{\e, \xi}})=\Psi(\varrho,\sL_{\varrho},y_0,\sL_{\xi})\no\\
&&+\int_0^tb_1(X_s^{\e}, \sL_{X_s^{\e}},Y_s^{\e, y_0, \sL_{\xi}}, \sL_{Y_s^{\e, \xi}})\cdot\p_x\Psi(X_s^{\e}, \sL_{X_s^{\e}},
 Y_s^{\e, y_0, \sL_{\xi}}, \sL_{Y_s^{\e, \xi}})\dif s\no\\
&&+\int_0^t\p_x\Psi(X_s^{\e}, \sL_{X_s^{\e}},Y_s^{\e, y_0, \sL_{\xi}}, \sL_{Y_s^{\e, \xi}})\cdot\s_1(X_s^{\e}, \sL_{X_s^{\e}})\dif B_s\no\\
&&+\frac12\int_0^tTr\left[\s_1\s_1^*(X_s^{\e}, \sL_{X_s^{\e}})\cdot\p_{x}^2\Psi(X_s^{\e}, \sL_{X_s^{\e}},Y_s^{\e, y_0, \sL_{\xi}}, \sL_{Y_s^{\e, \xi}})\right]\dif s\no\\
&&+\mathbb{\tilde{E}}\int_0^tb_1(\tilde{X}_s^{\e}, \sL_{X_s^{\e}},\tilde{Y}_s^{\e, y_0, \sL_{\xi}}, \sL_{Y_s^{\e, \xi}})
 \cdot\p_\mu\Psi(X_s^{\e}, \sL_{X_s^{\e}},Y_s^{\e, y_0, \sL_{\xi}}, \sL_{Y_s^{\e, \xi}})(\tilde{X}_s^{\e})\dif s\no\\
&&+\frac12\mathbb{\tilde{E}}\int_0^tTr\left[\s_1\s_1^*(\tilde{X}_s^{\e}, \sL_{X_s^{\e}})\cdot\p_{\tx}\p_\mu\Psi(X_s^{\e}, \sL_{X_s^{\e}},
 Y_s^{\e, y_0, \sL_{\xi}}, \sL_{Y_s^{\e, \xi}})(\tilde{X}_s^{\e})\right]\dif s\no\\
&&+\frac1{\e}\int_0^tb_2(Y_s^{\e, y_0, \sL_{\xi}},\sL_{Y_s^{\e, \xi}})\cdot\p_y\Psi(X_s^{\e}, \sL_{X_s^{\e}},
 Y_s^{\e, y_0, \sL_{\xi}}, \sL_{Y_s^{\e, \xi}})\dif s\no\\
&&+\frac1{2\e}\int_0^tTr\left[\s_2\s_2^*(Y_s^{\e, y_0, \sL_{\xi}}, \sL_{Y_s^{\e, \xi}})\cdot\p_{y}^2\Psi(X_s^{\e}, \sL_{X_s^{\e}},
 Y_s^{\e, y_0, \sL_{\xi}}, \sL_{Y_s^{\e, \xi}})\right]\dif s\no\\
&&+\frac1{\sqrt\e}\int_0^t\p_y\Psi(X_s^{\e}, \sL_{X_s^{\e}},Y_s^{\e, y_0, \sL_{\xi}}, \sL_{Y_s^{\e, \xi}})\cdot
\s_2(Y_s^{\e, y_0, \sL_{\xi}}, \sL_{Y_s^{\e, \xi}})\dif W_s\no\\
&&+\frac1\e\mathbb{\tilde{E}}\int_0^tb_2(\tilde{Y}_s^{\e, \tilde{\xi}}, \sL_{Y_s^{\e,\xi}})\cdot\p_\nu\Psi(X_s^{\e}, \sL_{X_s^{\e}},
 Y_s^{\e, y_0, \sL_{\xi}}, \sL_{Y_s^{\e, \xi}})(\tilde{Y}_s^{\e, \tilde{\xi}})\dif s\no\\
&&+\frac1{2\e}\mathbb{\tilde{E}}\int_0^tTr\left[\s_2\s_2^*(\tilde{Y}_s^{\e, \tilde{\xi}}, \sL_{Y_s^{\e,\xi}})\cdot\p_{\ty}\p_\nu\Psi(X_s^{\e}, \sL_{X_s^{\e}},Y_s^{\e, y_0, \sL_{\xi}}, \sL_{Y_s^{\e, \xi}})(\tilde{Y}_s^{\e, \tilde{\xi}})\right]\dif s\no\\
&=&\Psi(\varrho,\sL_{\varrho},y_0,\sL_{\xi})+N_t^{\e,1}+\frac1{\sqrt\e}N_t^{\e,2}\no\\
&&+\mathbb{\tilde{E}}\int_0^tb_1(\tilde{X}_s^{\e}, \sL_{X_s^{\e}},\tilde{Y}_s^{\e, y_0, \sL_{\xi}}, \sL_{Y_s^{\e, \xi}})
\cdot\p_\mu\Psi(X_s^{\e}, \sL_{X_s^{\e}},Y_s^{\e, y_0, \sL_{\xi}}, \sL_{Y_s^{\e, \xi}})(\tilde{X}_s^{\e})\dif s\no\\
&&+\frac12\mathbb{\tilde{E}}\int_0^tTr\left[\s_1\s_1^*(\tilde{X}_s^{\e}, \sL_{X_s^{\e}})\cdot\p_{\tx}\p_\mu\Psi(X_s^{\e}, \sL_{X_s^{\e}},
 Y_s^{\e, y_0, \sL_{\xi}}, \sL_{Y_s^{\e, \xi}})(\tilde{X}_s^{\e})\right]\dif s\no\\
&&+\int_0^t\cL_1\Psi(X_s^{\e}, \sL_{X_s^{\e}},Y_s^{\e, y_0, \sL_{\xi}}, \sL_{Y_s^{\e, \xi}})\dif s+\frac1{\e}\int_0^t\cL\Psi(X_s^{\e}, \sL_{X_s^{\e}},Y_s^{\e, y_0, \sL_{\xi}}, \sL_{Y_s^{\e, \xi}})\dif s,\no\\
\label{phiito}
\ee
where the process $(\tilde{X}_s^{\e},\tilde{Y}_s^{\e, y_0, \sL_{\xi}},\tilde{Y}_s^{\e, \tilde{\xi}})$ denotes a copy of the original process $(X_s^{\e},Y_s^{\e, y_0, \sL_{\xi}},$ $Y_s^{\e, \xi})$, constructed on a probability space $(\tilde{\Omega},\tilde{\sF},\{\tilde{\sF}_t\}_{t \in[0, T]}, \tilde{\mP})$,
which is an exact copy of the original probability space $(\Omega,\sF,\{\sF_t\}_{t \in[0, T]}, \mP)$,
\ce
&&N_t^{\e,1}:=\int_0^t\p_x\Psi(X_s^{\e}, \sL_{X_s^{\e}},Y_s^{\e, y_0, \sL_{\xi}}, \sL_{Y_s^{\e, \xi}})\cdot\s_1(X_s^{\e},\sL_{X_s^{\e}})\dif B_s,\\
&&N_t^{\e,2}:=\int_0^t\p_y\Psi(X_s^{\e}, \sL_{X_s^{\e}},Y_s^{\e, y_0, \sL_{\xi}}, \sL_{Y_s^{\e, \xi}})\cdot\s_2(Y_s^{\e, y_0, \sL_{\xi}}, \sL_{Y_s^{\e, \xi}})\dif W_s,
\de
and
\ce
\cL_1\Psi(x,\mu,y,\nu)
&:=&b_1(x,\mu,y,\nu)\cdot\p_x\Psi(x,\mu,y,\nu)\\
&&+\frac12Tr[\s_1\s_1^*(x,\mu)\cdot\p^2_{x}\Psi(x,\mu,y,\nu)].
\de
Then, we conclude that
\ce
&&\mE\left(\sup_{0\leq t\leq T}\left|\int_0^t\Big(b_1(X_s^{\e}, \sL_{X_s^{\e}},Y_s^{\e, y_0, \sL_{\xi}}, \sL_{Y_s^{\e, \xi}})
 -\bar{b}_1(X_s^{\e}, \sL_{X_s^{\e}})\Big)\dif s\right|^p\right)\\
&=&\mE\left(\sup_{0\leq t\leq T}\left|\int_0^t\cL\Psi(X_s^{\e}, \sL_{X_s^{\e}},Y_s^{\e, y_0, \sL_{\xi}}, \sL_{Y_s^{\e, \xi}})
\dif s\right|^p\right)\\
&\leq&C\e^p\mE\left(\sup_{0\leq t\leq T}\Big|\Psi(X_t^{\e}, \sL_{X_t^{\e}},Y_t^{\e, y_0, \sL_{\xi}}, \sL_{Y_t^{\e, \xi}})\Big|^p\right)
+C\e^p\mE\left|\Psi(\varrho,\sL_{\varrho},y_0,\sL_{\xi})\right|^p\\
&&+C\e^p\mE\mathbb{\tilde{E}}\Bigg(\sup_{0\leq t\leq T}\Big|\int_0^tb_1(\tilde{X}_s^{\e}, \sL_{X_s^{\e}},\tilde{Y}_s^{\e, y_0, \sL_{\xi}}, \sL_{Y_s^{\e, \xi}})\\
&&\qquad\qquad\qquad\qquad\qquad\qquad\cdot\p_\mu\Psi(X_s^{\e}, \sL_{X_s^{\e}},
Y_s^{\e, y_0, \sL_{\xi}}, \sL_{Y_s^{\e, \xi}})(\tilde{X}_s^{\e})\dif s\Big|^p\Bigg)\\
&&+C\e^p\mE\mathbb{\tilde{E}}\Bigg(\sup_{0\leq t\leq T}\Big|
\int_0^tTr\Big[\s_1\s_1^*(\tilde{X}_s^{\e}, \sL_{X_s^{\e}})\\
&&\qquad\qquad\qquad\qquad\qquad\qquad\cdot\p_{\tx}\p_\mu\Psi(X_s^{\e}, \sL_{X_s^{\e}},
 Y_s^{\e, y_0, \sL_{\xi}}, \sL_{Y_s^{\e, \xi}})(\tilde{X}_s^{\e})\Big]\dif s\Big|^p\Bigg)\\
&&+C\e^p\mE\left(\sup_{0\leq t\leq T}\Big|\int_0^t\cL_1\Psi(X_s^{\e}, \sL_{X_s^{\e}},Y_s^{\e, y_0, \sL_{\xi}}, \sL_{Y_s^{\e, \xi}})\dif s\Big|^p\right)\\
&&+C\e^p\mE\left(\sup_{0\leq t\leq T}|N_t^{\e,1}|^p\right)+C\e^{\frac p2}\mE\left(\sup_{0\leq t\leq T}|N_t^{\e,2}|^p\right)\\
&=&\sum_{i=1}^7\mathscr{J}_i.
\de

For $\mathscr{J}_1$ and $\mathscr{J}_2$, by Proposition \ref{poieqesti} and Lemma \ref{xtyt}, we have
\ce
\mathscr{J}_1
&\leq& C\e^p\left\{1+\mE\left(\sup_{0\leq t\leq T}|Y_t^{\e, y_0, \sL_{\xi}}|^p\right)+\mE\left(\sup_{0\leq t\leq T}|Y_t^{\e, \xi}|^p\right)\right\}\\
&\leq& C_T\e^{p-1}\{1+|y_0|^{3p}+\mE|\xi|^{3p}\},
\de
and
\ce
\mathscr{J}_2\leq C\e^p\{1+|y_0|^p+\mE|\xi|^p\}.
\de

For $\mathscr{J}_3$ and $\mathscr{J}_4$, based on \eqref{b1s1grow}, the H\"{o}lder inequality and Proposition \ref{poieqesti}, one can obtain that
\ce
\mathscr{J}_3&\leq&C\e^pT^{p-1}\mE\mathbb{\tilde{E}}\Big[\int_0^T|b_1(\tilde{X}_s^{\e}, \sL_{X_s^{\e}},\tilde{Y}_s^{\e, y_0, \sL_{\xi}}, \sL_{Y_s^{\e, \xi}})|^p\\
&&\qquad\qquad\qquad\qquad\quad\cdot\|\p_\mu\Psi(X_s^{\e}, \sL_{X_s^{\e}},Y_s^{\e, y_0, \sL_{\xi}}, \sL_{Y_s^{\e, \xi}})\|_{L^2(\mu)}^p \dif s\Big]\\
&\leq&C\e^pT^p\left\{1+\mE\left(\sup_{0\leq t\leq T}|X_t^{\e}|^{2p}\right)+\sup_{0\leq t\leq T}\mE|Y_t^{\e, y_0, \sL_{\xi}}|^{2p}+\sup_{0\leq t\leq T}\mE|Y_t^{\e, \xi}|^{2p}\right\},
\de
and
\ce
\mathscr{J}_4
&\leq&C\e^pT^{p-1}\mE\mathbb{\tilde{E}}\int_0^T\|\s_1\s_1^*(\tilde{X}_s^{\e}, \sL_{X_s^{\e}})\|^p\\
&&\qquad\qquad\qquad\qquad\cdot\|\p_{\tx}\p_\mu\Psi(X_s^{\e}, \sL_{X_s^{\e}},Y_s^{\e, y_0, \sL_{\xi}}, \sL_{Y_s^{\e, \xi}})\|_{L^2(\mu)}^p  \dif s\\
&\leq&C\e^pT^p\left\{1+\mE\left(\sup_{0\leq t\leq T}|X_t^{\e}|^{3p}\right)+\sup_{0\leq t\leq T}\mE|Y_t^{\e, y_0, \sL_{\xi}}|^{3p}
+\sup_{0\leq t\leq T}\mE|Y_t^{\e, \xi}|^{3p}\right\}.
\de

By similar arguments for $\mathscr{J}_3$ and $\mathscr{J}_4$, we obtain
\ce
\mathscr{J}_5\leq C\e^pT^p\left\{1+\mE\left(\sup_{0\leq t\leq T}|X_t^{\e}|^{3p}\right)+\sup_{0\leq t\leq T}\mE|Y_t^{\e, y_0, \sL_{\xi}}|^{3p}
+\sup_{0\leq t\leq T}\mE|Y_t^{\e, \xi}|^{3p}\right\}.
\de

For $\mathscr{J}_6$ and $\mathscr{J}_7$, by the H\"{o}lder inequality, the Burkholder-Davis-Gundy inequality, Remark \ref{condexpl} and Proposition \ref{poieqesti}, we have
\be
&&\mathscr{J}_6+\mathscr{J}_7\no\\
&\leq&C\e^p\mE\left[\int_0^T\|\p_x\Psi(X_s^{\e}, \sL_{X_s^{\e}},Y_s^{\e, y_0, \sL_{\xi}}, \sL_{Y_s^{\e, \xi}})\|^2
\|\s_1(X_s^{\e}, \sL_{X_s^{\e}})\|^2\dif s\right]^{\frac p2}\no\\
&&+C\e^{\frac p2}\mE\left[\int_0^T\|\p_y\Psi(X_s^{\e}, \sL_{X_s^{\e}}, Y_s^{\e, y_0, \sL_{\xi}}, \sL_{Y_s^{\e, \xi}})\|^2
\|\s_2(Y_s^{\e, y_0, \sL_{\xi}}, \sL_{Y_s^{\e, \xi}})\|^2\dif s\right]^{\frac p2}\no\\
&\leq&C\e^pT^{\frac p2-1}\int_0^T\left(1+\mE|X_s^{\e}|^{2p}+\mE|Y_s^{\e, y_0, \sL_{\xi}}|^{2p}+\mE|Y_s^{\e, \xi}|^{2p}\right)\dif s\no\\
&&+C\e^{\frac p2}T^{\frac p2-1}\int_0^T\left(1+\mE|Y_s^{\e, y_0, \sL_{\xi}}|^p
+\mE|Y_s^{\e, \xi}|^p\right)\dif s\no\\
&\leq&C(\e^p+\e^{\frac p2})T^{\frac p2}\left\{1+\mE\left(\sup_{0\leq t\leq T}|X_t^{\e}|^{2p}\right)+\sup_{0\leq t\leq T}\mE|Y_t^{\e, y_0, \sL_{\xi}}|^{2p}+\sup_{0\leq t\leq T}\mE|Y_t^{\e, \xi}|^{2p}\right\}.\no\\
\label{sJ6sJ7}
\ee

Collecting the above deduction, by Lemma \ref{xtyt} we conclude that
\be
&&\mE\left(\sup_{0\leq t\leq T}\left|\int_0^t\Big(b_1(X_s^{\e}, \sL_{X_s^{\e}},Y_s^{\e, y_0, \sL_{\xi}}, \sL_{Y_s^{\e, \xi}})
 -\bar{b}_1(X_s^{\e}, \sL_{X_s^{\e}})\Big)\dif s\right|^p\right)\no\\
&\leq& C_T\e^{\frac p2}(1+\mE\vert\varrho\vert^{3p}+\vert y_0\vert^{3p}+\mE\vert\xi\vert^{3p}),\label{b1barb1lp}
\ee
which together with \eqref{xebarxgrow} implies that
$$
\mE\left(\sup_{0 \leq t \leq T}\vert X_t^{\e}- \bar{X}_t\vert^p\right)\leq C_T\e^{\frac p2}(1+\mE\vert\varrho\vert^{3p}+\vert y_0\vert^{3p}+\mE\vert\xi\vert^{3p}).
$$
The proof is complete.

\section{Proof of theorem \ref{cltth}}\label{cltthproo}

In this section, we provide the proof of Theorem \ref{cltth}. The proof consists of three components. In the first subsection (Subsection \ref{auxipro}), we construct the auxiliary process $\vartheta_\cdot^{\e}$ and show that the difference $U_\cdot^{\e} - \vartheta_\cdot^{\e}$ converges to zero in the $L^2$ sense. In the following subsection (Subsection \ref{tight}), the tightness of $\vartheta_\cdot^{\e}$ in the space $C([0,T], \mR^n)$ is established. Finally, in Subsection \ref{cltthproode}, we prove the weak convergence of $\vartheta_\cdot^{\e}$ to the solution $U_\cdot$ of equation \eqref{uteq} in $C([0, T], \mR^n)$, thereby confirming the validity of Theorem \ref{cltth}.

\subsection{Construction of the auxiliary process}\label{auxipro}

In this subsection, we construct an auxiliary process $\vartheta_\cdot^{\e}$ and show that the difference $U_\cdot^{\e}-\vartheta_\cdot^{\e}$ converges to zero in the $L^2$ sense.

First of all, by \eqref{orieq} and \eqref{barxeq}, the deviation process $\{U_t^{\e}:=\frac{X_t^{\e}- \bar{X}_t}{\sqrt\e}\}_{t\geq0}$ satisfies that
\be\left\{\begin{array}{l}
\dif U_t^{\e}=\frac{1}{\sqrt{\e}}[b_1(X_t^{\e},\sL_{X_t^{\e}},Y_t^{\e,
y_0,\sL_{\xi}},\sL_{Y_t^{\e,\xi}})-\bar{b}_1(X_t^{\e},\sL_{X_t^{\e}})]\dif t\\
\qquad\quad+\frac{1}{\sqrt{\e}}[\bar{b}_1(X_t^{\e}, \sL_{X_t^{\e}})-\bar{b}_1(\bar{X}_t, \sL_{\bar{X}_t})] \dif t
+\frac{1}{\sqrt{\e}}[\s_1(X_t^{\e}, \sL_{X_t^{\e}})-\s_1(\bar{X}_t, \sL_{\bar{X}_t})] \dif B_t,\\
U_0^{\e}=0.
\label{ueeq}
\end{array}
\right.
\ee
According to Theorem \ref{xbarxp}, there exists a constant $C_T$ independent of $\e$ such that
\be
\mE\left(\sup_{0 \leq t \leq T}\vert U_t^{\e}\vert^4\right)\leq
C_T(1+\mE\vert\varrho\vert^{12}+\vert y_0\vert^{12}+\mE\vert\xi\vert^{12}).\label{uel4}
\ee

For any $\e\in(0,1)$, we consider the following auxiliary process $\vartheta_t^{\e}$:
\be\left\{\begin{array}{l}
\dif \vartheta_t^{\e}=\frac{1}{\sqrt{\e}}[b_1(X_t^{\e}, \sL_{X_t^{\e}},  Y_t^{\e, y_0, \sL_{\xi}}, \sL_{Y_t^{\e, \xi}})-\bar{b}_1(X_t^{\e}, \sL_{X_t^{\e}})] \dif t \\
\qquad\quad+\p_x\bar{b}_1(\bar{X}_t, \sL_{\bar{X}_t})\vartheta_t^{\e}\dif t+\mE[\p_\mu\bar{b}_1(u, \sL_{\bar{X}_t})(\bar{X}_t)\vartheta_t^{\e}]
|_{u=\bar{X}_t}\dif t\\
\qquad\quad+\p_x\s_1(\bar{X}_t, \sL_{\bar{X}_t})\vartheta_t^{\e}\dif B_t+\mE[\p_\mu\s_1(u, \sL_{\bar{X}_t})(\bar{X}_t)\vartheta_t^{\e}]
|_{u=\bar{X}_t}\dif B_t,\\
\vartheta_0^{\e}=0.
\end{array}
\right.
\label{varthetaeq}
\ee

\bp\label{uethetae}
Assume that $(\mathbf{H}_{b_1, \s_1}^1)$, $(\mathbf{H}_{b_1}^2)$ and $(\mathbf{H}_{b_2, \s_2}^1)$-$(\mathbf{H}_{b_2, \s_2}^3)$ hold, and $\s_1\in C_b^{2,(1,1)}(\mR^n\times\cP_2(\mR^n),\mR^{n \times d_1})$. It holds that
$$
\lim_{\e\rightarrow0}\mE\left(\sup_{0 \leq t \leq T}\vert U_t^{\e}-\vartheta_t^{\e}\vert^2\right)=0.
$$
\ep
\begin{proof}
Define $\chi_t^{\e}:=U_t^{\e}-\vartheta_t^{\e}$, we obtain
\ce
\dif \chi_t^{\e}
&=&\dif U_t^{\e}-\dif\vartheta_t^{\e}\\
&=&\frac{1}{\sqrt{\e}}\Big[\bar{b}_1(X_t^{\e}, \sL_{X_t^{\e}})-\bar{b}_1(\bar{X}_t, \sL_{\bar{X}_t})
-\p_x\bar{b}_1(\bar{X}_t, \sL_{\bar{X}_t})\sqrt\e U_t^{\e}\\
&&\qquad\qquad-\mE[\p_\mu\bar{b}_1(u, \sL_{\bar{X}_t})(\bar{X}_t)\sqrt\e U_t^{\e}]|_{u=\bar{X}_t}\Big] \dif t\\
&&+\p_x\bar{b}_1(\bar{X}_t, \sL_{\bar{X}_t})\chi_t^{\e}\dif t+\mE[\p_\mu\bar{b}_1(u, \sL_{\bar{X}_t})(\bar{X}_t)\chi_t^{\e}]|_{u=\bar{X}_t}\dif t\\
&&+\frac{1}{\sqrt{\e}}\Big[\s_1(X_t^{\e}, \sL_{X_t^{\e}})-\s_1(\bar{X}_t, \sL_{\bar{X}_t})-\p_x\s_1(\bar{X}_t, \sL_{\bar{X}_t})
\sqrt\e U_t^{\e}\\
&&\qquad\qquad-\mE[\p_\mu\s_1(u, \sL_{\bar{X}_t})(\bar{X}_t)\sqrt\e U_t^{\e}]|_{u=\bar{X}_t}\Big] \dif B_t\\
&&+\p_x\s_1(\bar{X}_t, \sL_{\bar{X}_t})\chi_t^{\e}\dif B_t+\mE[\p_\mu\s_1(u, \sL_{\bar{X}_t})(\bar{X}_t)\chi_t^{\e}]
|_{u=\bar{X}_t}\dif B_t.
\de
By the H\"{o}lder inequality and the Burkholder-Davis-Gundy inequality, it holds that
\be
&&\mE\left(\sup_{0 \leq t \leq T}\vert\chi_t^{\e}\vert^2\right)\no\\
&\leq&\frac{6T}{\e}\mE\int_0^T\Big|\bar{b}_1(X_t^{\e}, \sL_{X_t^{\e}})-\bar{b}_1(\bar{X}_t, \sL_{\bar{X}_t})
-\p_x\bar{b}_1(\bar{X}_t, \sL_{\bar{X}_t})\sqrt\e U_t^{\e}\no\\
&&\qquad\qquad\qquad-\mE[\p_\mu\bar{b}_1(u, \sL_{\bar{X}_t})(\bar{X}_t)\sqrt\e U_t^{\e}]|_{u=\bar{X}_t}\Big|^2\dif t\no\\
&&+6T\mE\int_0^T\left|\p_x\bar{b}_1(\bar{X}_t, \sL_{\bar{X}_t})\chi_t^{\e} \right|^2\dif t+6T\mE\int_0^T\left|\mE[\p_\mu\bar{b}_1(u, \sL_{\bar{X}_t})(\bar{X}_t)\chi_t^{\e}]|_{u=\bar{X}_t}\right|^2\dif t\no\\
&&+\frac{C}{\e}\mE\int_0^T\Big\|\s_1(X_t^{\e}, \sL_{X_t^{\e}})-\s_1(\bar{X}_t, \sL_{\bar{X}_t})
-\p_x\s_1(\bar{X}_t, \sL_{\bar{X}_t})\sqrt\e U_t^{\e}\no\\
&&\qquad\qquad\qquad-\mE[\p_\mu\s_1(u, \sL_{\bar{X}_t})(\bar{X}_t)\sqrt\e U_t^{\e}]|_{u=\bar{X}_t}\Big\|^2\dif t\no\\
&&+C\mE\int_0^T\left\|\p_x\s_1(\bar{X}_t, \sL_{\bar{X}_t})\chi_t^{\e} \right\|^2\dif t
+C\mE\int_0^T\left\|\mE[\p_\mu\s_1(u, \sL_{\bar{X}_t})(\bar{X}_t)\chi_t^{\e}]
|_{u=\bar{X}_t}\right\|^2\dif t\no\\
&=:&\sum_{k=1}^{6}I_k^\e(T).
\label{chiel2}
\ee

Next, we focus on proving the joint continuity of $\p_x\bar{b}_1(\cdot,\cdot)$ and $\p_\mu\bar{b}_1(x,\cdot)(\cdot)$, as well as the uniform boundedness of $\p_x\bar{b}_1(x,\mu)$ and $\p_\mu\bar{b}_1(x,\mu)(\tx)$. According to Lemma \ref{meb1}, we have for any $x\in \mR^n$, $\mu\in\cP_2\left(\mR^n\right)$,
$$
\bar{b}_1(x,\mu)=\lim_{t\rightarrow\infty}\mE b_1(x,\mu, Y_t^{y_0, \sL_{\xi}}, \sL_{Y_t^{\xi}}).
$$
Then it holds that
$$
\p_x\bar{b}_1(x,\mu)=\lim_{t\rightarrow\infty}\mE\big[\p_x b_1(x,\mu, Y_t^{y_0, \sL_{\xi}}, \sL_{Y_t^{\xi}})\big],
$$
and
$$
\p_\mu\bar{b}_1(x,\mu)(\tx)=\lim_{t\rightarrow\infty}\mE\big[\p_\mu b_1(x,\mu, Y_t^{y_0, \sL_{\xi}}, \sL_{Y_t^{\xi}})(\tx)\big].
$$
For any sequence $\{(\mu_n,\tx_n)\}_{n\geq1}\subset\cP_2(\mR^n)\times\mR^n$ such that $\mu_n\rightarrow\mu$ in $\mW_2$ and $|\tx_n-\tx|\rightarrow0$ as $n\rightarrow\infty$, we know that
\ce
&&\sup_{t\geq0}\mE\|\p_\mu b_1(x,\mu_n, Y_t^{y_0, \sL_{\xi}}, \sL_{Y_t^{\xi}})(\tx_n)-\p_\mu b_1(x,\mu, Y_t^{y_0, \sL_{\xi}},\sL_{Y_t^{\xi}})(\tx)\|\\
&\leq&\sup _{x\in \mR^n, y \in \mR^m, \nu \in \cP_2(\mR^m)}\|\p_\mu b_1(x,\mu_n, y,\nu)(\tx_n)-\p_\mu b_1(x,\mu, y,\nu)(\tx)\|.
\de
It follows from the definition of $b_1\in C_b^{2,(1,1),2,(1,1)}(\mR^n\times\cP_2(\mR^n)\times\mR^m\times
\cP_2(\mR^m),\mR^n)$ that $\p_\mu b_1(x,\mu,y,\nu)(\tx)$ is uniformly continuous on $\mR^n\times\cP_2(\mR^n)\times\mR^m\times\cP_2(\mR^m)
\times\mR^n$. Thus
\ce
\lim_{n\rightarrow\infty}\|\p_\mu\bar{b}_1(x,\mu_n)(\tx_n)-\p_\mu\bar{b}_1(x,\mu)(\tx)\|=0,
\de
which implies that $\p_\mu\bar{b}_1(x,\cdot)(\cdot)$ is jointly continuous. By similar argument, for any sequence $\{(x_n,\mu_n)\}_{n\geq1}\subset\mR^n\times\cP_2(\mR^n)$ satisfying $|x_n-x|\rightarrow0$ and $\mu_n\rightarrow\mu$ in $\mW_2$  as $n\rightarrow\infty$, we obtain
\ce
\lim_{n\rightarrow\infty}\|\p_x\bar{b}_1(x_n,\mu_n)-\p_x\bar{b}_1(x,\mu)\|=0,
\de
that is, $\p_x\bar{b}_1(\cdot,\cdot)$ is jointly continuous.
In addition, the uniform boundedness of $\p_x b_1$ and $\p_\mu b_1$ implies that
\ce
\sup _{x \in \mR^n, \mu \in \cP_2(\mR^n)}\|\p_x \bar{b}_1(x,\mu)\|<\infty,
\de
and
\ce
\sup _{x \in \mR^n, \mu \in \cP_2(\mR^n),\tx\in \mR^n}\|\p_\mu \bar{b}_1(x,\mu)(\tx)\|<\infty,
\de
So, we conclude that $\p_x\bar{b}_1(x,\mu)$ and $\p_\mu\bar{b}_1(x,\mu)(\tx)$ are uniformly bounded.

In the following, we estimate $I_k^\e(T)$, $k=1,\cdots,6$ individually.
For $I_1^\e(T)$, according to the mean value theorem (See \cite[Lemma 4.1]{kn}), there exists a constant $r\in(0,1)$ such that
\be
I_1^\e(T)
&=&6T\mE\int_0^T\Big|\p_x\bar{b}_1(\bar{X}_t+r\sqrt\e U_t^{\e}, \sL_{X_t^{\e}})U_t^{\e}\no\\
&&\qquad\qquad\qquad+\mE[\p_\mu\bar{b}_1(u, \sL_{\bar{X}_t+r\sqrt\e U_t^{\e}})(\bar{X}_t+r\sqrt\e U_t^{\e})U_t^{\e}]|_{u=\bar{X}_t}\no\\
&&\qquad\qquad\qquad-\p_x\bar{b}_1(\bar{X}_t, \sL_{\bar{X}_t})U_t^{\e}
-\mE[\p_\mu\bar{b}_1(u, \sL_{\bar{X}_t})(\bar{X}_t)U_t^{\e}]|_{u=\bar{X}_t}\Big|^2\dif t\no\\
&\leq&C_T\Bigg\{\int_0^T\Big[\mE\|\p_x\bar{b}_1(\bar{X}_t+r\sqrt\e U_t^{\e}, \sL_{X_t^{\e}})-\p_x\bar{b}_1(\bar{X}_t, \sL_{\bar{X}_t})\|^4\Big]^{\frac12}\dif t\no\\
&&\quad+\mE\int_0^T\mE\Big[|\p_\mu\bar{b}_1(u, \sL_{\bar{X}_t+r\sqrt\e U_t^{\e}})(\bar{X}_t+r\sqrt\e U_t^{\e})-\p_\mu\bar{b}_1(u, \sL_{\bar{X}_t})(\bar{X}_t)|^2\Big]\Big|_{u=\bar{X}_t}\dif t\Bigg\}\no\\
&&\qquad\cdot\left[\mE\left(\sup_{0 \leq t \leq T}\vert U_t^{\e}\vert^4\right)\right]^{\frac12}.\no
\ee
Note that $\p_x\bar{b}_1(\cdot,\cdot)$ and $\p_\mu\bar{b}_1(x,\cdot)(\cdot)$ are jointly continuous. Thus, by the boundedness of $\p_x\bar{b}_1$ and $\p_\mu\bar{b}_1$, and \eqref{uel4}, it holds that
\be
\lim_{\e\rightarrow0}I_1^\e(T)=0.\label{I1e}
\ee
Using similar arguments as above, we obtain
\be
\lim_{\e\rightarrow0}I_4^\e(T)=0.\label{I4e}
\ee
From the boundedness of $\p_x\bar{b}_1$, $\p_\mu\bar{b}_1$, $\p_x\s_1$ and $\p_\mu\s_1$, it follows that
\be
I_2^\e(T)+I_3^\e(T)+I_5^\e(T)+I_{6}^\e(T)\leq C_T\int_0^T\mE\left(\sup_{0 \leq s \leq t}\vert\chi_s^{\e}\vert^2\right)\dif t.
\label{I2356e}
\ee

Finally, by \eqref{chiel2}, \eqref{I2356e} and the Gronwall inequality, it immediately leads to
\ce
\mE\left(\sup_{0 \leq t \leq T}\vert\chi_t^{\e}\vert^2\right)\leq C_T(I_1^\e(T)+I_4^\e(T)).
\de
Combining this with \eqref{I1e} and \eqref{I4e} gives
$$
\lim_{\e\rightarrow0}\mE\left(\sup_{0 \leq t \leq T}\vert\chi_t^{\e}\vert^2\right)=0.
$$
The proof is complete.
\end{proof}

\subsection{Tightness}\label{tight}
 In this subsection, we intend to prove the tightness of the distribution for the solution $\vartheta^{\e}$ to Eq.\eqref{varthetaeq} in $C([0,T],\mR^n)$.

First of all, we write
$$\vartheta_t^{\e}=J_1^\e(t)+J_2^\e(t)+J_3^\e(t)+J_4^\e(t)+J_5^\e(t),
$$
where
\ce
&&J_1^\e(t):=\frac{1}{\sqrt{\e}}\int_0^t[b_1(X_s^{\e}, \sL_{X_s^{\e}},  Y_s^{\e, y_0, \sL_{\xi}}, \sL_{Y_s^{\e, \xi}})-\bar{b}_1(X_s^{\e}, \sL_{X_s^{\e}})] \dif s, \\
&&J_2^\e(t):=\int_0^t\p_x\bar{b}_1(\bar{X}_s, \sL_{\bar{X}_s})\vartheta_s^{\e}\dif s, \\
&&J_3^\e(t):=\int_0^t\mE[\p_\mu\bar{b}_1(u, \sL_{\bar{X}_s})(\bar{X}_s)\vartheta_s^{\e}]|_{u=\bar{X}_s}\dif s, \\
&&J_4^\e(t):=\int_0^t\p_x\s_1(\bar{X}_s, \sL_{\bar{X}_s})\vartheta_s^{\e}\dif B_s, \\
&&J_5^\e(t):=\int_0^t\mE[\p_\mu\s_1(u, \sL_{\bar{X}_s})(\bar{X}_s)\vartheta_s^{\e}]|_{u=\bar{X}_s}\dif B_s.
\de

\bl\label{tightlem}
Suppose that $(\mathbf{H}_{b_1, \s_1}^1)$, $(\mathbf{H}_{b_1}^2)$ and $(\mathbf{H}_{b_2, \s_2}^1)$-$(\mathbf{H}_{b_2, \s_2}^3)$ hold,
and $\s_1\in C_b^{2,(1,1)}(\mR^n\times\cP_2(\mR^n),\mR^{n \times d_1})$. Then $\{\Pi^{\e}:=(\vartheta^{\e}, J_1^\e, J_2^\e,$ $J_3^\e, J_4^\e, J_5^\e), \e\in(0,1)\}$ is tight in $C\left([0, T], \mR^{6n}\right)$.
\el
\begin{proof}
In view of the term $J_1^{\e}(t)$, from \eqref{phiito} we deduce that
\ce
J_1^\e(t)&=&\frac{1}{\sqrt{\e}}\int_0^t[b_1(X_s^{\e}, \sL_{X_s^{\e}},  Y_s^{\e, y_0, \sL_{\xi}}, \sL_{Y_s^{\e, \xi}})-\bar{b}_1(X_s^{\e}, \sL_{X_s^{\e}})] \dif s\\
&=&-\frac{1}{\sqrt{\e}}\int_0^t\cL\Psi(X_s^{\e}, \sL_{X_s^{\e}},Y_s^{\e, y_0, \sL_{\xi}}, \sL_{Y_s^{\e, \xi}})\dif s\\
&=&\sqrt{\e}\Big\{-\Psi(X_t^{\e}, \sL_{X_t^{\e}},Y_t^{\e, y_0, \sL_{\xi}}, \sL_{Y_t^{\e, \xi}})+\Psi(\varrho,\sL_{\varrho},y_0,\sL_{\xi})\no\\
&&+\mathbb{\tilde{E}}\int_0^tb_1(\tilde{X}_s^{\e}, \sL_{X_s^{\e}},\tilde{Y}_s^{\e, y_0, \sL_{\xi}}, \sL_{Y_s^{\e, \xi}})
\cdot\p_\mu\Psi(X_s^{\e}, \sL_{X_s^{\e}},Y_s^{\e, y_0, \sL_{\xi}}, \sL_{Y_s^{\e, \xi}})(\tilde{X}_s^{\e})\dif s\no\\
&&+\frac12\mathbb{\tilde{E}}\int_0^tTr\left[\s_1\s_1^*(\tilde{X}_s^{\e}, \sL_{X_s^{\e}})\cdot\p_{\tx}\p_\mu\Psi(X_s^{\e}, \sL_{X_s^{\e}},
 Y_s^{\e, y_0, \sL_{\xi}}, \sL_{Y_s^{\e, \xi}})(\tilde{X}_s^{\e})\right]\dif s\no\\
&&+\int_0^t\cL_1\Psi(X_s^{\e}, \sL_{X_s^{\e}},Y_s^{\e, y_0, \sL_{\xi}}, \sL_{Y_s^{\e, \xi}})\dif s+N_t^{\e,1}\Big\}+N_t^{\e,2}\\
&=:&J_{11}^\e(t)+N_t^{\e,2}.
\de
From analogous deduction to that in Subsection 4.4, it follows that
\be
\sup_{\e\in(0,1)}\mE\left(\sup_{t\in[0,T]}\vert J_{11}^{\e}(t)\vert^4\right)\leq C_T\e(1+\mE\vert\varrho\vert^{12}+\vert y_0\vert^{12}+\mE\vert\xi\vert^{12}),\label{Je11}
\ee
which implies that $J_{11}^\e$ converges to 0 in probability in $C\left([0, T], \mR^{n}\right)$. By Lemma 20.3 in \cite[P.185]{Huang}, for establishing the tightness of $\{\Pi^{\e}\}$, it is sufficient to verify the following conditions for $i=2,3,4,5$:
\be
\sup_{\e\in(0,1)}\mE\left(\sup_{t\in[0,T]}\vert J_i^{\e}(t)\vert^4\right)<\infty\label{Jeil4},
\ee
and
\be
\sup_{\e\in(0,1)}\mE|J_i^{\e}(t_2)-J_i^{\e}(t_1)|^{4}\leq C_T|t_2-t_1|^{2}\label{Jeit2t1l4},
\ee
where the same conditions remain valid if $J_i^{\e}(t)$ is replaced by $N_t^{\e,2}$.

We divide the remaining proof into two steps.

{\bf Step 1.} We aim to establish \eqref{Jeil4}.

Firstly, for $N_t^{\e,2}$, we deduce from \eqref{sJ6sJ7} and Lemma \ref{xtyt} that
\ce
\sup_{\e\in(0,1)}\mE\left(\sup_{t\in[0,T]}\vert N_t^{\e,2}\vert^4\right)\leq C_T(1+\vert y_0\vert^{12}+\mE\vert\xi\vert^{12}).
\de

Secondly, in order to investigate the uniform estimates of $\{J_i^{\e}\}$, $i=2,3,4,5$, we need to prove that
\be
\sup_{\e\in(0,1)}\mE\left(\sup_{t\in[0,T]}\vert \vartheta_t^{\e}\vert^4\right)
\leq C_T(1+\mE\vert\varrho\vert^{12}+\vert y_0\vert^{12}+\mE\vert\xi\vert^{12}).\label{varthetal4}
\ee
In fact, using the uniform boundedness of $\p_x \bar{b}_1$, $\p_\mu \bar{b}_1$, $\p_x \s_1$ and $\p_\mu\s_1$, together with the H\"{o}lder inequality and the Burkholder-Davis-Gundy inequality, one can conclude that
\ce
\mE\left(\sup_{t\in[0,T]}\vert \vartheta_t^{\e}\vert^4\right)
&\leq&C\mE\left(\sup_{t\in[0,T]}\vert J_1^{\e}(t)\vert^4\right)+C_T\int_0^T\mE\vert \vartheta_s^{\e}\vert^4\dif s\\
&&+C\mE\left(\int_0^T|\p_x\s_1(\bar{X}_s, \sL_{\bar{X}_s})\vartheta_s^{\e}|^2\dif s\right)^2\\
&&+C\mE\left(\int_0^T\left|\mE[\p_\mu\s_1(u, \sL_{\bar{X}_s})(\bar{X}_s)\vartheta_s^{\e}]|_{u=\bar{X}_s}\right|^2\dif s\right)^2\\
&\leq&C_T(1+\mE\vert\varrho\vert^{12}+\vert y_0\vert^{12}+\mE\vert\xi\vert^{12})+C_T\int_0^T\mE\vert \vartheta_s^{\e}\vert^4\dif s.
\de
As a result, \eqref{varthetal4} is derived by applying Gronwall's inequality. Consequently, it follows that $\{J_i^{\e}\}$, $i=2,3,4,5$ satisfy \eqref{Jeil4}.

{\bf Step 2.} We shall demonstrate \eqref{Jeit2t1l4}.

For $N_{t}^{\e,2}$, by the Burkholder-Davis-Gundy inequality, the H\"{o}lder inequality, Proposition \ref{poieqesti}, \eqref{b2s2grow} and Lemma \ref{xtyt}, we obtain
\ce
\mE|N_{t_2}^{\e,2}-N_{t_1}^{\e,2}|^4
&\leq&C\mE\left(\int_{t_1}^{t_2}\|\s_2(Y_s^{\e, y_0, \sL_{\xi}}, \sL_{Y_s^{\e, \xi}})\|^2\dif s\right)^2\\
&\leq&C(t_2-t_1)\int_{t_1}^{t_2}(1+\mE|Y_s^{\e, y_0, \sL_{\xi}}|^4+\mE|Y_s^{\e, \xi}|^4)\dif s\\
&\leq&C(1+\vert y_0\vert^{12}+\mE\vert\xi\vert^{12})|t_2-t_1|^2.
\de

Similarly, by combining \eqref{varthetal4} with the uniform boundedness of $\p_x \bar{b}_1$, $\p_\mu \bar{b}_1$, $\p_x \s_1$ and $\p_\mu \s_1$, it holds that
\ce
\mE|J_i^{\e}(t_2)-J_i^{\e}(t_1)|^4\leq C_T(1+\mE\vert\varrho\vert^{12}+\vert y_0\vert^{12}+\mE\vert\xi\vert^{12})|t_2-t_1|^2,\quad i=2,\ldots,5.
\de

As a consequence, the terms $N_{t}^{\e,2}$ and $\{J_i^{\e}\}$, $i=2,\ldots,5$, satisfy \eqref{Jeit2t1l4}. The proof is complete.
\end{proof}

\subsection{Proof of Theorem 3.2}\label{cltthproode}

In this subsection, we complete the proof of Theorem \ref{cltth}.
The main objective is to identify the weak limit of $\vartheta_t^{\e}$ in $C\left([0, T], \mR^{n}\right)$.

\bp\label{weak}
Assume that $(\mathbf{H}_{b_1, \s_1}^1)$, $(\mathbf{H}_{b_1}^2)$ and $(\mathbf{H}_{b_2, \s_2}^1)$-$(\mathbf{H}_{b_2, \s_2}^3)$ hold, and $\s_1\in C_b^{2,(1,1)}(\mR^n\times\cP_2(\mR^n),\mR^{n \times d_1})$. Then the solution $\vartheta^{\e}$ to Eq.\eqref{varthetaeq} weakly converges to the solution $U$ of Eq.\eqref{uteq} in $C\left([0, T] , \mR^{n}\right)$, as $\e\rightarrow0$.
\ep
\begin{proof}
By Lemma \ref{tightlem}, it is sufficient to demonstrate that, there exists a subsequence $\{\e_{n_k}\}_{k\geq1}$
such that $\vartheta^{\e_{n_k}}$ converges weakly to $U$ in $C\left([0, T],\mR^{n}\right)$. Let $B^{\e_{n_k}}=B$, $W^{\e_{n_k}}=W$ and assume that $\{(\vartheta^{\e_{n_k}}, J_1^{\e_{n_k}}, J_2^{\e_{n_k}}, J_3^{\e_{n_k}}, J_4^{\e_{n_k}},$ $J_5^{\e_{n_k}}, B^{\e_{n_k}}, W^{\e_{n_k}})\}_{k\geq1}$ convergence weakly to $(\vartheta, J_1, J_2, J_3, J_4, J_5, B, W)$. Then, by the Skorohod representation theorem (cf.\cite[Theorem C.1]{{bhr}}), there exists a new probability space denoted by $(\hat{\Omega},\hat{\sF},\hat{\mP})$ and a sequence
\ce
\{(\hat{\vartheta}^{\e_{n_k}}, \hat{J}_1^{\e_{n_k}}, \hat{J}_2^{\e_{n_k}}, \hat{J}_3^{\e_{n_k}}, \hat{J}_4^{\e_{n_k}}, \hat{J}_5^{\e_{n_k}},\hat{X}^{\e_{n_k}}, \hat{Y}^{\e_{n_k}, y_0, \sL_{\hat{\xi}}}, \hat{Y}^{\e_{n_k}, \hat{\xi}},\hat{\bar{X}},\hat{B}^{\e_{n_k}},\hat{W}^{\e_{n_k}})\}_{k\geq1},
\de
which has the same law as $\{(\vartheta^{\e_{n_k}}, J_1^{\e_{n_k}}, J_2^{\e_{n_k}}, J_3^{\e_{n_k}}, J_4^{\e_{n_k}}, J_5^{\e_{n_k}},X^{\e_{n_k}}, Y^{\e_{n_k}, y_0, \sL_{\xi}}, Y^{\e_{n_k}, \xi},\bar{X},$ $B^{\e_{n_k}},W^{\e_{n_k}})\}_{k\geq1}$ and the following properties hold:

(i) $(\hat{\vartheta}^{\e_{n_k}}, \hat{J}_1^{\e_{n_k}}, \hat{J}_2^{\e_{n_k}}, \hat{J}_3^{\e_{n_k}}, \hat{J}_4^{\e_{n_k}}, \hat{J}_5^{\e_{n_k}},\hat{B}^{\e_{n_k}},\hat{W}^{\e_{n_k}})$ converges to
$(\hat{\vartheta}, \hat{J}_1, \hat{J}_2, \hat{J}_3, \hat{J}_4, \hat{J}_5,\hat{B},\hat{W})$ in $C\left([0, T] ; \mR^{7n+m}\right)$, $\hat{\mP}$-a.s., as $k\rightarrow\infty$;

(ii) $(\hat{B}^{\e_{n_k}}(\hat{\o}),\hat{W}^{\e_{n_k}}(\hat{\o}))=(\hat{B}(\hat{\o}),\hat{W}(\hat{\o}))$ for all $\hat{\o}\in\hat{\Omega}$.\\
The process $(\hat{X}^{\e_{n_k}}, \hat{Y}^{\e_{n_k}, y_0, \sL_{\hat{\xi}}},\hat{Y}^{\e_{n_k},\hat{\xi}})$ is the solution to the system \eqref{orieq}, where $\e_{n_k}$ substitutes for $\e$, $(\hat{B},\hat{W})$ replaces $(B,W)$, and $(\hat{\varrho},\hat{\xi})$ replaces $(\varrho,\xi)$. Similarly, $\hat{\bar{X}}$ is the solution to \eqref{barxeq} with $\hat{B}$ replacing $B$ and $\hat{\varrho}$ replacing $\varrho$. Here, $\hat{\varrho}$ and $\hat{\xi}$ are random variables defined on $(\hat{\Omega},\hat{\sF},\hat{\mP})$, which have the same distributions as $\varrho$ and $\xi$ respectively.

Since the sequence $\{(\vartheta^{\e_{n_k}}, J_1^{\e_{n_k}}, J_2^{\e_{n_k}}, J_3^{\e_{n_k}}, J_4^{\e_{n_k}},$ $J_5^{\e_{n_k}},B,W)\}_{k\geq1}$ satisfies Eq.\eqref{varthetaeq}, it holds that
$$
\hat{\vartheta}_t^{\e_{n_k}}=\hat{J}_1^{\e_{n_k}}(t)+\hat{J}_2^{\e_{n_k}}(t)+\hat{J}_3^{\e_{n_k}}(t)+\hat{J}_4^{\e_{n_k}}(t)+\hat{J}_5^{\e_{n_k}}(t),
$$
where
\ce
&&\hat{J}_1^{\e_{n_k}}(t):=\frac{1}{\sqrt{\e_{n_k}}}\int_0^t[b_1(\hat{X}_s^{\e_{n_k}}, \sL_{\hat{X}_s^{\e_{n_k}}},  \hat{Y}_s^{\e_{n_k}, y_0, \sL_{\hat{\xi}}}, \sL_{\hat{Y}_s^{\e_{n_k},\hat{\xi}}})-\bar{b}_1(\hat{X}_s^{\e_{n_k}}, \sL_{\hat{X}_s^{\e_{n_k}}})] \dif s, \\
&&\hat{J}_2^{\e_{n_k}}(t):=\int_0^t\p_x\bar{b}_1(\hat{\bar{X}}_s, \sL_{\hat{\bar{X}}_s})\hat{\vartheta}_s^{\e_{n_k}}\dif s, \\
&&\hat{J}_3^{\e_{n_k}}(t):=\int_0^t\hat{\mE}[\p_\mu\bar{b}_1(u, \sL_{\hat{\bar{X}}_s})(\hat{\bar{X}}_s)\hat{\vartheta}_s^{\e_{n_k}}]
|_{u=\hat{\bar{X}}_s}\dif s, \\
&&\hat{J}_4^{\e_{n_k}}(t):=\int_0^t\p_x\s_1(\hat{\bar{X}}_s, \sL_{\hat{\bar{X}}_s})\hat{\vartheta}_s^{\e_{n_k}}\dif \hat{B}_s, \\
&&\hat{J}_5^{\e_{n_k}}(t):=\int_0^t\hat{\mE}[\p_\mu\s_1(u, \sL_{\hat{\bar{X}}_s})(\hat{\bar{X}}_s)\hat{\vartheta}_s^{\e_{n_k}}]
|_{u=\hat{\bar{X}}_s}\dif \hat{B}_s,
\de
and $\hat{\mE}$ denotes the expectation with respect to $\hat{\mP}$. We claim that the limits $\hat{J}_i$, corresponding to $\hat{J}_i^{\e_{n_k}}$ for $i=2,3,4,5$, take the following forms $\hat{\mP}$-a.s.:
\ce
&&\hat{J}_2(t)=\int_0^t\p_x\bar{b}_1(\hat{\bar{X}}_s, \sL_{\hat{\bar{X}}_s})\hat{\vartheta}_s\dif s, \\
&&\hat{J}_3(t)=\int_0^t\hat{\mE}[\p_\mu\bar{b}_1(u, \sL_{\hat{\bar{X}}_s})(\hat{\bar{X}}_s)\hat{\vartheta}_s]|_{u=\hat{\bar{X}}_s}\dif s, \\
&&\hat{J}_4(t)=\int_0^t\p_x\s_1(\hat{\bar{X}}_s, \sL_{\hat{\bar{X}}_s})\hat{\vartheta}_s\dif \hat{B}_s, \\
&&\hat{J}_5(t)=\int_0^t\hat{\mE}[\p_\mu\s_1(u, \sL_{\hat{\bar{X}}_s})(\hat{\bar{X}}_s)\hat{\vartheta}_s]|_{u=\hat{\bar{X}}_s}\dif \hat{B}_s.
\de

Indeed, from \eqref{varthetal4}, it follows that
\ce
\sup_{k\geq1}\hat{\mE}\left(\sup_{t\in[0,T]}\vert \hat{\vartheta}_t^{\e_{n_k}}\vert^4\right)
\leq C_T(1+\hat{\mE}\vert\hat{\varrho}\vert^{12}+\vert y_0\vert^{12}+\hat{\mE}\vert\hat{\xi}\vert^{12}).
\de
Since $\hat{\vartheta}^{\e_{n_k}}\rightarrow\hat{\vartheta}$ in $C\left([0, T], \mR^{n}\right)$, $\hat{\mP}$-a.s., by the Vitali convergence theorem we infer that
\be
\lim_{k\rightarrow\infty}\hat{\mE}\left(\sup_{t\in[0,T]}\vert \hat{\vartheta}_t^{\e_{n_k}}-\hat{\vartheta}_t\vert^2\right)=0.
\label{hatthetael2}
\ee

As for $\hat{J}_2^{\e_{n_k}}(t)$, combining the uniform boundedness of $\p_x\bar{b}_1(x,\mu)$ with \eqref{hatthetael2}, we conclude that
\ce
\lim_{k\rightarrow\infty}\hat{\mE}\left(\sup_{t\in[0,T]}\vert \hat{J}_2^{\e_{n_k}}(t)-\hat{J}_2(t)\vert^2\right)
&=&\lim_{k\rightarrow\infty}\hat{\mE}\left(\sup_{t\in[0,T]}\left|\int_0^t\p_x\bar{b}_1(\hat{\bar{X}}_s, \sL_{\hat{\bar{X}}_s})(\hat{\vartheta}_s^{\e_{n_k}}-\hat{\vartheta}_s)\dif s\right|^2\right)\\
&\leq&C_T\lim_{k\rightarrow\infty}\hat{\mE}\left(\sup_{t\in[0,T]}\vert \hat{\vartheta}_t^{\e_{n_k}}-\hat{\vartheta}_t\vert^2\right)
=0.
\de

For $\hat{J}_3^{\e_{n_k}}(t)$, owing to the uniform boundedness of $\p_\mu\bar{b}_1(x,\mu)(\tx)$ and \eqref{hatthetael2}, it holds that
\ce
\lim_{k\rightarrow\infty}\sup_{t\in[0,T]}\vert \hat{J}_3^{\e_{n_k}}(t)-\hat{J}_3(t)\vert^2
&=&\lim_{k\rightarrow\infty}\sup_{t\in[0,T]}\left|\int_0^t\hat{\mE}[\p_\mu\bar{b}_1(u, \sL_{\hat{\bar{X}}_s})(\hat{\bar{X}}_s)(\hat{\vartheta}_s^{\e_{n_k}}-\hat{\vartheta}_s)]|_{u=\hat{\bar{X}}_s}\dif s\right|^2\\
&\leq&C_T\lim_{k\rightarrow\infty}\hat{\mE}\left(\sup_{t\in[0,T]}\vert \hat{\vartheta}_t^{\e_{n_k}}-\hat{\vartheta}_t\vert^2\right)
=0.
\de

For $\hat{J}_4^{\e_{n_k}}(t)$ and $\hat{J}_5^{\e_{n_k}}(t)$, by the uniform boundedness of $\p_x\s_1(x,\mu)$, $\p_\mu\s_1(x,\mu)(\tx)$, \eqref{hatthetael2}, and the Burkholder-Davis-Gundy inequality, we obtain
\ce
\lim_{k\rightarrow\infty}\hat{\mE}\left(\sup_{t\in[0,T]}\vert \hat{J}_4^{\e_{n_k}}(t)-\hat{J}_4(t)\vert^2\right)
&\leq&C\lim_{k\rightarrow\infty}\hat{\mE}\left(\int_0^T\|\p_x\s_1(\hat{\bar{X}}_s, \sL_{\hat{\bar{X}}_s})(\hat{\vartheta}_s^{\e_{n_k}}-\hat{\vartheta}_s)\|^2\dif s\right)\\
&\leq&C_T\lim_{k\rightarrow\infty}\hat{\mE}\left(\sup_{t\in[0,T]}\vert \hat{\vartheta}_t^{\e_{n_k}}-\hat{\vartheta}_t\vert^2\right)
=0,
\de
and
\ce
&&\lim_{k\rightarrow\infty}\hat{\mE}\left(\sup_{t\in[0,T]}| \hat{J}_5^{\e_{n_k}}(t)-\hat{J}_5(t)|^2\right)\\
&\leq&C\lim_{k\rightarrow\infty}\hat{\mE}\left(\int_0^T\left|\hat{\mE}[\p_\mu\s_1(u, \sL_{\hat{\bar{X}}_s})(\hat{\bar{X}}_s)(\hat{\vartheta}_s^{\e_{n_k}}-\hat{\vartheta}_s)]|_{u=\hat{\bar{X}}_s}\right|^2
\dif s\right)\\
&\leq&C_T\lim_{k\rightarrow\infty}\hat{\mE}\left(\sup_{t\in[0,T]}\vert \hat{\vartheta}_t^{\e_{n_k}}-\hat{\vartheta}_t\vert^2\right)
=0.
\de

Next, we focus on identifying the limiting process of $\hat{J}_1^{\e_{n_k}}$ in $C\left([0, T], \mR^{n}\right)$. According to \eqref{Je11}, $J_{11}^{\e}$ converges to 0 in probability in $C\left([0, T],\mR^{n}\right)$. As a result, $\hat{J}_{1}$ must coincide with the limiting process $\hat{N}_t^2$ of $\hat{N}_t^{\e_{n_k},2}$ in $C\left([0, T],\mR^{n}\right)$, where
$$
\hat{N}_t^{\e_{n_k},2}=\int_0^t\p_y\Psi_{\s_2}(\hat{X}_s^{\e_{n_k}}, \sL_{\hat{X}_s^{\e_{n_k}}},  \hat{Y}_s^{\e_{n_k}, y_0, \sL_{\hat{\xi}}}, \sL_{\hat{Y}_s^{\e_{n_k},\hat{\xi}}})\mathrm{d}\hat{W}_s
$$
is a continuous martingale on the filtered probability space $(\hat{\Omega},\hat{\sF},\{\hat{\sF}_t\}_{t\geq0},\hat{\mP})$ with quadratic variational process
$$
\la\hat{N}_t^{\e_{n_k},2}\ra_t=\int_0^t(\p_y\Psi_{\s_2})(\p_y\Psi_{\s_2})^*(\hat{X}_s^{\e_{n_k}}, \sL_{\hat{X}_s^{\e_{n_k}}},  \hat{Y}_s^{\e_{n_k}, y_0, \sL_{\hat{\xi}}}, \sL_{\hat{Y}_s^{\e_{n_k},\hat{\xi}}})\dif s, \quad t\in[0,T],
$$
and $\hat{\sF}_t$ represents the natural filtration generated by $\{\hat{\varrho},\hat{\xi},\hat{B}_s,\hat{W}_s,s\leq t\}$.

Since $\hat{N}^{\e_{n_k},2}$ converges to $\hat{N}^2$, $\hat{\mP}$-a.s., we obtain that $\hat{N}_t^2$ is a continuous martingale on $(\hat{\Omega},\hat{\sF},\{\hat{\sF}_t\}_{t\geq0},\hat{\mP})$. Note that the quadratic variational process of $\hat{N}_t^2$ is
\be
\la\hat{N}^2\ra_t=\int_0^t\overline{(\p_y\Psi_{\s_2})(\p_y\Psi_{\s_2})^*}(\hat{\bar{X}}_s, \sL_{\hat{\bar{X}}_s})\dif s,\quad t\in[0,T].
\label{hatN2}
\ee
(See the detailed proof in the Appendix). Thus, by the martingale representation theorem (cf.\cite[Theorem 8.2]{{pz}}), there exists a probability space $(\check{\Omega},\check{\sF},\check{\mP})$, a filtration $\{\check{\sF}_t\}_{t\geq0}$ and a $n$-dimensional standard Brownian motion $\hat{V}$ defined on $(\hat{\Omega}\times\check{\Omega},\hat{\sF}\times\check{\sF},\hat{\mP}\times\check{\mP})$, which is
adapted to $\hat{\sF}_t\times\check{\sF}_t$, such that
$$
\hat{N}_t^2=\int_0^t\left(\overline{(\p_y\Psi_{\s_2})(\p_y\Psi_{\s_2})^*}\right)^{\frac12}(\hat{\bar{X}}_s, \sL_{\hat{\bar{X}}_s})\mathrm{d}\hat{V}_s, \quad t\in[0,T],
$$
where
$$
\hat{N}_t^2(\hat{\o},\check{\o})=\hat{N}_t^2(\hat{\o}),\quad\hat{\bar{X}}_s(\hat{\o},\check{\o})=\hat{\bar{X}}_s(\hat{\o}),\quad
(\hat{\o},\check{\o})\in\hat{\Omega}\times\check{\Omega}.
$$
Observe that $\la\hat{B},\hat{N}^{\e_{n_k},2}\ra_t=0$ and $\hat{N}^{\e_{n_k},2}$ converges to $\hat{N}^2$ almost surely under $\hat{\mP}\times\check{\mP}$. Consequently, $\la\hat{B},\hat{N}^2\ra_t=0$, which indicates that $\hat{V}$ is independent of $\hat{B}$.

Finally, based on the above discussion, $\hat{\vartheta}$ satisfies
\ce\left\{\begin{array}{l}
\mathrm{d}\hat{\vartheta}_t=\p_x\bar{b}_1(\hat{\bar{X}}_t, \sL_{\hat{\bar{X}}_t})\hat{\vartheta}_t\dif t
+\bar{\mE}[\p_\mu\bar{b}_1(u, \sL_{\hat{\bar{X}}_t})(\hat{\bar{X}}_t)\hat{\vartheta}_t]|_{u=\hat{\bar{X}}_t}\dif t
+\p_x\s_1(\hat{\bar{X}}_t, \sL_{\hat{\bar{X}}_t})\hat{\vartheta}_t\mathrm{d} \hat{B}_t \\
\qquad\quad+\bar{\mE}[\p_\mu\s_1(u, \sL_{\hat{\bar{X}}_t})(\hat{\bar{X}}_t)\hat{\vartheta}_t]|_{u=\hat{\bar{X}}_t}\mathrm{d} \hat{B}_t+\left(\overline{(\p_y\Psi_{\s_2})(\p_y\Psi_{\s_2})^*}\right)^{\frac12}(\hat{\bar{X}}_t, \sL_{\hat{\bar{X}}_t})\mathrm{d}\hat{V}_t,\\
\hat{\vartheta}_0=0,
\end{array}
\right.
\de
where $\bar{\mE}$ is the expectation with respect to $\hat{\mP}\times\check{\mP}$ and
$$
\hat{\vartheta}_t(\hat{\o},\check{\o})=\hat{\vartheta}_t(\hat{\o}),\quad\hat{B}_t(\hat{\o},\check{\o})=\hat{B}_t(\hat{\o}).
$$
According to \cite[Section 5.3]{{hlls}}, Eq.\eqref{uteq} has the weak uniqueness. Then $\hat{\vartheta}$ coincides with the solution $U$ to Eq.\eqref{uteq} in the sense of distribution. The proof is complete.
\end{proof}

Now, we are able to complete the proof of our main result.

{\bf Proof of Theorem \ref{cltth}:} Combining Proposition \ref{uethetae} and \ref{weak}, we immediately obtain that Theorem 3.2 holds. The proof is complete.

\section{An example}\label{example}

Now let us present an example that demonstrates the applicability of our results.
\bx\label{ex}
Consider the following slow-fast system:
\be\left\{\begin{array}{l}
\dif X_t^{\e}=\Big[sin(aX_t^{\e})+\int_{\mR}\check{x}\sL_{X_t^{\e}}(\dif\check{x})+cos(bY_t^{\e, y_0, \sL_{\xi}})
+\int_{\mR}cos(q\check{y})\sL_{Y_t^{\e, \xi}}(\dif\check{y}) \Big]\dif t\\
\qquad\qquad+\int_{\mR}\frac{|\check{x}|^3}{1+\check{x}^2}\sL_{X_t^{\e}}(\dif\check{x}) \dif B_t, \\
X_0^{\e}=\varrho, \quad 0\leq t\leq T, \\
\dif Y_t^{\e, \xi}=\frac{1}{\e} \int_{\mR}(-kY_t^{\e, \xi}+m\check{y})\sL_{Y_t^{\e, \xi}}(\dif\check{y}) \dif t+\frac{1}{\sqrt{\e}}\dif W_t, \\
Y_0^{\e, \xi}=\xi, \quad 0\leq t\leq T, \\
\dif Y_t^{\e, y_0, \sL_{\xi}}=\frac{1}{\e}\int_{\mR}(-kY_t^{\e, y_0, \sL_{\xi}}+m\check{y})\sL_{Y_t^{\e, \xi}}(\dif\check{y}) \dif t
+\frac{1}{\sqrt{\e}} \dif W_t, \\
Y_0^{\e, y_0, \sL_{\xi}}=y_0, \quad 0\leq t\leq T,
\end{array}\right.
\label{exorieq}
\ee
where $\left(B_t\right),\left(W_t\right)$ are $1$-dimensional standard Brownian motions, respectively, defined on the complete filtered probability space $(\Omega,\sF,\{\sF_t\}_{t \in[0, T]}, \mP)$ and are mutually independent. $\varrho, \xi$ are $\sF_0$-measurable Gaussian random variables and $a,b,q,k,m$ are positive constants.
Let
\ce
&&b_1(x, \mu, y, \nu)=sin(ax)+\int_{\mR}\check{x}\mu(\dif\check{x})+cos(by)+\int_{\mR}cos(q\check{y})\nu(\dif\check{y}),\\
&&\s_1(x,\mu)=\int_{\mR}\frac{|\check{x}|^3}{1+\check{x}^2}\mu(\dif\check{x}),
\de
and
$$
b_2(y, \nu)=\int_{\mR}(-ky+m\check{y})\nu(\mathrm{d}\check{y}),\quad \s_2(y, \nu)=1.
$$
For $x_i \in \mR, \mu_i \in \cP_2\left(\mR\right), y_i \in \mR$, $\nu_i \in \cP_2\left(\mR\right), i=1,2$, choosing $\pi^*\in\mathscr{C}(\mu_1,\mu_2)$ and $\bar{\pi}\in\mathscr{C}(\nu_1,\nu_2)$ such that
\ce
&&|b_1(x_1, \mu_1, y_1, \nu_1)-b_1(x_2, \mu_2, y_2, \nu_2)|^2\no\\
&\leq&4|sin(ax_1)-sin(ax_2)|^2+4\left|\int_{\mR}\check{x}_1\mu_1(\dif\check{x}_1)-\int_{\mR}\check{x}_2\mu_2(\dif\check{x}_2)\right|^2\no\\
&&+4|cos(by_1)-cos(by_2)|^2+4\left|\int_{\mR}cos(q\check{y}_1)\nu_1(\dif\check{y}_1)-\int_{\mR}cos(q\check{y}_2)\nu_2(\dif\check{y}_2)\right|^2\no\\
&\leq&4a^2|x_1-x_2|^2+4\left|\int_{\mR\times\mR}\check{x}_1\pi^*(\dif\check{x}_1,\dif\check{x}_2)
-\int_{\mR\times\mR}\check{x}_2\pi^*(\dif\check{x}_1,\dif\check{x}_2)\right|^2\no\\
&&+4b^2|y_1-y_2|^2+4\left|\int_{\mR\times\mR}cos(q\check{y}_1)\bar{\pi}(\dif\check{y}_1,\dif\check{y}_2)
-\int_{\mR\times\mR}cos(q\check{y}_2)\bar{\pi}(\dif\check{y}_1,\dif\check{y}_2)\right|^2\no\\
&\leq&4a^2|x_1-x_2|^2+4\int_{\mR\times\mR}|\check{x}_1-\check{x}_2|^2\pi^*(\dif\check{x}_1,\dif\check{x}_2)+4b^2|y_1-y_2|^2\no\\
&&+4q^2\int_{\mR\times\mR}|\check{y}_1-\check{y}_2|^2 \bar{\pi}(\dif\check{y}_1,\dif\check{y}_2).
\de
By taking the infimum over $\sC(\mu_1,\mu_2)$ and $\sC(\nu_1,\nu_2)$, and using the characterization of the $L^2$-Wasserstein metric, we deduce that
\ce
&&|b_1(x_1, \mu_1, y_1, \nu_1)-b_1(x_2, \mu_2, y_2, \nu_2)|^2\\
&\leq&4a^2|x_1-x_2|^2+4\mW_2^2(\mu_1, \mu_2)+4b^2|y_1-y_2|^2+4q^2\mW_2^2(\nu_1, \nu_2).
\de
By straightforward estimates one checks that $(\mathbf{H}^1_{b_{1}, \s_{1}})$ holds. Moreover, for $y_i \in \mR$, $\nu_i \in \cP_2(\mR), i=1,2$,
$$
|b_2(y_1, \nu_1)-b_2(y_2, \nu_2)|^2+|\s_2(y_1, \nu_1)-\s_2(y_2, \nu_2)|^2\leq 2(k^2+m^2)(|y_1-y_2|^2+\mW_2^2(\nu_1, \nu_2)),
$$
and
\ce
&&2\la y_1-y_2, b_2(y_1, \nu_1)-b_2(y_2, \nu_2)\ra+(3p-1)\|\s_2(y_1, \nu_1)-\s_2(y_2, \nu_2)\|^2\\
&\leq&2\left\la y_1-y_2,\int_{\mR}(-ky_1+m\check{y}_1)\nu_1(\dif\check{y}_1)-\int_{\mR}(-ky_2+m\check{y}_2)\nu_2(\dif\check{y}_2)\right\ra\\
&\leq&-2k|y_1-y_2|^2+2|y_1-y_2|\left|\int_{\mR}m\check{y}_1\nu_1(\dif\check{y}_1)-\int_{\mR}m\check{y}_2\nu_2(\dif\check{y}_2)\right|\\
&\leq&-k|y_1-y_2|^2+\frac{m^2}k\left|\int_{\mR}\check{y}_1\nu_1(\dif\check{y}_1)-\int_{\mR}\check{y}_2\nu_2(\dif\check{y}_2)\right|^2,
\de
which yields that
\ce
&&2\la y_1-y_2, b_2(y_1, \nu_1)-b_2(y_2, \nu_2)\ra+(3p-1)|\s_2(y_1, \nu_1)-\s_2(y_2, \nu_2)|^2\\
&\leq& -k|y_1-y_2|^2+\frac{m^2}k\mW_2^2(\nu_1, \nu_2).
\de
Let $k=\frac{1}{24p}$, $m=\frac{1}{48p}$, $(\mathbf{H}^1_{b_{2}, \s_{2}})$ and $(\mathbf{H}^2_{b_{2}, \s_{2}})$ are satisfied.
Note that
\ce
&&\p_\mu b_1(x, \mu, y, \nu)(\tx)=1,\quad\p_\nu b_1(x, \mu, y, \nu)(\ty)=-qsin(q\ty), \\
&&\p_{\tx}\p_\mu b_1(x, \mu, y, \nu)(\tx)=0, ~~ \p_\mu\p_y b_1(x, \mu, y, \nu)(\tx)=0,~~\p_{\ty}\p_\nu b_1(x, \mu, y, \nu)(\ty)=-q^2cos(q\ty),\\
&&\p_\mu\s_1(x,\mu)(\tx)=\frac{3\tx^2+\tx^4}{(1+\tx^2)^2}\cdot sgn(\tx),\quad\p_{\tx}\p_\mu\s_1(x,\mu)(\tx)=\frac{2|\tx|(3-\tx^2)}{(1+\tx^2)^3},\\
&&\p_\nu b_2(y, \nu)(\ty)=m,\quad\p_\nu\s_2(y, \nu)(\ty)=0,\\
&&\p_{\ty}\p_\nu b_2(y, \nu)(\ty)=\p_{\ty}\p_\nu\s_2(y, \nu)(\ty)=0.
\de
A straightforward verification shows that $(\mathbf{H}^2_{b_{1}})$ and $(\mathbf{H}^3_{b_{2}, \s_{2}})$ hold, and $\s_1\in C_b^{2,(1,1)}(\mR\times\cP_2(\mR),\mR)$. Consequently, the conclusions of Theorems \ref{xbarxp} and \ref{cltth} follow. 

Besides, we remind that for any $\tx>0$, 
$$
\p_{\tx}^2\p_\mu\s_1(x,\mu)(\tx)=\frac{6(\tx^4-6\tx^2+1)}{(1+\tx^2)^4},
$$
and for any $\tx<0$,
$$
\p_{\tx}^2\p_\mu\s_1(x,\mu)(\tx)=-\frac{6(\tx^4-6\tx^2+1)}{(1+\tx^2)^4}.
$$
That is, $\p_{\tx}^2\p_\mu\s_1(x,\mu)(\tx)$ is discontinuous at $\tx=0$. So, $\s_1$ does not satisfy the assumption in \cite[Theorem 2.3]{lx}. Therefore, our conditions for the central limit theorem are more general.
\ex

\section{Appendix}\label{append}

In this section, we prove \eqref{parxb1s0} and \eqref{hatN2}.

\subsection{Proof of \eqref{parxb1s0}}\label{parxb1sopro}

Let $C([0,\infty),\mR^m)$ be the collection of continuous functions from $[0,\infty)$ to $\mR^m$, equipped with the uniform convergence topology. Set
\ce
&&\tilde{\Omega}:=C([0,\infty),\mR^m)\times C([0,\infty),\cP_2\left(\mR^m\right)),\\
&&\tilde{\sF}:=\sB\Big(C([0,\infty),\mR^m)\Big)\times\sB\Big(C([0,\infty),\cP_2\left(\mR^m\right))\Big),\\
&&\tilde{\sF}_t:=\s(M_r:0\leq r\leq t),\quad t\geq0.
\de
By a similar argument to that in \cite[Theorem 4.11]{rrw}, there exists a unique probability measure $\tilde{\mP}$ on $(\tilde{\Omega},\tilde{\sF})$ such that $M_\cdot$ is a Markov process with respect to $(\tilde{\sF}_t)$ and has the transition function $\{\sL_{Y_{t}^{y,\nu}}\times \d_{\sL_{Y_{t}^{\xi}}}:t\geq0,(y,\nu)\in\mR^m\times\cP_2(\mR^m)\}$. Note that
$$
\sL_{M_{t}}^{\tilde{\mP}}:=\tilde{\mP}\circ M_t^{-1}=\sL_{Y_{t}^{y,\nu}}\times \d_{\sL_{Y_{t}^{\xi}}}.
$$
Thus, it holds that
\be
\tilde{b}_{1,s_0}(x,\mu,y,\nu,s)
&=&\hat{b}_1(x,\mu,y,\nu,s)-\hat{b}_1(x,\mu,y,\nu,s+s_0)\no\\
&=&\hat{b}_1(x,\mu,y,\nu,s)-\mE b_1(x,\mu, Y_{s+s_0}^{y,\nu}, \sL_{Y_{s+s_0}^{\xi}})\no\\
&=&\hat{b}_1(x,\mu,y,\nu,s)-\tilde{\mE}[\tilde{\mE}[b_1(x,\mu,M_{s+s_0})|\tilde{\sF}_{s_0}]]\no\\
&=&\hat{b}_1(x,\mu,y,\nu,s)-\tilde{\mE}[\tilde{\mE}^{M_{s_0}}[b_1(x,\mu,M_{s})]]\no\\
&=&\hat{b}_1(x,\mu,y,\nu,s)-\mE \hat{b}_1(x,\mu, Y_{s_0}^{y,\nu}, \sL_{Y_{s_0}^{\xi}},s).\label{tildeb1s0}
\ee
Then we obtain
\be
\p_x\tilde{b}_{1,s_0}(x,\mu,y,\nu,s)=\p_x\hat{b}_1(x,\mu,y,\nu,s)-\p_x\mE \hat{b}_1(x,\mu, Y_{s_0}^{y,\nu}, \sL_{Y_{s_0}^{\xi}},s).\label{parxtildeb1s0}
\ee
By ($\mathbf{H}_{b_1}^2$), \eqref{w2lyxi1lyxi2} and Lemma \ref{y1nu1y2nu2}, we get
\be
&&|\p_x\hat{b}_1(x,\mu,y_1,\nu_1,s)-\p_x\hat{b}_1(x,\mu,y_2,\nu_2,s)|\no\\
&\leq&\mE\|\p_xb_1(x,\mu, Y_s^{y_1,\nu_1},\sL_{Y_s^{\xi_1}})-\p_xb_1(x,\mu, Y_s^{y_2,\nu_2}, \sL_{Y_s^{\xi_2}})\|\no\\
&\leq&C\mE\Big(|Y_s^{y_1,\nu_1}-Y_s^{y_2,\nu_2}|+\mW_2(\sL_{Y_s^{\xi_1}}, \sL_{Y_s^{\xi_2}})\Big)\no\\
&\leq&C\Big[|y_1-y_2|e^{-\frac{\b_1}2s}+\mW_2(\nu_1, \nu_2)e^{-\frac{\b_1-\b_2}2s}\Big],\label{parxhatb1ynv}
\ee
where $\nu_1=\sL_{\xi_1}$, $\nu_2=\sL_{\xi_2}$. Combining \eqref{parxtildeb1s0}, \eqref{parxhatb1ynv} and Lemma \ref{ymul2}, we obtain
\ce
&&\|\p_x\tilde{b}_{1,s_0}(x,\mu,y,\nu,s)\|\\
&\leq&C\mE\Big[|y-Y_{s_0}^{y,\nu}|e^{-\frac{\b_1}2s}+(\mE|\xi-Y_{s_0}^{\xi}|^2)^{\frac{1}2}e^{-\frac{\b_1-\b_2}2s}\Big]\\
&\leq&Ce^{-\frac{\b_1-\b_2}2s}(1+|y|^2+\nu(|\cdot|^2))^{1/2}.
\de

\subsection{Proof of \eqref{hatN2}}\label{hatN2pro}
Without loss of generality, we will use $(\Omega,\sF,\{\sF_t\}_{t\geq0},\mP)$ instead of $(\hat{\Omega},\hat{\sF},\{\hat{\sF}_t\}_{t\geq0},$ $\hat{\mP})$, and $(N^{\e,2},N^2,X^{\e},Y^{\e, y_0, \sL_{\xi}},Y^{\e,\xi},\bar{X})$ instead of $(\hat{N}^{\e_{n_k},2},\hat{N}^2,\hat{X}^{\e_{n_k}},$ $\hat{Y}^{\e_{n_k}, y_0, \sL_{\hat{\xi}}},\hat{Y}^{\e_{n_k},\hat{\xi}},\hat{\bar{X}})$ below.

We define the functions
\ce
&&H(x,\mu,y,\nu)=[(\p_y\Psi_{\s_2})(\p_y\Psi_{\s_2})^*](x,\mu,y,\nu),\\
&&\bar{H}(x,\mu)=\int_{\mR^m\times\cP_2(\mR^m)}H(x,\mu,y,\nu)(\eta\times\d_\eta)(\dif y, \dif\nu).
\de
Since $N^{\e,2}$ converges to $N^2$, $\mP$-a.s., we aim to prove
\ce
\la N^2\ra_t=\int_0^t\overline{(\p_y\Psi_{\s_2})(\p_y\Psi_{\s_2})^*}(\bar{X}_s, \sL_{\bar{X}_s})\dif s,\quad t\in[0,T].
\de
We claim:
\be
\lim_{\e\rightarrow0}\mE\left\|\int_0^tH(X_s^{\e}, \sL_{X_s^{\e}},  Y_s^{\e, y_0, \sL_{\xi}}, \sL_{Y_s^{\e, \xi}})\dif s
-\int_0^t\bar{H}(\bar{X}_s, \sL_{\bar{X}_s})\dif s\right\|^2=0.\label{HbarHl2}
\ee
Actually, if \eqref{HbarHl2} holds, there exists a subsequence $\{\e_k\}_{k\geq1}$
 tending to 0 such that
\ce
\int_0^tH(X_s^{\e_k}, \sL_{X_s^{\e_k}},  Y_s^{\e_k, y_0, \sL_{\xi}}, \sL_{Y_s^{\e_k, \xi}})\dif s
\rightarrow\int_0^t\bar{H}(\bar{X}_s, \sL_{\bar{X}_s})\dif s,\quad \mP\text{-a.s., as} \quad k\rightarrow\infty.
\de
 Moreover, note that $N^{\e_k,2}_t\otimes N^{\e_k,2}_t-\int_0^tH(X_s^{\e_k}, \sL_{X_s^{\e_k}},  Y_s^{\e_k, y_0, \sL_{\xi}}, \sL_{Y_s^{\e_k, \xi}})\dif s$ is a matrix-valued martingale. Using the Vitali convergence theorem, we conclude that for any $0\leq r\leq t\leq T$,
\ce
&&\mE\Big[\Big(N^2_t\otimes N^2_t-\int_0^t\bar{H}(\bar{X}_s, \sL_{\bar{X}_s})\dif s\Big)-\Big(N^2_r\otimes N^2_r-\int_0^r\bar{H}(\bar{X}_s, \sL_{\bar{X}_s})\dif s\Big)\Big|\sF_r\Big]\\
&=&\lim_{k\rightarrow\infty}\mE\Big[\Big(N^{\e_k,2}_t\otimes N^{\e_k,2}_t-\int_0^tH(X_s^{\e_k}, \sL_{X_s^{\e_k}},  Y_s^{\e_k, y_0, \sL_{\xi}}, \sL_{Y_s^{\e_k, \xi}})\dif s\Big)\\
&&-\Big(N^{\e_k,2}_r\otimes N^{\e_k,2}_r-\int_0^rH(X_s^{\e_k}, \sL_{X_s^{\e_k}},  Y_s^{\e_k, y_0, \sL_{\xi}}, \sL_{Y_s^{\e_k, \xi}})\dif s\Big)\Big|\sF_r\Big]\\
&=&0,\quad\mP\text{-a.s.},
\de
which implies that $N^2_t\otimes N^2_t-\int_0^t\bar{H}(\bar{X}_s, \sL_{\bar{X}_s})\dif s$ is a matrix-valued martingale.
Thus, we conclude that
\ce
\la N^2\ra_t=\int_0^t\overline{(\p_y\Psi_{\s_2})(\p_y\Psi_{\s_2})^*}(\bar{X}_s, \sL_{\bar{X}_s})\dif s,\quad t\in[0,T].
\de

\subsubsection{Proof of the claim \eqref{HbarHl2}}\label{HbarHl2pro}

We divide the proof into three steps. The first step is to prove the local Lipschitz continuity of $H(x,\mu,y,\nu)$. The second step is to demonstrate \eqref{HbarHl2}. The third step is to supplement the proof of estimate \eqref{su1} required in the second step.

{\bf Step 1.} We prove that $H(x,\mu,y,\nu)$ is locally Lipschitz continuous, and $\bar{H}(x,\mu)$ is Lipschitz continuous.

\eqref{parphic} implies that $\p_y\Psi(x,\mu,y,\nu)$ is bounded and Lipschitz continuous w.r.t. $(x,\mu,y,\nu)$, thus
\be
&&\|\p_y\Psi_{\s_2}(x_1,\mu_1,y_1,\nu_1)-\p_y\Psi_{\s_2}(x_2,\mu_2,y_2,\nu_2)\|\no\\
&\leq&\|\p_y\Psi(x_1,\mu_1,y_1,\nu_1)-\p_y\Psi(x_2,\mu_2,y_2,\nu_2)\|\cdot\|\s_2(y_1,\nu_1)\|\no\\
&&+\|\p_y\Psi(x_2,\mu_2,y_2,\nu_2)\|\cdot\|\s_2(y_1,\nu_1)-\s_2(y_2,\nu_2)\|\no\\
&\leq&C(|x_1-x_2|+\mW_2(\mu_1, \mu_2)+|y_1-y_2|+\mW_2(\nu_1, \nu_2))(1+|y_1|^2+\nu_1(|\cdot|^2))^{1/2},\label{parpsis2loclip}
\ee
and
\be
\|\p_y\Psi_{\s_2}(x,\mu,y,\nu)\|\leq \|\p_y\Psi(x,\mu,y,\nu)\|\cdot\|\s_2(y,\nu)\|\leq C(1+|y|^2+\nu(|\cdot|^2))^{1/2}.\label{parpsis2grow}
\ee
Similarly, we have
\be
&&\|\p_y\Psi_{\s_2}(x_1,\mu_1,y_1,\nu_1)-\p_y\Psi_{\s_2}(x_2,\mu_2,y_2,\nu_2)\|\no\\
&\leq& C(|x_1-x_2|+\mW_2(\mu_1, \mu_2)+|y_1-y_2|+\mW_2(\nu_1, \nu_2))(1+|y_2|^2+\nu_2(|\cdot|^2))^{1/2}.\label{parpsis2loclip2}
\ee
By \eqref{parpsis2loclip}-\eqref{parpsis2loclip2}, it holds that
\be
&&\|H(x_1,\mu_1,y_1,\nu_1)-H(x_2,\mu_2,y_2,\nu_2)\|\no\\
&\leq&\|\p_y\Psi_{\s_2}(x_1,\mu_1,y_1,\nu_1)-\p_y\Psi_{\s_2}(x_2,\mu_2,y_2,\nu_2)\|\cdot\|\p_y\Psi_{\s_2}(x_1,\mu_1,y_1,\nu_1)\|\no\\
&&+\|\p_y\Psi_{\s_2}(x_2,\mu_2,y_2,\nu_2)\|\cdot\|\p_y\Psi_{\s_2}(x_1,\mu_1,y_1,\nu_1)-\p_y\Psi_{\s_2}(x_2,\mu_2,y_2,\nu_2)\|\no\\
&\leq& C(1+|y_1|^2+\nu_1(|\cdot|^2)+|y_2|^2+\nu_2(|\cdot|^2))\no\\
&&\cdot(|x_1-x_2|+\mW_2(\mu_1, \mu_2)+|y_1-y_2|+\mW_2(\nu_1, \nu_2)),\label{Hloclip}
\ee
i.e. $H(x,\mu,y,\nu)$ is locally Lipschitz continuous.
By \eqref{b2s2grow} and the uniform boundedness of $\p_y\Psi(x,\mu,y,\nu)$, we have
\be
&&\|H(x,\mu,y,\nu)\|=\|[(\p_y\Psi_{\s_2})(\p_y\Psi_{\s_2})^*](x,\mu,y,\nu)\|\leq C(1+|y|^2+\nu(|\cdot|^2)).\label{Hgrow}
\ee

Combining Lemma \ref{ymul2}, Lemma \ref{y1nu1y2nu2} and \eqref{Hloclip}, we obtain
\be
&&\left\|\mE H(x,\mu, Y_t^{y_0, \sL_{\xi}}, \sL_{Y_t^{\xi}})-\bar{H}(x,\mu)\right\|\no\\
&=&\left\|\mE H(x,\mu, Y_t^{y_0, \sL_{\xi}}, \sL_{Y_t^{\xi}})-\int_{\mR^m\times\cP_2(\mR^m)}H(x,\mu,y,\nu)(\eta\times\d_\eta)(\dif y, \dif\nu)\right\|\no\\
&\overset{\sL_{\zeta}=\nu}{=}&\left\|\mE H(x,\mu, Y_t^{y_0, \sL_{\xi}}, \sL_{Y_t^{\xi}})
-\int_{\mR^m\times\cP_2(\mR^m)}\mE H(x,\mu, Y_t^{y,\nu}, \sL_{Y_t^{\zeta}})(\eta\times\d_\eta)(\dif y, \dif\nu)\right\|\no\\
&\leq&\int_{\mR^m\times\cP_2(\mR^m)}\mE\left\|H(x,\mu, Y_t^{y_0, \sL_{\xi}}, \sL_{Y_t^{\xi}})-H(x,\mu, Y_t^{y,\nu}, \sL_{Y_t^{\zeta}})\right\|(\eta\times\d_\eta)(\dif y, \dif\nu)\no\\
&\overset{\sL_{\tilde{\zeta}}=\eta}{=}&\int_{\mR^m}\mE\left\|H(x,\mu, Y_t^{y_0, \sL_{\xi}}, \sL_{Y_t^{\xi}})-H(x,\mu, Y_t^{y,\eta}, \sL_{Y_t^{\tilde{\zeta}}})\right\|\eta(\dif y)\no\\
&\leq&C\int_{\mR^m}\Big[1+(\mE|Y_t^{y_0,\sL_{\xi}}|^4)^{\frac12}+\mE|Y_t^{\xi}|^2+(\mE|Y_t^{y,\eta}|^4)^{\frac12}+\mE|Y_t^{\tilde{\zeta}}|^2\Big]\no\\
&&\qquad\qquad\qquad\qquad\cdot\Big[(\mE|Y_t^{y_0, \sL_{\xi}}-Y_t^{y,\eta}|^2)^{\frac12}+(\mE|Y_t^{\xi}-Y_t^{\tilde{\zeta}}|^2)^{\frac12}
\Big]\eta(\dif y)\no\\
&\leq&C\int_{\mR^m}\Big[1+|y_0|^2+(\mE|\xi|^4)^{\frac12}+|y|^2+(\mE|\tilde{\zeta}|^4)^{\frac12}\Big]\no\\
&&\qquad\qquad\qquad\cdot\Big[|y_0-y|e^{-\frac{\b_1}{2}t}+(\mE|\xi-\tilde{\zeta}|^2)^{\frac12}e^{-\frac{\b_1-\b_2}{2}t}\Big]\eta(\dif y)\no\\
&\leq&Ce^{-\frac{\b_1-\b_2}{2}t}(1+|y_0|^4+\mE|\xi|^4).\label{HbarH}
\ee
By integrating \eqref{Hloclip}-\eqref{HbarH} and Lemma \ref{ymul2}, we attain
\ce
&&\|\bar{H}(x_1,\mu_1)-\bar{H}(x_2,\mu_2)\|\\
&\leq&\|\mE H(x_1,\mu_1, Y_t^{y_0, \sL_{\xi}}, \sL_{Y_t^{\xi}})-\bar{H}(x_1,\mu_1)\|+\|\mE H(x_2,\mu_2, Y_t^{y_0, \sL_{\xi}}, \sL_{Y_t^{\xi}})-\bar{H}(x_2,\mu_2)\|\\
&&+\|\mE H(x_1,\mu_1, Y_t^{y_0, \sL_{\xi}}, \sL_{Y_t^{\xi}})-\mE H(x_2,\mu_2, Y_t^{y_0, \sL_{\xi}}, \sL_{Y_t^{\xi}})\|\\
&\leq&Ce^{-\frac{\b_1-\b_2}{2}t}(1+|y_0|^4+\mE|\xi|^4)\\
&&+C\Big[\Big(1+|y_0|^2e^{-\a_1 t}+\mE|\xi|^2e^{-(\a_1-\a_2) t}\Big)\cdot(|x_1-x_2|+\mW_2(\mu_1, \mu_2))\Big],
\de
and
\ce
\|\bar{H}(x,\mu)\|
&=&\left\|\int_{\mR^m\times\cP_2(\mR^m)}H(x,\mu,y,\nu)(\eta\times\d_\eta)(\dif y, \dif\nu)\right\|\\
&\overset{\sL_{\zeta}=\nu}{=}&\left\|\int_{\mR^m\times\cP_2(\mR^m)}\mE H(x,\mu, Y_t^{y,\nu}, \sL_{Y_t^{\zeta}})(\eta\times\d_\eta)(\dif y, \dif\nu)\right\|\\
&\overset{\sL_{\tilde{\zeta}}=\eta}{\leq}&\int_{\mR^m}\mE \left\|H(x,\mu, Y_t^{y,\eta}, \sL_{Y_t^{\tilde{\zeta}}})\right\|\eta(\dif y)\\
&\leq&C\int_{\mR^m}\Big(1+|y|^2e^{-\a_1 t}+\mE|\tilde{\zeta}|^2e^{-(\a_1-\a_2) t}\Big)\eta(\dif y)\\
&\leq&C+Ce^{-(\a_1-\a_2) t}(1+\mE|\xi|^2).
\de
Letting $t\rightarrow\infty$, we obtain
\be
\|\bar{H}(x_1,\mu_1)-\bar{H}(x_2,\mu_2)\|\leq C(|x_1-x_2|+\mW_2(\mu_1, \mu_2)),\label{barHlip}
\ee
and
\be
\|\bar{H}(x,\mu)\|\leq C.\label{barHboun}
\ee

{\bf Step 2.} We prove \eqref{HbarHl2}.

Note that
\ce
X_t^{\e}-X_{t(\triangle)}^{\e}=\int_{t(\triangle)}^tb_1(X_s^{\e}, \sL_{X_s^{\e}}, Y_s^{\e, y_0, \sL_{\xi}}, \sL_{Y_s^{\e, \xi}}) \dif s
+\int_{t(\triangle)}^t\s_1(X_s^{\e}, \sL_{X_s^{\e}}) \dif B_s,
\de
where $t(\triangle):=\[\frac{t}{\triangle}\]\triangle$ and $\[\frac{t}{\triangle}\]$ denotes the integer part of $\frac{t}{\triangle}$.
By \eqref{b1s1grow}, the H\"{o}lder inequality, the Burkholder-Davis-Gundy inequality and Lemma \ref{xtyt}, we have
\be
&&\mE|X_t^{\e}-X_{t(\triangle)}^{\e}|^4\no\\
&\leq&C\triangle^3\mE\int_{t(\triangle)}^t\left|b_1(X_s^{\e}, \sL_{X_s^{\e}}, Y_s^{\e, y_0, \sL_{\xi}}, \sL_{Y_s^{\e, \xi}})\right|^4 \dif s
+C\mE\left[\int_{t(\triangle)}^t\left\|\s_1(X_s^{\e}, \sL_{X_s^{\e}})\right\|^2 \dif s\right]^2\no\\
&\leq&C(\triangle^3+\triangle)\int_{t(\triangle)}^t\left[1+\mE|X_s^{\e}|^4+\mE|Y_s^{\e, y_0, \sL_{\xi}}|^4+\mE|Y_s^{\e, \xi}|^4\right]\dif s\no\\
&\leq&C_T(\triangle^4+\triangle^2)(1+\mE\vert\varrho\vert^{12}+\vert y_0\vert^{12}+\mE\vert\xi\vert^{12}).\label{xextriel4}
\ee

Notice that
\be
&&\int_0^t\left[H(X_s^{\e}, \sL_{X_s^{\e}}, Y_s^{\e, y_0, \sL_{\xi}}, \sL_{Y_s^{\e, \xi}})-\bar{H}(\bar{X}_s, \sL_{\bar{X}_s})\right]\dif s\no\\
&=&\int_0^t\left[H(X_s^{\e}, \sL_{X_s^{\e}}, Y_s^{\e, y_0, \sL_{\xi}}, \sL_{Y_s^{\e, \xi}})-H(X_{s(\triangle)}^{\e}, \sL_{X_{s(\triangle)}^{\e}}, Y_s^{\e, y_0, \sL_{\xi}}, \sL_{Y_s^{\e, \xi}})\right]\dif s\no\\
&&+\int_0^t\left[H(X_{s(\triangle)}^{\e}, \sL_{X_{s(\triangle)}^{\e}}, Y_s^{\e, y_0, \sL_{\xi}}, \sL_{Y_s^{\e, \xi}})
-\bar{H}(X_{s(\triangle)}^{\e}, \sL_{X_{s(\triangle)}^{\e}})\right]\dif s\no\\
&&+\int_0^t\left[\bar{H}(X_{s(\triangle)}^{\e}, \sL_{X_{s(\triangle)}^{\e}})-\bar{H}(X_s^{\e}, \sL_{X_s^{\e}})\right]\dif s\no\\
&&+\int_0^t\left[\bar{H}(X_s^{\e}, \sL_{X_s^{\e}})-\bar{H}(\bar{X}_s, \sL_{\bar{X}_s})\right]\dif s\no\\
&=:&\mathscr{O}_1(t)+\mathscr{O}_2(t)+\mathscr{O}_3(t)+\mathscr{O}_4(t).\label{Hyey0xibarH}
\ee

For $\mathscr{O}_1(t)$, based on the H\"{o}lder inequality, \eqref{Hloclip}, \eqref{xextriel4} and Lemma \ref{xtyt}, one can obtain that
\be
&&\mE\|\mathscr{O}_1(t)\|^2\no\\
&\leq&C_T\mE\int_0^t\left(1+|Y_s^{\e, y_0, \sL_{\xi}}|^2+\sL_{Y_s^{\e, \xi}}(|\cdot|^2)\right)^2\left(|X_s^{\e}-X_{s(\triangle)}^{\e}|
+\mW_2(\sL_{X_s^{\e}},\sL_{X_{s(\triangle)}^{\e}})\right)^2\dif s\no\\
&\leq&C_T\int_0^t\left(1+\mE|Y_s^{\e, y_0, \sL_{\xi}}|^8+\mE|Y_s^{\e, \xi}|^8\right)^{\frac12}
\cdot\left(\mE|X_s^{\e}-X_{s(\triangle)}^{\e}|^4\right)^{\frac12}\dif s\no\\
&\leq&C_T(\triangle^2+\triangle)(1+\mE\vert\varrho\vert^{12}+\vert y_0\vert^{12}+\mE\vert\xi\vert^{12}).\label{so1l2}
\ee

As for $\mathscr{O}_3(t)$ and $\mathscr{O}_4(t)$,  by the H\"{o}lder inequality, \eqref{barHlip}, \eqref{xextriel4} and Theorem \ref{xbarxp}, we get
\be
\mE\|\mathscr{O}_3(t)\|^2+\mE\|\mathscr{O}_4(t)\|^2
&\leq&C_T\int_0^t\left[\mE|X_s^{\e}-X_{s(\triangle)}^{\e}|^2+\mE|X_s^{\e}-\bar{X}_s|^2\right]\dif s\no\\
&\leq&C_T(\triangle^2+\triangle+\e)(1+\mE\vert\varrho\vert^{12}+\vert y_0\vert^{12}+\mE\vert\xi\vert^{12}).\label{so34l2}
\ee

Now we focus on the term $\mathscr{O}_2(t)$. Observe that
\be
&&\mE\|\mathscr{O}_2(t)\|^2\no\\
&=&\mE\Bigg\|\sum_{k=0}^{[t/\triangle]-1}\int_{k\triangle}^{(k+1)\triangle}
\left[H(X_{s(\triangle)}^{\e}, \sL_{X_{s(\triangle)}^{\e}}, Y_s^{\e, y_0, \sL_{\xi}}, \sL_{Y_s^{\e, \xi}})
-\bar{H}(X_{s(\triangle)}^{\e}, \sL_{X_{s(\triangle)}^{\e}})\right]\dif s\no\\
&&\qquad+\int_{t(\triangle)}^t
\left[H(X_{s(\triangle)}^{\e}, \sL_{X_{s(\triangle)}^{\e}}, Y_s^{\e, y_0, \sL_{\xi}}, \sL_{Y_s^{\e, \xi}})
-\bar{H}(X_{s(\triangle)}^{\e}, \sL_{X_{s(\triangle)}^{\e}})\right]\dif s\Bigg\|^2\no\\
&\leq&\frac{C_T}{\triangle}\sum_{k=0}^{[t/\triangle]-1}\mE\left\|\int_{k\triangle}^{(k+1)\triangle}
\left[H(X_{s(\triangle)}^{\e}, \sL_{X_{s(\triangle)}^{\e}}, Y_s^{\e, y_0, \sL_{\xi}}, \sL_{Y_s^{\e, \xi}})
-\bar{H}(X_{s(\triangle)}^{\e}, \sL_{X_{s(\triangle)}^{\e}})\right]\dif s\right\|^2\no\\
&&+2\mE\left\|\int_{t(\triangle)}^t\left[H(X_{s(\triangle)}^{\e}, \sL_{X_{s(\triangle)}^{\e}}, Y_s^{\e, y_0, \sL_{\xi}}, \sL_{Y_s^{\e, \xi}})
-\bar{H}(X_{s(\triangle)}^{\e}, \sL_{X_{s(\triangle)}^{\e}})\right]\dif s\right\|^2\no\\
&=:&\mathscr{U}_1(t)+\mathscr{U}_2(t).\label{so2l2}
\ee

For $\mathscr{U}_1$, it holds that
\be
\mathscr{U}_1(t)\leq C_T(1+\vert y_0\vert^{12}+\mE\vert\xi\vert^{12})\frac{\e}{\triangle}.\label{su1}
\ee

For $\mathscr{U}_2$, combining the H\"{o}lder inequality, \eqref{Hgrow}, \eqref{barHboun} and Lemma \ref{xtyt}, we obtain
\be
\mathscr{U}_2(t)
&\leq&C\triangle\mE\int_{t(\triangle)}^t\left(\|H(X_{s(\triangle)}^{\e}, \sL_{X_{s(\triangle)}^{\e}}, Y_s^{\e, y_0, \sL_{\xi}}, \sL_{Y_s^{\e, \xi}})\|^2+\|\bar{H}(X_{s(\triangle)}^{\e}, \sL_{X_{s(\triangle)}^{\e}})\|^2\right)\dif s\no\\
&\leq&C\triangle\mE\int_{t(\triangle)}^t\left(1+|Y_s^{\e, y_0, \sL_{\xi}}|^4+\mE|Y_s^{\e, \xi}|^4\right)\dif s\no\\
&\leq&C\triangle^2(1+\vert y_0\vert^{12}+\mE\vert\xi\vert^{12}).\label{su2}
\ee

Finally, \eqref{Hyey0xibarH}-\eqref{su2} yield that
\ce
&&\mE\left\|\int_0^tH(X_s^{\e}, \sL_{X_s^{\e}},  Y_s^{\e, y_0, \sL_{\xi}}, \sL_{Y_s^{\e, \xi}})\dif s-\int_0^t\bar{H}(\bar{X}_s, \sL_{\bar{X}_s})\dif s\right\|^2\\
&\leq&C_T(\triangle^2+\triangle+\e+\frac{\e}{\triangle})(1+\mE\vert\varrho\vert^{12}+\vert y_0\vert^{12}+\mE\vert\xi\vert^{12}).
\de
Taking $\triangle=\e^{1/2}$, we infer that \eqref{HbarHl2} holds.

{\bf Step 3.} We now prove the estimate \eqref{su1}.

Note that
\be
&&\mathscr{U}_1(t)\no\\
&\leq&\frac{C_T}{\triangle}\sum_{k=0}^{[T/\triangle]-1}\mE\left\|\int_{k\triangle}^{(k+1)\triangle}
\left[H(X_{k\triangle}^{\e}, \sL_{X_{k\triangle}^{\e}}, Y_s^{\e, y_0, \sL_{\xi}}, \sL_{Y_s^{\e, \xi}})
-\bar{H}(X_{k\triangle}^{\e}, \sL_{X_{k\triangle}^{\e}})\right]\dif s\right\|^2\no\\
&\leq&\frac{C_T}{\triangle^2}\sup_{0\leq k\leq [T/\triangle]-1}\mE\left\|\int_{0}^{\triangle}
\left[H(X_{k\triangle}^{\e}, \sL_{X_{k\triangle}^{\e}}, Y_{s+k\triangle}^{\e, y_0, \sL_{\xi}}, \sL_{Y_{s+k\triangle}^{\e, \xi}})
-\bar{H}(X_{k\triangle}^{\e}, \sL_{X_{k\triangle}^{\e}})\right]\dif s\right\|^2\no\\
&=&\frac{C_T\e^2}{\triangle^2}\sup_{0\leq k\leq [T/\triangle]-1}\mE\left\|\int_{0}^{\frac{\triangle}{\e}}
\left[H(X_{k\triangle}^{\e}, \sL_{X_{k\triangle}^{\e}}, Y_{s\e+k\triangle}^{\e, y_0, \sL_{\xi}}, \sL_{Y_{s\e+k\triangle}^{\e, \xi}})
-\bar{H}(X_{k\triangle}^{\e}, \sL_{X_{k\triangle}^{\e}})\right]\dif s\right\|^2\no\\
&\leq&\frac{C_T\e^2}{\triangle^2}\sup_{0\leq k\leq [T/\triangle]-1}\int_{0}^{\frac\triangle\e}\int_{r}^{\frac\triangle\e}\Theta(s,r)\dif s\dif r,
\label{su1esti}
\ee
where for any $0< r< s\leq\triangle/{\e}$,
\ce
\Theta(s,r):&=&\mE\Big\la
H(X_{k\triangle}^{\e}, \sL_{X_{k\triangle}^{\e}}, Y_{s\e+k\triangle}^{\e, y_0, \sL_{\xi}}, \sL_{Y_{s\e+k\triangle}^{\e, \xi}})-\bar{H}(X_{k\triangle}^{\e}, \sL_{X_{k\triangle}^{\e}}),\\
&&\qquad\qquad\qquad H(X_{k\triangle}^{\e}, \sL_{X_{k\triangle}^{\e}}, Y_{r\e+k\triangle}^{\e, y_0, \sL_{\xi}}, \sL_{Y_{r\e+k\triangle}^{\e, \xi}})-\bar{H}(X_{k\triangle}^{\e}, \sL_{X_{k\triangle}^{\e}})\Big\ra.
\de

In the following, we estimate $\Theta(s,r)$. For any $r>0$, $\zeta\in L^2(\Omega, \sF_r, \mP; \mR^m)$, $y\in\mR^m$, we introduce the auxiliary system:
\ce\left\{\begin{array}{l}
\check{Y}_s^{\e,r,\zeta}=\zeta+\frac1\e\int_r^sb_2(\check{Y}_u^{\e,r,\zeta}, \sL_{\check{Y}_u^{\e,r,\zeta}}) \dif u+ \frac1{\sqrt\e}\int_r^s\s_2(\check{Y}_u^{\e,r,\zeta}, \sL_{\check{Y}_u^{\e,r,\zeta}}) \dif W_u, \\
\check{Y}_s^{\e,r,y,\sL_{\zeta}}=y+\frac1\e\int_r^sb_2(\check{Y}_u^{\e,r,y,\sL_{\zeta}}, \sL_{\check{Y}_u^{\e,r,\zeta}}) \dif u+\frac1{\sqrt\e}\int_r^s\s_2(\check{Y}_u^{\e,r,y,\sL_{\zeta}},\sL_{\check{Y}_u^{\e,r,\zeta}})\dif W_u.
\end{array}\right.
\de
Then it holds that
\ce
Y_{s\e+k\triangle}^{\e,\xi}=\check{Y}_{s\e+k\triangle}^{\e,k\triangle,Y_{k\triangle}^{\e,\xi}},\quad
Y_{s\e+k\triangle}^{\e, y_0, \sL_{\xi}}=\check{Y}_{s\e+k\triangle}^{\e,k\triangle,Y_{k\triangle}^{\e, y_0, \sL_{\xi}},\sL_{Y_{k\triangle}^{\e,\xi}}},\quad s\in[0,\triangle/{\e}].
\de
Moreover, $ X_{k\triangle}^{\e}$, $Y_{k\triangle}^{\e, y_0, \sL_{\xi}}$ are $\sF_{k\triangle}$-measurable, and for any $y \in \mR^m$, $t>k\triangle$, $\check{Y}_t^{\e,k\triangle,y,\sL_{Y_{k\triangle}^{\e,\xi}}}$  is independent of $\sF_{k\triangle}$. Thus,
\be
&&\Theta(s,r)\no\\
&=&\mE\Big\la H(X_{k\triangle}^{\e}, \sL_{X_{k\triangle}^{\e}},\check{Y}_{s\e+k\triangle}^{\e,k\triangle,Y_{k\triangle}^{\e, y_0, \sL_{\xi}},\sL_{Y_{k\triangle}^{\e,\xi}}},\sL_{\check{Y}_{s\e+k\triangle}^{\e,k\triangle,Y_{k\triangle}^{\e,\xi}}})
-\bar{H}(X_{k\triangle}^{\e}, \sL_{X_{k\triangle}^{\e}}),\no\\
&&\qquad H(X_{k\triangle}^{\e}, \sL_{X_{k\triangle}^{\e}}, \check{Y}_{r\e+k\triangle}^{\e,k\triangle,Y_{k\triangle}^{\e, y_0, \sL_{\xi}},\sL_{Y_{k\triangle}^{\e,\xi}}},\sL_{\check{Y}_{r\e+k\triangle}^{\e,k\triangle,Y_{k\triangle}^{\e,\xi}}})
-\bar{H}(X_{k\triangle}^{\e}, \sL_{X_{k\triangle}^{\e}})\Big\ra\no\\
&=&\mE\Big[\mE\Big[\Big\la H(X_{k\triangle}^{\e}, \sL_{X_{k\triangle}^{\e}},\check{Y}_{s\e+k\triangle}^{\e,k\triangle,Y_{k\triangle}^{\e, y_0, \sL_{\xi}},\sL_{Y_{k\triangle}^{\e,\xi}}},\sL_{\check{Y}_{s\e+k\triangle}^{\e,k\triangle,Y_{k\triangle}^{\e,\xi}}})
-\bar{H}(X_{k\triangle}^{\e}, \sL_{X_{k\triangle}^{\e}}),\no\\
&&\qquad H(X_{k\triangle}^{\e}, \sL_{X_{k\triangle}^{\e}}, \check{Y}_{r\e+k\triangle}^{\e,k\triangle,Y_{k\triangle}^{\e, y_0, \sL_{\xi}},\sL_{Y_{k\triangle}^{\e,\xi}}},\sL_{\check{Y}_{r\e+k\triangle}^{\e,k\triangle,Y_{k\triangle}^{\e,\xi}}})
-\bar{H}(X_{k\triangle}^{\e}, \sL_{X_{k\triangle}^{\e}})\Big\ra\Big|\sF_{k\triangle}\Big]\Big]\no\\
&=&\mE\Bigg[\mE\Bigg[\Big\la H(x,\mu,\check{Y}_{s\e+k\triangle}^{\e,k\triangle,y,\nu},\sL_{\check{Y}_{s\e+k\triangle}^{\e,k\triangle,Y_{k\triangle}^{\e,\xi}}})
-\bar{H}(x,\mu),\no\\
&&\qquad H(x,\mu,\check{Y}_{r\e+k\triangle}^{\e,k\triangle,y,\nu},\sL_{\check{Y}_{r\e+k\triangle}^{\e,k\triangle,Y_{k\triangle}^{\e,\xi}}})
-\bar{H}(x,\mu)\Big\ra\Big|_{(x,\mu,y,\nu)=(X_{k\triangle}^{\e}, \sL_{X_{k\triangle}^{\e}},Y_{k\triangle}^{\e, y_0, \sL_{\xi}},\sL_{Y_{k\triangle}^{\e,\xi}})}\Bigg]\Bigg].\no\\
\label{HtribarHtril2}
\ee
Here, we investigate $\check{Y}_{s\e+k\triangle}^{\e,k\triangle,y,\nu}$.
\ce
\check{Y}_{s\e+k\triangle}^{\e,k\triangle,y,\nu}&=&y+\frac1\e\int_{k\triangle}^{s\e+k\triangle}b_2(\check{Y}_u^{\e,k\triangle,y,\nu}, \sL_{\check{Y}_u^{\e,k\triangle,Y_{k\triangle}^{\e,\xi}}}) \dif u+\frac1{\sqrt\e}\int_{k\triangle}^{s\e+k\triangle}\s_2(\check{Y}_u^{\e,k\triangle,y,\nu},\sL_{\check{Y}_u^{\e,k\triangle,Y_{k\triangle}^{\e,\xi}}})\dif W_u\\
&=&y+\frac1\e\int_{0}^{s\e}b_2(\check{Y}_{r+k\triangle}^{\e,k\triangle,y,\nu}, \sL_{\check{Y}_{r+k\triangle}^{\e,k\triangle,Y_{k\triangle}^{\e,\xi}}}) \dif r+\frac1{\sqrt\e}\int_{0}^{s\e}\s_2(\check{Y}_{r+k\triangle}^{\e,k\triangle,y,\nu},\sL_{\check{Y}_{r+k\triangle}^{\e,k\triangle,Y_{k\triangle}^{\e,\xi}}})\dif \check{W}_r\\
&=&y+\int_{0}^{s}b_2(\check{Y}_{v\e+k\triangle}^{\e,k\triangle,y,\nu}, \sL_{\check{Y}_{v\e+k\triangle}^{\e,k\triangle,Y_{k\triangle}^{\e,\xi}}}) \dif v+\int_{0}^{s}\s_2(\check{Y}_{v\e+k\triangle}^{\e,k\triangle,y,\nu},\sL_{\check{Y}_{v\e+k\triangle}^{\e,k\triangle,Y_{k\triangle}^{\e,\xi}}})\dif \tilde{\check{W}}_v,
\de
where $\check{W}_r:=W_{r+k\triangle}-W_{k\triangle}$, and $\tilde{\check{W}}_v:=\frac1{\sqrt\e}\check{W}_{v\e}$ are two $d_2$-dimensional standard Brownian motions. Then $(\check{Y}_{s\e+k\triangle}^{\e,k\triangle,y,\nu},\sL_{\check{Y}_{s\e+k\triangle}^{\e,k\triangle,Y_{k\triangle}^{\e,\xi}}})$ and $(Y_s^{y,\nu},\sL_{Y_s^\zeta})$ have the same distribution, where $\nu=\sL_\zeta$. 

Next, we investigate the property of $(Y_s^{y,\nu},\sL_{Y_s^\zeta})$. Let $C\left([0, \triangle/\e], \mR^m\right)$ be the collection of continuous functions from $[0, \triangle/\e]$ to $\mR^m$, endowed with the uniform convergence topology. Set
\ce
&&\tilde{\Omega}:=C\left([0, \triangle/\e], \mR^m\right) \times C\left([0, \triangle/\e],
\cP_2\left(\mR^m\right)\right), \\
&&\tilde{\sF}:=\sB\Big(C\left([0, \triangle/\e], \mR^m\right)\Big) \times \sB\Big(C\left([0, \triangle/\e], \cP_2\left(\mR^m\right)\right)\Big),\\
&&\tilde{\sF}_s:=\s\left(M_r: 0 \leq r \leq s\right), \quad 0 \leq s \leq \triangle/\e,
\de
where $M_\cdot$  denotes the coordinate process.
Following a similar reasoning to that in \cite[Theorem 4.11]{{rrw}}, there exists a unique probability measure $\tilde{\mP}$ on $(\tilde{\Omega},\tilde{\sF})$ such that $M_\cdot$ is a homogeneous Markov process with respect to $(\tilde{\sF}_s)$, and its transition function is given by $\{P_s(y,\nu;\cdot)=\sL_{Y_s^{y,\nu}}\times\d_{\sL_{Y_s^\zeta}}: 0\leq s\leq\triangle/\e, (y,\nu)\in\mR^m\times\cP_2\left(\mR^m\right)\}$ and $\sL_{M_0}^{\tilde{\mP}}=\d_y\times\d_\nu$. Note that
$$
\sL_{M_s}^{\tilde{\mP}}:=\tilde{\mP}\circ M_s^{-1}=\int_{\mR^m\times\cP_2(\mR^m)}P_s(y^\prime,\nu^\prime;\cdot)(\d_y\times\d_\nu)(\dif y^\prime,\dif\nu^\prime)=\sL_{Y_s^{y,\nu}}\times \d_{\sL_{Y_s^\zeta}}.
$$
Set
$$
(\tilde{Y}_s^{y,\nu},\sL^{\tilde{\mP}}_{\tilde{Y}_s^{\tilde{\zeta}}}):=M_s.
$$
Therefore, it follows that for $0< r< s\leq\triangle/\e$,
\be
&&\mE\Big\la H(x,\mu,\check{Y}_{s\e+k\triangle}^{\e,k\triangle,y,\nu},\sL_{\check{Y}_{s\e+k\triangle}^{\e,k\triangle,Y_{k\triangle}^{\e,\xi}}})
-\bar{H}(x,\mu), H(x,\mu,\check{Y}_{r\e+k\triangle}^{\e,k\triangle,y,\nu},\sL_{\check{Y}_{r\e+k\triangle}^{\e,k\triangle,Y_{k\triangle}^{\e,\xi}}})
-\bar{H}(x,\mu)\Big\ra\no\\
&=&\tilde{\mE}\Big\la H(x,\mu,\tilde{Y}_s^{y,\nu},\sL^{\tilde{\mP}}_{\tilde{Y}_s^{\tilde{\zeta}}})-\bar{H}(x,\mu), H(x,\mu,\tilde{Y}_r^{y,\nu},\sL^{\tilde{\mP}}_{\tilde{Y}_r^{\tilde{\zeta}}})-\bar{H}(x,\mu)\Big\ra\no\\
&=&\tilde{\mE}\Big[\tilde{\mE}\Big[\Big\la H(x,\mu,\tilde{Y}_s^{y,\nu},\sL^{\tilde{\mP}}_{\tilde{Y}_s^{\tilde{\zeta}}})-\bar{H}(x,\mu), H(x,\mu,\tilde{Y}_r^{y,\nu},\sL^{\tilde{\mP}}_{\tilde{Y}_r^{\tilde{\zeta}}})-\bar{H}(x,\mu)\Big\ra\Big|
\tilde{\sF}_r\Big]\Big]\no\\
&=&\tilde{\mE}\Big[\Big\la\tilde{\mE}\left[H(x,\mu,\tilde{Y}_s^{y,\nu},\sL^{\tilde{\mP}}_{\tilde{Y}_s^{\tilde{\zeta}}})-\bar{H}(x,\mu)
\Big|\tilde{\sF}_r\right], H(x,\mu,\tilde{Y}_r^{y,\nu},\sL^{\tilde{\mP}}_{\tilde{Y}_r^{\tilde{\zeta}}})-\bar{H}(x,\mu)\Big\ra\Big]\no\\
&=&\tilde{\mE}\Big[\Big\la\tilde{\mE}\Big[H(x,\mu,\tilde{Y}_{s-r}^{\hat{y},\hat{\nu}},\sL^{\tilde{\mP}}_{\tilde{Y}_{s-r}^{\hat{\zeta}}})\Big]
\Big|_{(\hat{y},\hat{\nu},\hat{\zeta})=(\tilde{Y}_r^{y,\nu},\sL^{\tilde{\mP}}_{\tilde{Y}_r^{\tilde{\zeta}}},\tilde{Y}_r^{\tilde{\zeta}})}
-\bar{H}(x,\mu), H(x,\mu,\tilde{Y}_r^{y,\nu},\sL^{\tilde{\mP}}_{\tilde{Y}_r^{\tilde{\zeta}}})-\bar{H}(x,\mu)\Big\ra\Big]\no\\
&\leq&\left(\tilde{\mE}\Big\|\tilde{\mE}\Big[H(x,\mu,\tilde{Y}_{s-r}^{\hat{y},\hat{\nu}},\sL^{\tilde{\mP}}_{\tilde{Y}_{s-r}^{\hat{\zeta}}})\Big]
\Big|_{(\hat{y},\hat{\nu},\hat{\zeta})=(\tilde{Y}_r^{y,\nu},\sL^{\tilde{\mP}}_{\tilde{Y}_r^{\tilde{\zeta}}},\tilde{Y}_r^{\tilde{\zeta}})}
-\bar{H}(x,\mu)\Big\|^2\right)^{\frac12}\no\\
&&\qquad\qquad\qquad\qquad\qquad\cdot\left(\tilde{\mE}\Big\|H(x,\mu,\tilde{Y}_r^{y,\nu},\sL^{\tilde{\mP}}_{\tilde{Y}_r^{\tilde{\zeta}}})-\bar{H}(x,\mu)\Big\|^2\right)^{\frac12},\label{ceas}
\ee
where $\tilde{\mE}$ is the expectation on $(\tilde{\Omega},\tilde{\sF},\{\tilde{\sF_t}\}_{t\geq0},\tilde{\mP})$. Due to \eqref{HbarH} and Lemma \ref{ymul2}, we obtain that
\be
&&\left(\tilde{\mE}\Big\|\tilde{\mE}\Big[H(x,\mu,\tilde{Y}_{s-r}^{\hat{y},\hat{\nu}},\sL^{\tilde{\mP}}_{\tilde{Y}_{s-r}^{\hat{\zeta}}})\Big]
\Big|_{(\hat{y},\hat{\nu},\hat{\zeta})=(\tilde{Y}_r^{y,\nu},\sL^{\tilde{\mP}}_{\tilde{Y}_r^{\tilde{\zeta}}},\tilde{Y}_r^{\tilde{\zeta}})}
-\bar{H}(x,\mu)\Big\|^2\right)^{\frac12}\no\\
&\leq&Ce^{-\frac{\b_1-\b_2}{2}(s-r)}\left(1+\tilde{\mE}|\tilde{Y}_r^{y,\nu}|^{8}+\tilde{\mE}|\tilde{Y}_r^{\tilde{\zeta}}|^{8}\right)^{\frac12}\no\\
&\leq&Ce^{-\frac{\b_1-\b_2}{2}(s-r)}\left(1+|y|^{12}+\mE|\zeta|^{12}\right)^{\frac12}.\label{mEH}
\ee
By \eqref{Hgrow}, \eqref{barHboun} and Lemma \ref{ymul2}, we have
\ce
&&\left(\tilde{\mE}\Big\|H(x,\mu,\tilde{Y}_r^{y,\nu},\sL^{\tilde{\mP}}_{\tilde{Y}_r^{\tilde{\zeta}}})-\bar{H}(x,\mu)\Big\|^2
\right)^{\frac12}\\
&\leq&2\left(\tilde{\mE}\Big\|H(x,\mu,\tilde{Y}_r^{y,\nu},\sL^{\tilde{\mP}}_{\tilde{Y}_r^{\tilde{\zeta}}})\Big\|^2
+\left\|\bar{H}(x,\mu)\right\|^2\right)^{\frac12}\\
&\leq&C\left(1+\tilde{\mE}|\tilde{Y}_r^{y,\nu}|^4+\tilde{\mE}|\tilde{Y}_r^{\tilde{\zeta}}|^4\right)^{\frac12}\\
&\leq&C\left(1+|y|^{12}+\mE|\zeta|^{12}\right)^{\frac12},
\de
which together with \eqref{HtribarHtril2}-\eqref{mEH} and Lemma \ref{xtyt} yields that
\ce
\Theta(s,r)
&\leq&Ce^{-\frac{\b_1-\b_2}{2}(s-r)}\left(1+\mE|Y_{k\triangle}^{\e, y_0, \sL_{\xi}}|^{12}+\mE|Y_{k\triangle}^{\e, \xi}|^{12}\right)\\
&\leq&Ce^{-\frac{\b_1-\b_2}{2}(s-r)}\left(1+|y_0|^{12}+\mE|\xi|^{12}\right).
\de
Inserting the above inequality in \eqref{su1esti}, we get
\ce
\mathscr{U}_1(t)
&\leq&\frac{C_T\e^2}{\triangle^2}\int_{0}^{\frac\triangle\e}\int_{r}^{\frac\triangle\e}e^{-\frac{\b_1-\b_2}{2}(s-r)}\left(1+|y_0|^{12}+\mE|\xi|^{12}\right)
\dif s\dif r\\
&\leq&C_T(1+\vert y_0\vert^{12}+\mE\vert\xi\vert^{12})\frac{\e}{\triangle}.
\de
Therefore, the inequality \eqref{su1} holds.


\begin{thebibliography}{999}
\bibitem{brw}
J. Bao, P. Ren and F-Y. Wang: Bismut formula for Lions derivative of distribution-path dependent SDEs,
{\it J. Differential Equations}, 282(2021)285-329.


\bibitem{bhr}
Z. Brze\'{z}niak, E. Hausenblas and  P.A. Razafimandimby: Stochastic reaction-diffusion equations driven by jump processes, {\it Potential Anal}, 49(2018)131-201.

\bibitem{blpr}
R. Buckdahn, J. Li, S. Peng and C. Rainer: Mean-field stochastic differential equations
and associated PDEs, {\it Ann. Probab.}, 45(2017)824-878.

\bibitem{cm} D. Crisan and E. McMurray: Smoothing properties of McKean-Vlasov SDEs, {\it Probab. Theory Related Fields}, 171(2018)97-148.

\bibitem{pz}
G. Da Prato and J. Zabczyk: Stochastic equations in infinite dimensions, Cambridge University Press, Cambridge (1992).

\bibitem{dq1}
X. Ding and H. Qiao: Euler-Maruyama approximations for stochastic McKean-Vlasov equations
with non-Lipschitz coefficients, {\it J. Theoret. Probab.}, 34(2021)1408-1425.

\bibitem{dq2}
X. Ding and H. Qiao: Stability for stochastic McKean-Vlasov equations with non-Lipschitz coefficients, {\it SIAM J. Control Optim.}, 59(2021)887-905.


\bibitem{dsxz}
Z. Dong, X. Sun, H. Xiao and J. Zhai: Averaging principle for one dimensional stochastic Burgers equation, {\it J. Differential Equations}, 265(2018)4749-4797.

\bibitem{fyy} X. Fan, T. Yu and C. Yuan: Asymptotic behaviors for distribution dependent SDEs driven by fractional Brownian motions, {\it Stochastic Processes and their Applications}, 164(2023)383-415.

\bibitem{flqz} K. Fang, W. Liu, H. Qiao and F. Zhu: Asymptotic behaviors of small perturbation for multivalued McKean-Vlasov stochastic differential equations, {\it Applied Mathematics and Optimization}, 88(2023)22.

\bibitem{hlls}
W. Hong, S. Li, W. Liu and X. Sun: Central limit type theorem and large deviation principle for multi-scale McKean-Vlasov SDEs, {\it Probab. Theory Related Fields}, 187(2023)133-201.

\bibitem{hll}
W. Hong, S. Li and W. Liu: Strong convergence rates in averaging principle for slow-fast McKean-Vlasov
SPDEs, {\it J. Differential Equations}, 316(2022)94-135.

\bibitem{hw}
X. Huang and F.-Y. Wang: Distribution dependent SDEs with singular coefficients, {\it Stochastic Process. Appl.}, 129(2019)4747-4770.

\bibitem{hrw}
X. Huang, P. Ren and F.-Y. Wang: Distribution dependent stochastic differential equations,
{\it Front. Math. China}, 16(2021)257-301.

\bibitem{Huang}
Z. Huang: {\it Basis of stochastic analysis (in Chinese)}, 2nd ed., Science Press, Beijing (2001).

\bibitem{rk}
R. Khasminskii: On the principle of averaging the It\^{o} stochastic differential equations, {\it Kibernetika}, 4(1968)260-279.

\bibitem{kn}
C. Kumar and Neelima: On explicit Milstein-type scheme for McKean-Vlasov stochastic differential equations with super-linear drift coefficient, {\it Electron. J. Probab.}, 26(2021)1-32.

\bibitem{lx}
Y. Li and L. Xie: Functional law of large numbers and central limit theorem for slow-fast McKean-Vlasov equations, {\it Discrete Contin. Dyn. Syst. Ser. S}, 16(2023)846-877.

\bibitem{lrsx}
W. Liu, M. R\"{o}ckner, X. Sun and Y. Xie: Averaging principle for slow-fast stochastic differential equations with time dependent locally Lipschitz coefficients, {\it J. Differential Equations}, 268(2020)2910-2948.

\bibitem{hpm}
H.P. McKean: A class of Markov processes associated with nonlinear parabolic equations, {\it Proc. Natl. Sci. U.S.A.}, 56(1966)1907-1911.

\bibitem{pix}
B. Pei, Y. Inahama and Y. Xu: Averaging principle for fast-slow system driven by mixed fractional Brownian rough path, {\it J. Differential Equations}, 301(2021)202-235.

\bibitem{qw}
 H. Qiao and W. Wei: Strong approximation of nonlinear filtering for multiscale McKean-Vlasov stochastic systems, http://arxiv.org/abs/2206.05037.

\bibitem{rrw}
P. Ren, M. R\"{o}ckner and F.-Y. Wang: Linearization of nonlinear Fokker-Planck equations and applications, {\it J. Differential Equations}, 322(2022)1-37.

\bibitem{rw}
P. Ren and F.-Y. Wang: Space-Distribution PDEs for path independent additive functionals of McKean-Vlasov SDEs, {\it Infin. Dimens. Anal. Quantum Probab. Relat. Top.}, 23(2020)2050018.

\bibitem{rsx}
 M. R\"{o}ckner, X. Sun and Y. Xie: Strong convergence order for slow-fast McKean-Vlasov stochastic differential equations, {\it Ann. Inst. H. Poincar\'{e} Probab. Statist.}, 57(2021)547-576.

\bibitem{rz}
M. R\"{o}ckner and X. Zhang: Well-posedness of distribution dependent SDEs with singular drifts, {\it Bernoulli}, 27(2021)1131-1158.

\bibitem{Shen1}
G.  Shen, J. Xiang and J.-L. Wu: Averaging principle for distribution dependent stochastic differential equations driven by fractional Brownian motion and standard Brownian motion, {\it J. Differential Equations}, 321(2022)381-414.

\bibitem{Shen2}
G. Shen, J. Yin and J.-L. Wu: Stochastic averaging principle for two-time-scale SDEs with distribution dependent coefficients driven by fractional Brownian motion, {\it Comm. Math. Statist.}, 2023, https://doi.org/10.1007/s40304-023-00364-4.

\bibitem{sy}
Y. Suo and C. Yuan: Central limit theorem and moderate deviation principle for McKean-Vlasov SDEs, {\it Acta Appl. Math.}, 175(2021)1-19.


\bibitem{Wang}
F.-Y. Wang: Distribution dependent SDEs for Landau type equations, {\it Stochastic Process. Appl.}, 128(2018)595-621.

\bibitem{wy}
F. Wu and G. Yin: Fast-slow-coupled stochastic functional differential equations, {\it J. Differential Equations}, 323(2022)1-37.

\bibitem{xllm}
J. Xu, J. Liu, J. Liu and Y. Miao: Strong averaging principle for two-time-scale stochastic McKean-Vlasov equations, {\it Appl. Math. Optim.}, 84(2021)837-867.
\end{thebibliography}
\end{document}